\documentclass[11pt]{amsart}
                             
                        \usepackage{amsmath}
                             \theoremstyle{remark}
                             \newcommand{\re}[1]{(\ref{#1})}

                             \selectfont
                             \newtheorem{thm}{\textbf{Theorem}}[section]
                             \newtheorem{lem}{\textbf{Lemma}}[section]

                             \theoremstyle{plain}

                             \newtheorem{defn}{\textbf{Definition}}[section]

                             \numberwithin{equation}{section}
                             \numberwithin{figure}{section}
                             \newenvironment{prev}[1][\underline{Proof}]{\noindent\textbf{#1.} }{\ \rule{0.6em}{0.6em}}

\begin{document}
                     \title[Optimal control stochastic AC-NS]{ Optimal boundary control for the stochastic Allen-Cahn Navier-Stokes system in two dimensions}
                     \author{R. D. AYISSI$^{*}$, G. Deugou\'{e}$^{**}$, J. Ngandjou Zangue$^{*,**}$  and  T. Tachim Medjo$^{***}$}
                     \dedicatory{\vspace{-10pt}\normalsize{
                     $^{*}$ Department of Mathematics, Faculty of Science, University of Yaoundé, P. O. BOX 812, Yaoundé, Cameroon \\
                     $^{**}$  Department of Mathematics and Computer Science, University of Dschang, P. O. BOX 67, Dschang, Cameroon \\
                      $^{***}$ Department of Mathematics and Statistics, Florida International University, MMC, Miami, FL $33199$, USA }}
                             
                      \begin{abstract}
                      In this work, we study an optimal boundary control for the stochastic Allen–Cahn Navier–Stokes system. The governing system of nonlinear partial differential equations consists of the stochastic Navier-Stokes equations with non homogeneous Navier-slip boundary condition coupled with a phase-field equation, which is the convective Allen-Cahn equation type. We investigate the well-posedness of the nonlinear system. More precisely, the existence and uniqueness of global strong  solution in dimension two is established and  we prove the existence of an optimal solution to the control problem.                    
                       \end{abstract}
                     \keywords{ Optimal stochastic control; Allen-Hilliard Navier-Stokes system; Diffuse interface model; Non-homogeneous boundary conditions.}
                     \date{\today}
                           \maketitle
                      \section{Introduction}
                       A well-known system which models the flows of an incompressible binary mixture of immiscible fluids subject to phase separation is the so-called model "\textbf{H}" (see \cite{A25}, see also \cite{A30}). Modeling the bahavior of immiscible binary fluids is a very active area of research because of its importance, for instance, in Biology and materials sciences. This paper is devoted to the study of a control problem for the stochastic Allen-Cahn Navier-Stokes (AC-NS) system with non-homogeneous Navier--slip boundary conditions in a bounded domain  $\text{D}\subset\mathbb{R}^{2}$. The system is considered as a model describing the motion of binary fluid mixtures. 
                       
                       Let $(\Omega, \mathcal{F},\mathbb{P},\{\mathcal{F}_{t}\}_{t\geq 0})$ be a filtered probability space satisfying the usual hypotheses i.e., $\{\mathcal{F}_{t}\}_{t\geq 0}$ is a right-continuous filtration such that $\mathcal{F}_{0}$ contains all the $\mathbb{P}$-null subsets of $(\Omega,\mathcal{F})$. In this study, we consider the stochastic AC-NS system with multiplicative noise given by: 
                            \begin{equation}\label{NS1}
                            \begin{cases}
                          	d\text{v} +(-\nu\Delta \text{v}+(\text{v}\cdot\nabla)\text{v}+\nabla \text{p}-\mathcal{K}\mu\nabla\phi)dt = \text{g}(t,\text{v})dW,\quad \mbox{in}\ \   \text{D}\times(0,T),\\
                           	\text{div v}=0,\ \ \mbox{in}\ \   \text{D}\times(0,T),\\
                             \frac{\partial \phi}{\partial t}+\text{v}\cdot\nabla\phi +\mu= 0,\quad \mbox{in}\ \   \text{D}\times(0,T),\\
                               \mu= -\delta\Delta\phi +\xi\text{f}(\phi),\ \ \mbox{in}\ \   \text{D}\times(0,T),\\
                            \frac{\partial \phi}{\partial n} =0, \ \ \text{on} \ \ \partial \text{D}\times(0,T),\\
                             \text{v}\cdot n=a,\ \ [2D(\text{v}) n+\alpha \text{v}]\cdot\tau=b,\ \ \text{on}\quad \Gamma_{T}=\partial \text{D}\times(0,T),\\
                            (\text{v},\phi)(0)=(\text{v}_{0},\phi_{0}), \quad \text{in} \ \ \text{D}.
                             \end{cases}
                             \end{equation}
           In $\eqref{NS1}$, $T>0$ fixed, $\text{v}=(\text{v}_{1},\text{v}_{2})$ is the unknown functions are the velocity of the fluid, its pressure $\text{p}$ and the order (phase) parameter $\phi$. The quantity $\mu$ is the chemical potential of the binary mixture which is given as the variational derivative of the following free energy functional
          \begin{equation*}
            \mathcal{F}(\phi)=\int_{\text{D}}(\frac{\delta}{2}\left|\nabla\phi\right|^{2}+\xi\text{F}(\phi))dx,
            \end{equation*}
     where $\text{F}(r)=\int^{r}_{0}f(\zeta)d\zeta$ is the suitable double-well potential of the regular type. A typical example of potential  of regular F is 
      $\text{F}(r)=(r^{2}-1)^{2}$,\ \ $r\in\mathbb{R}$.  The quantities $\nu$ and $\mathcal{K}$ are positive constants that correspond to the kinematic viscosity of the fluid and capillarity (stress) coefficient, respectively. Here, $\delta$ and $\xi$ are two positive parameters describing the interactions between the two phases. In particular, $\delta$ is related to the thickness of the interface separating the two fluids. $(\text{v}_{0},\phi_{0})$ is the initial condition such that $\text{v}_{0}$ verifies
                                     \begin{equation}
                                     \text{div v}_{0}=0 \ \ \mbox{in}\ \   \text{D}.
                                     \end{equation}
     Here
                        \begin{math}
                        D(\text{v})=\frac{1}{2}[\nabla\text{v}+(\nabla\text{v})^{T}]
                        \end{math}
       is the rate-of-strain tensor; $n$ is the external unit normal to the boundary $\partial \text{D}\in\mathcal{C}^{2}$ of the domain $\text{D}$ and $\tau$ is the tangent unit vector to $\partial \text{D}$, such that $(n,\tau)$ forms a standard orientation in $\mathbb{R}^{2}$. The positive constant $\alpha$ is the so-called friction coefficient. The quantity $a$ corresponds to the inflow and outflow fluids trough $\partial \text{D}$, satisfying the compatibility condition
                        \begin{equation}
                        \int_{\partial \text{D}}a(t,x)d\gamma=0, \ \ \mbox{for any}\ \ t\in[0,T],
                        \end{equation}
    which is a consequence of the fact that the divergence-free vector $\text{v}$ satisfies the boundary condition $\text{v}\cdot n=a$ on $\Gamma_{T}.$
   This condition means that the quality of inflow fluid should coincide with the quantity of outflow fluid. The boundary function $a$ and $b$ will be considered as the control variables for the physical system $\eqref{NS1}$. The term $\text{g}(t,\text{v})dW$ is a multiplicative noise depending eventually on $\text{v}$. This work is the first contribution to the analytic study of the boundary control problem for the stochastic AC-NS system with non-homogeneous Navier-slip boundary conditions. 
   
   In the literature, the mathematical study of stochastic two-phase flows models with the Dirichlet boundary condition were studied by several authors (see, for instance \cite{Ta1, TD5, jid, Ta33, GTC, TW, Zhang, TAL, GTLD}) for the stochastic CH-NS system and (see, for instance \cite{A251, A20, A21} and the references therein) for the stochastic AC-NS system. Let us mention that the main difference between the Cahn–Hilliard equations and the Allen–Cahn equations is the fourth-order dissipative operator that appears in the Cahn–Hilliard equations, while the Allen–Cahn equations enjoy a dissipative operator of order only two. Because of this important difference many results that are available for the CH-NS model may not be easily replicated for the AC-NS system, even in the deterministic case. More precisely, the phase parameter $\phi$ enjoys more regularity (in space) in the CH-NS model than in the AC-NS system. This weak regularity in space makes the mathematical analysis of the capillarity term $\mu\nabla\phi$ (which is the main nonlinear term) in the AC-NS model more difficult than in the CH-NS system. We should note that even if the Dirichlet boundary condition is widely accepted as an appropriate boundary condition at the surface of the contact between a fluid and a solid, it is important to recall that there is another physically relevant boundary conditions is the so-called Navier-slip boundary condition, which reads as 
   \begin{equation}\label{mn}
   \text{v}\cdot n=0, \ \ [2D(\text{v})n+\alpha\text{v}]=0.
   \end{equation}
   Let us recall that in 1823, Navier's proposed a slip-with-friction boundary condition (see $\eqref{mn}$) and claimed that the component of the fluid velocity tangent to the surface should be proportional to the rate of strain at the surface \cite{NA}. This condition was rigorously justified in \cite{AD3} on the bases of miscroscopic boundary condition for the Boltzmann equations. Since this justification, several papers have been published (see, for instance \cite{AD1, AD2, S.H} and the references therein) for the deterministic Navier-Stokes equations and we refer to \cite{KIS} for the stochastic second grade fluid equations.
   
    Motivated by the study of optimal control problems for single fluid flows, we are interested to study the optimal control problem for the binary mixture described by the stochastic AC-NS system with Navier-slip boundary conditions. Optimal control problems associated with the stochastic AC-NS system have many applications where one is interested in influencing the phase separation process, such as in the separation of binary alloys, the formation of polymeric membranes, etc. The boundary control of the fluid flows is of main importance in several branches of the industry, for instance, in the aviation industry. Extensive research has been carried out concerning the implementation of injection-suction device to control the motion of the fluids (see, for instance \cite{S.H, TLB, ALB}). 
   
   Optimal control problems associated with the deterministic two-phase flow models with Dirichlet boundary condition have been widely studied in recent years, see for instance \cite{TDM1, TDM2, YB, YB0} and the references therein. In the stochastic case, not many results are available for two-phase flow models.
   Optimal control problems for the stochastic Navier-Stokes equations have been studied by several mathematicians and some of the results can be found for instance in  \cite{H.B1, H.B3, PB1, PB2, PB3, Moh}. In \cite{H.B3}, the author proved the existence of optimal and $\varepsilon-$optimal controls for the stochastic Navier-Stokes equation in a two dimensional bounded domain, controlled by different external forces  and assuming that the set of admissible controls is weakly sequentially compact. In \cite{Moh}, the authors studied the existence of optimal control for the stochastic Navier-Stokes equation in a two dimensional bounded domain. The proof is based on the dynamic programming and some properties of the transition semigroup corresponding to the stochastic 2D Navier-Stokes equation.  In \cite{PB1}, using an approach  based on semi-group  theory, the authors proved the well-posedeness of local mild solution to the stochastic Navier-Stokes equations with multiplicative Lévy noise. Then, they proved the existence and uniqueness of optimal control, by analyzing some cost functional depending on stopping time. They studied in \cite{PB2}, the optimality system and derived the necessary and sufficient optimality conditions with a multiplicative noise driven by a Wiener process. For optimal control problem with the stochastic Cahn-Hilliard equations, we can refer to \cite{MA10, A250}. Let us note that, in \cite{MA10}, a distributed optimal control problem associated to the stochastic pure Cahn-Hilliard equation was investigated. The control is represented by  a source-term in the equation of the chemical potential. The existence of an optimal control is then established thanks to some probability  and analytical compactness arguments. Moreover, the first order necessary conditions of optimality was also established by means of the monotonicity and compactness arguments. To the best of our knowledge,  There are not
   work written on problem control for the stochastic Allen-Cahn equations.  In \cite{josue}, the authors proved the existence of optimal and $\epsilon$-optimal controls for the stochastic CH-NS equation in a two dimensional bounded domain, controlled by different external forces  and assuming that the set of admissible controls is weakly sequentially compact.      
                     
    Our main objectives in this paper are: first to prove the existence and uniqueness of a strong solution (in the stochastic sense) for system $\eqref{NS1}$ by using a Galerkin approximation, as employed in previous studies. We deduce uniform estimates for the approximate solution that allows us to pass to the limit with respect to weak topology. In order to shown that the limit process is a solution, we employ the methods developed in \cite{Ta1}, then to control the solution of the system $\eqref{NS1}$ by the boundary condition $(a,b)$, which is a predictable stochastic process belonging to the space $\mathcal{A}$ of admissible controls defined in Section $\ref{N03}$. The cost functional is given by
                  \begin{equation}\label{Ns2}
                 \begin{aligned}
                 \mathcal{J}(a,b,\text{v},\phi)=\mathbb{E}\int_{0}^{T}\int_{\text{D}}\frac{1}{2}\left( \left|\text{v}-\text{v}_{d} \right|^{2}+\left|\phi-\phi_{d} \right|^{2}\right)dxdt
                 +\mathbb{E}\int_{\Gamma_{T}}\frac{1}{2}\left(\lambda_{1}\left| a\right|^{2}+\lambda_{2}\left| b\right|^{2} \right) d\gamma dt.
                 \end{aligned}
                 \end{equation}
                 The optimal control problem is formulated as follows:
                 \begin{equation*}
                 \textbf{(OCP)} \ \ \min_{(a,b)\in \mathcal{U}}\mathcal{J}(a,b,\text{v},\phi),
                 \end{equation*}
      subject to $\eqref{NS1}$. Here $\text{v}_{d}$ and $\phi_{d}$ are the target functions such that $\text{v}_{d}, \phi_{d}\in L^{2}(\Omega\times \text{D}\times(0,T))$. So motivated of the problem \textbf{(OCP)} is that we want to find the best control $(a,b)$ from the set of admissible controls such that the corresponding optimal solution of $\eqref{NS1}$ is as close as possible to the target state.  
                 
    As it is necessary to know the existence and uniqueness of strong solution to discuss the optimality principle, we restrict ourselves to the case of dimension two since the regularity of strong solution is important to study the optimal control problems. 
     
   The remainder of this article is organized as follows: in section, we present the general setting by introducing the appropriate functional spaces and some necessary classical inequalities. The existence and uniqueness of a strong  stochastic solution to \re{NS1} is formulated in section 3. In section 4, we establish the well-posedness of the model. In section 5, we prove the existence of an optimal solution to the  control problem.
     
                   \section{General setting}     
    Before we state the  main results of this paper, we fix the notation used through the text.\newline
   Let $X$ real Hilbert space with inner product $(.,.)_{X}$, we will denote the induced norm by  $\left\| .\right\|_{X} $, while $X^{*}$ will indicate its dual. Let $\text{D}\subset\mathbb{R}^{2}$ be a boundary domain with smooth boundary $\partial\text{D}$ and let $\left|\text{D} \right| $ denote the measure of the set $\text{D}$. By $\mathbb{L}^{p}(\text{D})=(L^{p}(\text{D}))^{2}$, with $0\leq p\leq \infty$, we denote the standard Lebesgue spaces and $\mathbb{H}^{p}(\text{D})=(H^{P}(\text{D}))^{2}$, the Sobolev spaces. The norms on $\mathbb{L}^{p}(\text{D})$  and $\mathbb{H}^{p}(\text{D})$ are indicated by $\left\|. \right\|_{\mathbb{L}^{p}} $ and $\left\|. \right\|_{\mathbb{H}^{p}} $. Let us  introduce the following Hilbert spaces:
           \begin{equation*}
               \begin{aligned}
               \mathbb{H}_{div}&=\{ u\in\mathbb{L}^{2}(\text{D}): div\, u=0  \ \ \mbox{in} \ \ \mathcal{D}'(\text{D}),  \ \ u.n=0  \ \ \mbox{in} \ \ \mathbb{H}^{-1/2}(\Gamma)\},\\
                \mathbb{V}_{div}&=\{ u\in\mathbb{H}^{1}(\text{D}): div\, u=0 \ \ \mbox{a.e.} \ \ \text{D}, \ \ u.n=0 \ \ \mbox{in}  \ \ \mathbb{H}^{1/2}(\Gamma) \}.
                 \end{aligned}
                  \end{equation*}
      We denote $\left|. \right|_{L^{2}} $  and $(.,.)$ the norm and the scalar product, respectively on $\mathbb{H}_{div}$.  We now define
                           \begin{equation*}
                           (Dv,Dz)=\int_{D}Dv\cdot Dz dx.
                             \end{equation*}
    On the space $\mathbb{V}_{div}$, we define the following inner product
                              \begin{equation*}
                              ((v,z))=2(Dv,Dz)+ \alpha\int_{\Gamma}v.z,
                             \end{equation*}
    and the corresponding norm $\left\|v \right\|=\left|\nabla v \right|_{L^{2}}=((v,v))^{1/2} $. Let us notice that, by Poincaré's and Korn's inequalities $\left|\nabla v \right|_{L^{2}}$ and $\left|D v \right|_{L^{2}}$ are equivalent norms in $\mathbb{V}_{div}$.  Since $\text{D}$ is bounded, the embedding of $\mathbb{V}_{div}\subset \mathbb{H}_{div}\equiv \mathbb{H}^{*}_{div}\subset\mathbb{V}^{*}_{div}$ is compact.\\
      Hereafter, we assume that $f\in\mathcal{C}^{1}(\mathbb{R})$ satisfies
                               \begin{equation}\label{f}
                              \begin{cases}
                              \lim\limits_{|r|\rightarrow\infty}f'(r)>0,\\
                                |f^{(i)}(r)|\leq c_{f}(1+|r|^{4-i}),\quad \forall\, r\in\mathbb{R}, \ \ i=0,1,
                               \end{cases}
                              \end{equation}
     where $c_{f}$ is some positive constant. 
      Note that the derivation of the typical double-well potential $f$ satisfies conditions similar to $(\ref{f})$. It follows from $(\ref{f})$ that, we can find $\vartheta>0$ such that: 
                    \begin{equation}\label{kl}
                    \lim\limits_{|r|\rightarrow\infty}f'(r)>2\vartheta>0.
                    \end{equation}
     Now, we introduce the linear positive unbounded operator $\text{A}_{\vartheta}$ on $\text{L}^{2}(\text{D})$ by: 
                               \begin{equation}\label{N7}
                           \text{A}_{\vartheta}\phi=-\Delta\phi+\vartheta\phi, \quad \forall\phi\in \text{D}(\text{A}_{\vartheta})=\{\phi\in \text{H}^{2}(	\text{D}), \partial_{n}\phi=0, \quad \text{on}\ \ \partial D\}.
                         \end{equation}
   Note that $\text{A}_{\vartheta}$ is self-adjoint, positive definite, unbounded linear operator on $\text{L}^{2}(\text{D})$ with compact inverse. More generally, we can define $\text{A}^{s}_{\vartheta}$, for any $s\in\mathbb{R}$, noting that $|\text{A}^{s/2}_{\vartheta}.|_{L^{2}}$, $s>0$, is an equivalent norm to the canonical $\text{H}^{s}$-norm on $\text{D}(\text{A}^{s/2}_{\vartheta})\subset\text{H}^{s}(D)$. Moreover, we set $\text{H}^{-s}(\text{D})=(\text{H}^{s}(\text{D}))^{*}$, wherever $s<0$.
    We also set 
    \begin{equation}
    f_{\vartheta}(r)=f(r)-\xi^{-1}\delta r,
    \end{equation}
 and observe that $f_{\vartheta}$   still satisfies $(\ref{f})$ with $\vartheta$   in place of $2\vartheta$ since $\delta\leq \xi$. Also, its primitive   $\text{F}_{\vartheta}(r)=\int^{r}_{0}f_{\vartheta}(\zeta)d\zeta$ is bounded from below.\\             
    Hereafter, we set 
        \begin{equation}\label{N9}
                 \text{H}=L^{2}(\text{D})=L^{2}(\text{D},\mathbb{R}), \ \  \text{V}_{s}=H^{s}(\text{D})=H^{s}(\text{D},\mathbb{R})=\text{D}(\text{A}^{s/2}_{\vartheta}).
                  \end{equation}
      The norms in $\text{H}$ and $\text{V}_{s}$ are denoted respectively by $|.|_{L^{2}}$ and $||.||_{s},$ where $||\psi||_{s}=|\text{A}^{s/2}_{\vartheta}\psi|_{L^{2}}$.\newline
      We denote $L^{p}(0,T; X)$ as the space of $X-$ valued measurable $p-$integrable functions defined on $[0,T]$ for $p\geq 1$.
       For $p,r\geq 1$, let $L^{p}(\Omega, L^{r}(0,T; X))$ be the space of the processes $v=v(\omega,t)$ with valued in $X$ defined on $\Omega\times [0,T]$, adapted to the filtration $\{\mathcal{F}_{t}\}_{t\in[0,T]}$, and endowed with the norms
                     \begin{equation*}
                     \left\|v \right\|_{L^{p}(\Omega, L^{r}(0,T; X))} =\left(\mathbb{E}\left(\int_{0}^{t} \left\|v \right\|^{r}_{X} dt\right)^{p/r}  \right)^{1/p}
                    \end{equation*}
           and
                              \begin{equation*}
                          \left\|v \right\|_{L^{p}(\Omega, L^{\infty}(0,T; X))} =\left(\mathbb{E}\sup\limits_{t\in[0,T]} \left\|v \right\|^{p}_{X}\right)^{1/p} , \ \ \mbox{if} \ \ r=\infty,
                              \end{equation*}
   where $\mathbb{E}$ is the mathematical expectation with respect to the probability measure $\mathbb{P}$. As usual, in the notation for processes $v=v(\omega,t)$, generally omit the dependence on $\omega\in \Omega$.\newline  
   We mention some useful interpolation inequalities, which we use frequently in the sequel. Let us first define a version of the Gagliardo-Nirenberg inequality which holds true for all $\text{u}\in \mathbb{W}^{1,p}_{0}(\text{D},\mathbb{R}^{n}), p\geq 1.$
         \begin{lem}(Gagliardo-Nirenberg inequality)\cite{TDM2}\label{laaz}. Let $\text{D}\subset\mathbb{R}^{n}$ and $\text{u}\in \mathbb{W}^{1,p}_{0}(\text{D},\mathbb{R}^{n}), p\geq 1.$ Then for any fixed number $1\leq q,r\leq \infty$, there exists a constant $C> 0$ depending only on $n,p,q$ such that:
                   \begin{equation}
               \left\|\text{u} \right\|_{\mathbb{L}^{r}}\leq C \left\|\nabla\text{u} \right\|^{\theta}_{\mathbb{L}^{p}}\left\|\text{u} \right\|^{1-\theta}_{\mathbb{L}^{q}},\quad \theta\in\left[0,1 \right],
                 \end{equation}
          where the numbers $p,q,r$ and $\theta$ satisfy the relation
                                    \begin{equation*}
                                 \theta=\left(\frac{1}{q}-\frac{1}{r}\right)\left(\frac{1}{n}-\frac{1}{p}+\frac{1}{q}\right)^{-1}.
                                  \end{equation*}
                                \end{lem}
        The particular cases of Lemma $\ref{laaz}$ are well known inequalities, thanks to the Ladyzhenskaya, which is given below.
                       \begin{lem}(Ladyzhenskaya's inequality)\cite{TDM2}. For  $\text{u}\in \mathcal{C}^{\infty}(\text{D},\mathbb{R}^{n}), n=2,3$ , there exists a constant $C> 0$ such that
                           \begin{equation}
                         \left\|\text{u} \right\|_{\mathbb{L}^{4}}\leq C \left\|\nabla\text{u} \right\|^{\frac{n}{4}}\left\|\text{u} \right\|^{1-\frac{n}{4}},\quad \text{for}\ \ n=2,3,
                          \end{equation}
          where $C=2,4$ for $n=2,3$ respectively.
          \begin{lem}(Gagliardo-Nirenberg interpolation inequality)\cite{TDM2}\label{laazo}. Let $\text{D}\subset\mathbb{R}^{n}$ and $\text{u}\in \mathbb{W}^{m,p}_{0}(\text{D},\mathbb{R}^{n}), p\geq 1.$ Then for any fixed number $1\leq q,r\leq \infty$ and a natural number $m$. Suppose also that a real number $\theta$ and a natural number $j$ are such that 
           \begin{equation*}
          \theta=\left(\frac{j}{n}+\frac{1}{q}-\frac{1}{r}\right)\left(\frac{m}{n}-\frac{1}{p}+\frac{1}{q}\right)^{-1}.
           \end{equation*}
          and $\frac{j}{m}\leq \theta < 1$. Then for any $\text{u}\in \mathbb{W}^{m,p}_{0}(\text{D},\mathbb{R}^{n})$, we have 
                   \begin{equation}
                         \left\|\nabla^{j}\text{u} \right\|_{\mathbb{L}^{r}}\leq C \left( \left\|\nabla^{m}\text{u} \right\|^{\theta}_{\mathbb{L}^{p}}\left\|\text{u} \right\|^{1-\theta}_{\mathbb{L}^{q}}+\left\|\text{u} \right\|_{\mathbb{L}^{s}}\right), 
                           \end{equation}
               where $s> 0$ is arbitrary and the constant $C> 0$ depends upon the domain $\text{D},m,n$.
                                          \end{lem}
                          \end{lem}
         \begin{lem}(Agmon's inequality)\cite{TDM2}.
        For any $u\in \mathbb{H}^{s_{2}}(\text{D},\mathbb{R}^{n})$, choose $s_{2}$ and $s_{2}$ such that $s_{1}< \frac{n}{2}< s_{2}$. Then, if $0<\delta< 1$ and $\frac{n}{2}=\delta s_{1}+(1-\delta)s_{2}$, the following inequality holds
                          \begin{equation}
                          \left\|\text{u} \right\|_{\mathbb{L}^{\infty}}\leq C \left\|\text{u} \right\|^{\delta}_{\mathbb{H}^{s_{1}}}\left\|\text{u} \right\|^{1-\delta}_{\mathbb{H}^{s_{2}}}.
                         \end{equation}
        For $u\in\mathbb{H}^{2}(\text{D})\cap\mathbb{H}^{1}_{0}(\text{D}) $, the Agmon inequality in $2\text{D}$ states that there exists a constant $C> 0$ such that:
                         \begin{equation}
                         \left\|\text{u} \right\|_{\mathbb{L}^{\infty}}\leq C \left\|\text{u} \right\|^{1/2}\left\|\text{u} \right\|^{1/2}_{\mathbb{H}^{2}}\leq C\left\|\text{u} \right\|_{\mathbb{H}^{2}}.
                            \end{equation}
                          \end{lem}
        In the following, we will make repeated use of Young's inequality
                          \begin{equation*}
                         ab\leq \varepsilon a^{2}+ \frac{1}{4\varepsilon}b^{2},\ \ \forall a, b\in\mathbb{R}\ \ \mbox{and},\ \ \varepsilon>0.
                           \end{equation*}
         We define by:
                               \begin{equation}
                               v_{D}=\int_{D}vdx.
                               \end{equation}
          It follows that for any vector $v\in \mathbb{V}_{div}$, we have $v_{D}=0$, since
                               \begin{equation}
                               \int_{D}v_{j}dx=\int_{D}div(vx_{j})dx=\int_{\Gamma}x_{j}(v.n)d\gamma=0, \ \ j=1,2.
                               \end{equation}
           Then, we have from \cite{NSC}, the following result.
                             \begin{lem}
           For any $v\in \mathbb{H}^{1}(D)$ and any $q\geq 2$, the Gagliardo-Nirenberg inequality
                             \begin{equation}\label{s1}
                             \left\|v-v_{D} \right\|_{\mathbb{L}^{q}}\leq C \left|v \right|^{2/q}_{\mathbb{L}^{2}} \left|\nabla v \right|^{1-2/q}_{\mathbb{L}^{2}},
                             \end{equation}
                             and the trace  interpolation inequality
                             \begin{equation}\label{s10}
                              \left\|v-v_{D} \right\|_{\mathbb{L}^{q}(\Gamma)}\leq C \left|v \right|^{2/q}_{\mathbb{L}^{2}} \left|\nabla v \right|^{1-2/q}_{\mathbb{L}^{2}},
                           \end{equation}
          are valid. Moreover, any $v\in \mathbb{V}_{div} $ satisfies Korn's inequality
                           \begin{equation}\label{s12}
                           \left\|v \right\|_{\mathbb{H}^{1}} \leq C\left\|v \right\|,
                           \end{equation}
       that is the norms $\left\|. \right\|_{\mathbb{H}^{1}}$ and $\left\|.\right\|$ are equivalent.
                             \end{lem}
                           \begin{prev}
                           See \cite[p. 3]{NSC}.
                           \end{prev}
          
       Next, we consider the formula
                           \begin{equation*}
                           -\int_{D}\Delta v.z dx=2\int_{D} Dv.Dz dx-\int_{\Gamma}2(n.v).zd\gamma,
                           \end{equation*}
      which follows from Green's theorem and is valid for any $v\in \mathbb{H}^{2}(D)\cap\mathbb{V}_{div}$ and $z\in \mathbb{H}^{1}$.
       From the boundary conditions, it follows that
                           \begin{equation*}
                           -\int_{D}\Delta v.z dx=2((v ,z))-\int_{\Gamma}b(z.\tau)d\gamma, \ \ \mbox{for any }\quad v\in \mathbb{V}_{div}\quad \mbox{and} \ \ z\in \mathbb{H}_{div}.
                            \end{equation*}
      Throughout the article, we often use the continuous embedding results
                           \begin{equation*}
                           \mathbb{H}^{1}(0,T)\subset \mathcal{C}([0,T]), \ \ \mathbb{H}^{1}(D)\subset L^{2}(\Gamma).
                           \end{equation*}
         Since $\mathbb{H}_{div}$ is a Hilbert space with inner product $(.,.)$, then for a vector
                           \begin{equation*}
                           h=(h_{1},h_{2},...,h_{m})\in\mathbb{H}^{m}_{div}=\underbrace{\mathbb{H}_{div}\times...\times\mathbb{H}_{div}}_{m-times},
                           \end{equation*}
      we introduce the norm and the absolute value of the inner product of $h$ with a fixed $v\in \mathbb{H}_{div}$ as
                           \begin{equation*}
                           \left|h \right|_{L^{2}}=\sum_{k=1}^{m}\left|h_{k} \right|_{L^{2}}, \ \ \mbox{and}\ \ \left|(h,v) \right|=\left(\sum_{k=1}^{m}(h_{k},v)^{2} \right)^{1/2}.
                           \end{equation*}
           We suppose that the stochastic noise is represented by
                           \begin{equation*}
                           g(t,y)dW_{t}=\sum_{k=1}^{m}g^{k}(t,y)dW^{k}_{t}
                           \end{equation*}
      where $g(t,y)=(g^{1}(t,y),...,g^{m}(t,y))$ has suitable growth assumptions, as defined in the following, and $W_{t}=(W^{1}_{t},...,W^{m}_{t})$ is a standard $\mathbb{R}^{m}$-valued Wiener process defined on a probability space.\\ 
      We will also denote by $C$ a generic positive constant that depend on the domain $\text{D}$ and on a given time $T$.
     To simplify the notations, we set (without loss of generality) $\mathcal{K}=\nu=\alpha=\delta=1$.\newline
                             We set
                             \begin{equation}
                             \mathbb{Y}=\mathbb{H}_{div}\times\text{V}_{1}=\mathbb{H}_{div}\times H^{1}(\text{D}),\ \ \mathbb{Y}_{1}=\mathbb{L}^{2}(\text{D})\times\text{V}_{1}.
                             \end{equation}
          The space $\mathbb{Y}$ is a complete metric space with respect to the norm
                             \begin{equation}\label{N16}
                              \left\| (\text{w},\psi)\right\| ^{2}_{\mathbb{Y}}= |\text{w}|^{2}_{L^{2}} + \left\|\psi\right\| ^{2}_{1}=|\text{w}|^{2}_{L^{2}}+|\nabla\psi|^{2}_{L^{2}}. 
                              \end{equation}
            Also, we define another Hilbert space $\mathbb{V}$ by:
                                                       \begin{equation}\label{N17}
                                                        \mathbb{V}=\mathbb{V}_{div}\times \text{V}_{2},\ \ \mathbb{V}_{1}=\mathbb{H}^{1}(\text{D})\times \text{V}_{2}.
                                                          \end{equation}
              endowed with the scalar product whose associated norm is given by
                                             \begin{equation}\label{N18}
                                                  ||(\text{w},\psi)||^{2}_{\mathbb{V}}= ||\text{w}||^{2} + \left\| \psi\right\| ^{2}_{2}.
                                                 \end{equation}
               Let us introduce the space of functions
                                     \begin{math}
                                     \mathcal{H}_{p}(\Gamma)=\{(a,b): \left\|(a,b) \right\|_{\mathcal{H}_{p}(\Gamma)}<\infty  \}
                                     \end{math}
                                     with\\ the norm
                                     \begin{equation*}
                                     \left\|(a,b) \right\|_{\mathcal{H}_{p}(\Gamma)}=\left\| a\right\|_{W_{p}^{1-\frac{1}{p}}(\Gamma)} +\left\|\partial_{t} a \right\|_{\mathbb{H}^{1/2}(\Gamma)} +\left\| b\right\|_{W_{p}^{-1/p}(\Gamma)} +\left\|b \right\|_{\mathbb{L}^{2}(\Gamma)} +\left\| \partial_{t} b\right\|_{\mathbb{H}^{-1/2}(\Gamma) } .
                                     \end{equation*}
               Now, we introduce some assumptions used throughout the paper:\\
                                     \textbf{(H1)} We assume that $g: [0,T]\times \mathbb{H}_{div}\rightarrow \mathbb{H}^{m}_{div}$  is Lipschitz in the second variable $\text{v}$ and satisfies a linear growth
                                     \begin{equation}\label{l1}
                                     \begin{aligned}
                                     \left|g(t,v)-g(t,w) \right|^{2}_{L^{2}}&\leq K\left|v-w \right|^{2}_{L^{2}},\\
                                       \left|g(t,v) \right|^{2}_{L^{2}}&\leq K(1+\left|v \right|^{2}_{L^{2}}), \ \ \forall v,w\in\mathbb{H}_{div}, t\in [0,T]
                                     \end{aligned}
                                     \end{equation}
                 for some positive constant $K$. \\
             \textbf{(H2)} We consider the data $(a,b)$ and $(u_{0},\phi_{0})$ that belong to the following Banach spaces
                                     \begin{equation}\label{l2}
                                     \begin{aligned}
                                     (a,b)&\in L^{2}(\Omega\times(0,T); \mathcal{H}_{p}(\Gamma) ) \ \ \mbox{for given}\ \ p\in (2,\infty),\\
                                     (\text{u}_{0},\phi_{0})&\in L^{2}(\Omega,\mathbb{Y}).
                                     \end{aligned}
                                     \end{equation}
             In addition, we suppose that $(a,b)$  is a pair of predictable stochastic processes. 
          Hereafter, for any $(\text{u},\phi)\in \mathbb{Y}$, we set
                                           \begin{equation}
                                          \mathcal{E}(\text{u},\phi)=\left\| (\text{u},\phi)\right\| ^{2}_{\mathbb{Y}}+ \int_{D}F(\phi).
                                            \end{equation}
          We can check that (see \cite{Ta1}) there exists a monotone non-decreasing function (independent on time and the initial condition) such that:
                                        \begin{equation}
                                        \left\| (\theta,\psi)\right\| ^{2}_{\mathbb{Y}}\leq \mathcal{E}(\theta,\psi)\leq \mathcal{Q}_{0}(\left\| (\theta,\psi)\right\| ^{2}_{\mathbb{Y}}), \quad \forall\, (\theta,\psi)\in \mathbb{Y}.
                                        \end{equation} 
           \section{Existence of the strong probabilistic solution}
          In this section, we investigate the existence of the strong probabilistic solution of problem $\eqref{NS1}$.  As in \cite{NSC}, we introduce a suitable change of variables based on the solution of the non homogeneous linear Stokes equation, which allows to write the state in terms of a vector field satisfying a homogeneous Navier-slip boundary conditions. 
                    \subsection{Well-posedness of Stokes equation}  We consider the stokes equation with non homogeneous Navier-slip boundary condition given by:
                    \begin{equation}\label{S1}
                   \begin{cases}
                   -\Delta \textbf{a}+\nabla p=0, \ \
                   div\ \  \textbf{a}=0 \ \ in\ \  D,\\
                   \textbf{a}\cdot n=a,\ \ \left[2D(\textbf{a})n+\alpha\textbf{a} \right]\cdot\tau=b \ \ on \ \ \Gamma.
                   \end{cases}
                    \end{equation}
         The following theorem gives the result the existence and uniqueness of the solution of problem $\eqref{S1}$.
                    \begin{thm}\label{th01}
        Let $(a,b)$ be given pair of functions satisfying $\eqref{l2}_{1}$. Then problem $\eqref{S1}$ admits a unique solution
                    \begin{equation}\label{hj1}
                    \textbf{a}\in L^{2}(\Omega, \mathbb{H}^{1}((0,T)\times \text{D}))\cap L^{2}(\Omega\times(0,T), W^{1}_{p}(\text{D})),
                    \end{equation}
                    such that
                    \begin{equation}\label{hj}
                  \left\| \textbf{a}\right\|_{W_{p}^{1}(\text{D})} +\left\|\partial_{t} \textbf{a} \right\|_{\mathbb{L}^{2}(\text{D})}\leq C \left\|(a,b) \right\|_{\mathcal{H}_{p}(\Gamma)}\ \ a.e \ \ in\ \ \Omega \times(0,T).
                  \end{equation}
         In particular, we have
                  \begin{equation*}
                            \textbf{a}\in L^{2}(\Omega, \mathcal{C}([0,T]; \mathbb{L}^{2}(\text{D})))\cap L^{2}(\Omega\times(0,T); \mathcal{C}(\overline{D})\cap\mathbb{H}^{1}(\text{D})).
                            \end{equation*}
                    \end{thm}
                    \begin{prev}
                    See \cite[Lemma 3.1, p. 5]{NSC}.
                    \end{prev}
          
      With the help of the solution of the non homogeneous Stokes equation, we give the definition of variation solution of  problem $\eqref{NS1}$.
                           \begin{defn}
       We assume that $(a,b)$ and $(\text{u}_{0},\phi_{0})$ satisfies the regularity $\eqref{l2}$, and \textbf{a} be the corresponding solution of $\eqref{S1}$. A process $(\text{u}+\textbf{a},\phi)=(\text{v},\phi)$  with $(\text{u},\phi)\in \mathcal{C}([0,T]; \mathbb{Y})\cap L^{2}(0,T; \mathbb{V})$, $\mathbb{P}-$a.e. in $\Omega$, is a strong (in the probabilistic sense) solution of the stochastic AC-NS system $\eqref{NS1}$ with $\text{v}_{0}=\text{u}_{0}+\textbf{a}(0)$, if $\mathbb{P}-$a.e. in $\Omega$, the following equalities hold:
                                              \begin{equation}\label{Ton2}
                                              \begin{cases}
                                              (\text{v}(t),\text{w})=(\text{v}_{0},\text{w}) +\int_{0}^{t}\left[  -((\text{v}(s),\text{w}))+\int_{\Gamma}b(\text{w}\cdot\tau)d\gamma - \left\langle (\text{v}(s)\cdot\nabla\text{v}(s)),\text{w}\right\rangle  \right]   ds\\
                                              +\int_{0}^{t}\left\langle (\mu(s)\nabla\phi(s), \text{w}\right\rangle ds
                                              +\int^{t}_{0}(\text{g}(s,\text{v}(s)),\text{w})dW_{s},\\
                                              (\phi(t),\psi)=(\phi_{0},\psi) -\int_{0}^{t}\left( \mu(s),\psi\right)  ds- \int_{0}^{t}\left\langle (\text{v}(s)\cdot\nabla\phi(s)),\psi\right\rangle   ds,\\
                                                \ \
                                               \mu=\text{A}_{\vartheta}\phi + \text{f}_{\vartheta}(\phi),
                                                \end{cases}
                                               \end{equation}
     for  all $t\in[0,T]$, $\text{w}\in\mathbb{V}_{div}$, $\psi\in\text{V}_{1}$, where the stochastic integral is understood  in the It\^{o} sense and defined by:
                            \begin{equation*}
                            \int^{t}_{0}(\text{g}(s,\text{v}(s)),\text{w})dW_{s}=\sum_{k=1}^{m}\int^{t}_{0}(\text{g}^{k}(s,\text{v}(s)),\text{w})dW^{k}_{s}.
                            \end{equation*}
                            \end{defn}
         In addition, we suppose that
                                          \begin{equation}\label{n1}
                                          \mathbb{E}[\mathcal{E}(\text{u}_{0},\phi_{0})]^{2}<\infty, \ \ (a,b)\in L^{4}(\Omega\times(0,T); \mathcal{H}_{p}(\Gamma)).
                                          \end{equation}
                                          The main result of this section is stated as follows:
                               \begin{thm}\label{Theo2}
      Suppose that $(a,b)$ and $(\text{u}_{0},\phi_{0})$ satisfy the regularity $\eqref{l2}$ and $\eqref{n1}$. Then, there exists, a unique solution $(\text{v},\phi)$ with $\text{v}=\text{u}+\textbf{a}$ to the system $\eqref{NS1}$ such that:
                                 \begin{equation*}
                                              (\text{u},\phi)\in\mathcal{C}([0,T]; \mathbb{Y})\cap L^{4}(0,T;\mathbb{V}), \ \ \mathbb{P}- \mbox{a.e. in}\ \ \Omega,
                                            \end{equation*}
                      and any $t\in[0,T]$, the following estimate holds
                                           \begin{equation}\label{g50a}
                                         \begin{aligned}
                                            \mathbb{E}\sup\limits_{s\in[0,t]}\mathcal{G}^{2}(s)\left\|(\text{u},\phi)(s) \right\|^{2}_{\mathbb{Y}}&+\mathbb{E}\int_{0}^{t}\mathcal{G}^{2}(s)\left[ \left\|(\text{u},\phi)(s) \right\|^{2}_{\mathbb{V}} \right] ds \\
                                         &\leq C\left(\mathbb{E}\mathcal{E}(\text{u}_{0},\phi_{0})+\mathbb{E}\int_{0}^{t}\mathcal{G}^{2}(s)\Lambda(s)ds \right),
                                        \end{aligned}
                                        \end{equation}
                                        and
                                      \begin{equation}\label{n2b0}
                                       \begin{aligned}
                                      \mathbb{E}\sup\limits_{s\in[0,t]}\mathcal{G}^{4}(s)\left\|(\text{u},\phi)(s) \right\|^{4}_{\mathbb{Y}} &+\mathbb{E}\left( \int_{0}^{ t}\mathcal{G}^{2}(s)\left[ \left\|(\text{u},\phi)(s)\right\|^{2}_{\mathbb{V}}\right] ds\right) ^{2}\\
                                      &\leq C\mathbb{E}\left[ \mathcal{E}(\text{u}_{0},\phi_{0})\right]^{2}
                                      +C\mathbb{E}\int_{0}^{t} \mathcal{G}^{4}(s)B(s)ds,
                                       \end{aligned}
                                     \end{equation}
               where the functions $\mathcal{G}$ and  $B,\Lambda$ are defined by $\eqref{g4}$ and by $\eqref{g6}$, $\eqref{g6a}$, respectively. Here $C$ is a positive constant. 
                                         \end{thm}   
                                          
                    \section{Proof of Theorem $\ref{Theo2}$ and Uniqueness of the solution  }\label{ft0}
     In this section,we prove Theorem $\ref{Theo2}$ using the Galerkin approximation and some ideas from the following papers \cite{Ta1}, \cite{Zhang}.
                      \subsection{Galerkin approximation}   
           Since the injection $\mathbb{V}\hookrightarrow\mathbb{Y}$ is a compact, there exists a basis $e_{i}=\{(w_{i},\psi_{i})\}\subset \mathbb{V}$ of eigenfunction verifying the property
                                     \begin{equation}\label{g1}
                                     ((y,e_{i}))=\lambda_{i}((y,e_{i})), \ \ \forall y\in \mathbb{V}, i\in\mathbb{N},
                                     \end{equation}
         which is an orthonormal basis of $\mathbb{Y}$, and the corresponding sequence $\{\lambda_{i}\}$ of eigen-values verifies $\lambda_{i}>0$ , $\forall\, i\in\mathbb{N}$ and $\lambda_{i}\rightarrow \infty$ as $i\rightarrow\infty$. For the  detail we refer in (\cite[ Theorem 1, p. 355]{EVAN}). Moreover, the ellipticity of the equation $\eqref{g1}$ and the regularity $\Gamma\in\mathcal{C}^{2}$ imply that $e_{i}=\{(w_{i},\psi_{i})\}\subset \mathcal{C}\cap\mathbb{V}$.
       For any fixed $n\in\mathbb{N}$, we consider the subspace $\mathbb{V}_{n}=span\{(w_{1},\psi_{1}),...,(w_{n},\psi_{n})\}$ of $\mathbb{V}$. Taking into account the relation $\eqref{g1}$, the approximation finite dimensional problem is: for $\mathbb{P}$- a.e. in $\Omega$ to find $(\text{v}_{n},\phi_{n})$ in the from
                                     \begin{equation*}
                                     \begin{aligned}
                                     \text{v}_{n}=\text{u}_{n}+\textbf{a} \ \ \mbox{with}\ \ \text{u}_{n}=\sum_{j=1}^{n}\beta_{j}^{n}w_{i},\ \
                                     \phi_{n}=\sum_{j=1}^{n}\chi_{j}^{n}\psi_{i}, \mu_{n}=\sum_{j=1}^{n}\chi_{j}^{n}\psi_{i} \ \ \mbox{with}\ \ t\in [0,T],
                                     \end{aligned}
                                     \end{equation*}
        as the solution of the following finite dimensional stochastic differential equation
                                     \begin{equation}\label{g2}
                                     \begin{cases}
                                     d(\text{v}_{n},\varphi)=\left[-((\text{v}_{n},\varphi))+\int_{\Gamma}b(\varphi\cdot\tau)d\gamma-\left\langle (\text{v}_{n}\cdot\nabla)\text{v}_{n},\varphi\right\rangle  \right]dt \\
                                     \left\langle \mu_{n}\nabla\phi_{n},\varphi\right\rangle dt+(g(t,\text{v}_{n}), \varphi)dW(t),\\
                                     d(\phi_{n},\psi)=\left[\left\langle \mu_{n}-\text{v}_{n}\cdot\nabla\phi_{n},\psi \right\rangle \right]dt,\ \  \forall t\in(0,T), (\varphi,\psi)\in\mathbb{V}_{n},\\
                                     \mu_{n}=\text{A}_{\vartheta}\phi_{n}+f_{\vartheta}(\phi_{n}), \ \ (\text{u}_{n},\phi_{n})(0)=(\text{u}_{n,0},\phi_{n,0}),
                                     \end{cases}
                                     \end{equation}
               where $(\text{u}_{n,0},\phi_{n,0})=\sum_{j=1}^{n}((\text{u}_{0},\phi_{0}),e_{j})e_{j}$ is the orthogonal projection of $(\text{u}_{0},\phi_{0})\in \mathbb{Y}$ into the space $\mathbb{V}_{n}$. Using the Parseval identity we deduce that
                                     \begin{equation}\label{g3}
                                     \left\| (\text{u}_{n,0},\phi_{n,0})\right\|_{\mathbb{Y}}< \left\|(\text{u}_{0},\phi_{0}) \right\|_{\mathbb{Y}} \ \ \mbox{and} \ \  (\text{u}_{n,0},\phi_{n,0})\rightarrow (\text{u}_{0},\phi_{0}) \ \ \mbox{strongly in }\ \ \mathbb{Y}.
                                     \end{equation}
            Since the problem $\eqref{g2}$ defines a system of $n$ stochastic ordinary differential equations with locally and globally Lipschitz nonlinearities. Then, the problem $\eqref{g2}$ has a unique local strong $(\text{u}_{n},\phi_{n})\in\mathcal{C}([0,T_{n}]; \mathbb{V}_{n})$, $\mathbb{P}-$a.e. in $\Omega$ (see for example \cite{KIS}). Next, we derive the a-priori energy estimates satisfied by the system $\eqref{g2}$ (see Lemma $\ref{La1}$ and $\ref{La2}$).   
            \begin{lem}\label{La1}
         Assume that $(a,b)$ and $(\text{u}_{0},\phi_{0})$  satisfy the regularity $\eqref{l2}$ and $\mathbb{E}[\mathcal{E}(\text{u}_{0},\phi_{0})]<\infty$. Then the system $\eqref{g2}$ has a solution $(\text{v}_{n},\phi_{n})$ with $\text{v}_{n}=\text{u}_{n}+\textbf{a}$ such that
                                       \begin{equation*}
                                       (\text{u}_{n},\phi_{n})\in\mathcal{C}([0,T]; \mathbb{Y})\cap L^{2}(0,T;\mathbb{V}), \ \ \mathbb{P}-  \mbox{a.e. in}\ \ \Omega.
                                       \end{equation*}
             Moreover, the exists a positive constant $C_{0}$, such that for the function
                                       \begin{equation}\label{g4}
                                       \mathcal{G}(t)= \exp\left(-C_{0}t-C_{0}\int_{0}^{t}\tilde{f}(s) ds \right) \ \ \mathbb{P}- \mbox{a.e. in}\ \ \Omega,
                                       \end{equation}
                and any $t\in[0,T]$, the following estimate holds
                                       \begin{equation}\label{g5}
                                       \begin{aligned}
                                       \mathbb{E}\sup\limits_{s\in[0,t]}\mathcal{G}^{2}(s)\left\|(\text{u}_{n},\phi_{n})(s) \right\|^{2}_{\mathbb{Y}}&+\mathbb{E}\int_{0}^{t}\mathcal{G}^{2}(s) \left\|(\text{u}_{n},\phi_{n})(s)\right\|^{2}_{\mathbb{V}} ds \\
                                       &\leq C\left(\mathbb{E}\mathcal{E}(\text{u}_{0},\phi_{0})+\mathbb{E}\int_{0}^{t}\mathcal{G}^{2}(s)\Lambda(s)ds \right),
                                       \end{aligned}
                                       \end{equation}
           where
                                       \begin{equation}\label{g6}
                                       \tilde{f}(s)=\left(\left|\textbf{a} \right|^{2}_{L^{2}} \left\|\textbf{a} \right\|^{2}  +\Lambda\right) \in L^{1}(\Omega\times(0,T)), \ \ \Lambda=\left\|(a,b) \right\|^{2}_{\mathcal{H}_{p}(\Gamma)}+1\in L^{1}(\Omega\times(0,T)),
                                       \end{equation}
       and the positive constant $C$ and $C_{0}$ are independent of the parameter $n$, which may depend on the regularity of the boundary $\Gamma$ and the physical constant $\alpha.$
                                       \end{lem}
                                       \begin{prev}
                Let $\mathcal{G}$ be the function defined by $\eqref{g4}$  with a constant $C_{0}$ to be selected later. 
                Now, taking $\varphi=w_{i}$ for each $i=1,...,n$ in the equation $\eqref{g2}_{1}$ and using $\text{v}_{n}=\text{u}_{n}+\textbf{a}$, we derive that
                                           \begin{equation*}
                                           \begin{aligned}
                                           d(\text{u}_{n},w_{i})&=\left[-((\text{u}_{n}+\textbf{a},w_{i}))+\int_{\Gamma}b(w_{i}\cdot\tau)d\gamma-((\text{u}_{n}+\textbf{a})\cdot\nabla)(\text{u}_{n}+\textbf{a}),w_{i}) \right]dt \\
                                            &+(-\partial_{t}\textbf{a}+\mu_{n}\nabla\phi_{n},w_{i})dt+(g(t,\text{v}_{n}), w_{i})dW(t).
                                           \end{aligned}
                                           \end{equation*}
             Using It\^{o}'s formula, we derive that
                  \begin{equation}\label{A09}
                  \begin{aligned}
                  d \left|\text{u}_{n}\right|^{2}_{L^{2}}&+2\left\|\text{u}_{n}\right\|^{2}dt\\ &=\left[-2((\textbf{a},\text{u}_{n}))+\int_{\Gamma}\left[-a(\text{u}_{n}\cdot\tau)^{2}+ 2b(\text{u}_{n}\cdot\tau)\right] d\gamma-2(((\text{u}_{n}+\textbf{a})\cdot\nabla)\textbf{a},\text{u}_{n}) \right]dt \\
                 &+2(-\partial_{t}\textbf{a}+\mu_{n}\nabla\phi_{n},\text{u}_{n})dt+2(g(t,\text{v}_{n}), \text{u}_{n})dW_{t}+ \sum\limits^{n}_{i=1}\left|(g(t,\text{v}_{n}), w_{i}) \right|^{2}_{L^{2}}dt,
                  \end{aligned}
                  \end{equation}
         where we used the fact that \begin{math}
                                 \int_{D}\left[ (\text{u}_{n}+\textbf{a})\cdot\nabla\text{u}_{n}\right] \text{u}_{n}dx=\int_{\Gamma}\frac{a}{2}(\text{u}_{n}\cdot\tau)^{2}d\gamma.\\
                                 \end{math}
         Taking inner product of $\eqref{g2}_{2}$ with $\mu_{n}$, we obtain
                                 \begin{equation}\label{A10}
                                 \begin{aligned}
                                 d\left( \left|\nabla\phi_{n}\right|^{2}_{L^{2}}+2\int_{D}F_{\vartheta}(\phi_{n}))\right) +2\left|\mu_{n} \right|^{2}_{L^{2}}dt  &=2\left[-((\text{u}_{n}+\textbf{a})\cdot\nabla\phi_{n},\mu_{n}) \right]dt.
                                 \end{aligned}
                                 \end{equation} 
            It follows from $\eqref{A09}-\eqref{A10}$ that
                           \begin{equation}\label{A11}
                            \begin{aligned}
                             &\mathcal{E}(\text{u}_{n},\phi_{n})(t)+2\int_{0}^{t}\left[ \left\|\text{u}_{n}(s)\right\|^{2}+\left|\mu_{n}(s) \right|^{2}_{L^{2}}\right] ds=\mathcal{E}(\text{u}_{n},\phi_{n})(0)\\
                            &+\int_{0}^{t}\left[-2((\textbf{a},\text{u}_{n}(s)))+\int_{\Gamma}\left[-a(\text{u}_{n}(s)\cdot\tau)^{2}+ 2b(\text{u}_{n}(s)\cdot\tau)\right] d\gamma\right]ds\\
                            &-2\int_{0}^{t}\left[\partial_{t}\textbf{a}+ (((\text{u}_{n}(s)+\textbf{a})\cdot\nabla)\textbf{a},\text{u}_{n}(s)) \right]ds 
                          -2\int_{0}^{t}(\textbf{a}\cdot\nabla\phi_{n},\mu_{n})ds\\
                          &+2\int_{0}^{t}(g(s,\text{v}_{n}(s)), \text{u}_{n})dW_{s}
                          + \sum\limits^{n}_{i=1}\int_{0}^{t}\left|(g(s,\text{v}_{n}(s)), w_{i}) \right|^{2}_{L^{2}}ds,\\
                          &=\mathcal{E}(\text{u}_{n},\phi_{n})(0)
                          +\int_{0}^{t}\left[ I_{1}(s) +I_{2}(s)+I_{3}(s)+I_{4}(s)\right]ds\\ &+2\int_{0}^{t}(g(s,\text{v}_{n}(s)), \text{u}_{n})dW_{s}.
                          \end{aligned}
                          \end{equation}
         where we used the fact that $((\text{u}_{n}\cdot\nabla\phi_{n},\mu_{n})=-(\mu_{n}\nabla\phi_{n},\text{u}_{n})$.\\
    We estimate each term in the right hand side of the equation $\eqref{A11}$ using the H\"{o}lder, Ladyzhenskaya, Young, Poincaré, Gagliardo--Nirenberg, Agmon inequalities, Sobolev embedding and Theorem $\ref{th01}$     
    \begin{equation}\label{es0}
                    \begin{aligned}
                    I_{1}&\leq 2\left\|\textbf{a} \right\|\left\|\text{u}_{n}(s)\right\|+\left|a \right|_{L^{\infty}(\Gamma)} \left|\text{u}_{n}(s) \right|^{2}_{L^{2}(\Gamma)}+ 2\left|b \right|_{L^{2}(\Gamma)} \left|\text{u}_{n}(s) \right|_{L^{2}(\Gamma)} \\
                    &\leq \frac{1}{2}\left\|\text{u}_{n}(s)\right\|^{2}+ C\left\|a \right\|^{2}_{W_{p}^{1-\frac{1}{p}}} \left|\text{u}_{n}(s) \right|^{2}_{L^{2}}+C(\left\|\textbf{a} \right\|^{2}+\left|b \right|_{L^{2}(\Gamma)})\\
                    &\leq \frac{1}{2}\left\|\text{u}_{n}(s)\right\|^{2}+C\Lambda\left|\text{u}_{n}(s) \right|^{2}_{L^{2}}+C\left\|(a,b) \right\|^{2}_{\mathcal{H}_{p}(\Gamma)},
                    \end{aligned}
                    \end{equation}
                    \begin{equation}\label{es1}
                   \begin{aligned}
                   I_{2}&\leq 2\left(\left|\partial_{t}\textbf{a} \right|_{L^{2}} +\left\|\textbf{a} \right\|_{\mathcal{C}(\overline{D})} \left| \nabla\textbf{a}\right|_{L^{2}}  \right)\left|\text{u}_{n}(s) \right|_{L^{2}}+2\left| \nabla\textbf{a}\right|_{L^{2}}\left\|\text{u}_{n}(s) \right\|^{2}_{\mathbb{L}^{4}}\\
                    &\leq \frac{1}{2}\left\|\text{u}_{n}(s)\right\|^{2}+C\Lambda\left|\text{u}_{n}(s) \right|^{2}_{L^{2}},
                     \end{aligned}
                    \end{equation}
                         \begin{equation}\label{est3}
                         \begin{aligned}
                         &I_{3}\leq \left\|\textbf{a} \right\|_{\mathbb{L}^{4}} \left|\mu_{n} \right|_{L^{2}}\left\|\nabla\phi_{n} \right\|_{L^{4}}
                         \leq \frac{1}{2}\left|\mu_{n} \right|^{2}_{L^{2}}+C\left\|\textbf{a} \right\|^{2}_{\mathbb{L}^{4}}\left\|\nabla\phi_{n} \right\|^{2}_{L^{4}}\\
                         &\leq \frac{1}{2}\left|\mu_{n} \right|_{L^{2}}+C\left|\textbf{a} \right|_{L^{2}} \left\|\textbf{a} \right\|\left( \left\|\phi_{n} \right\|_{2} \left|\nabla\phi_{n} \right|_{L^{2}}+\left|\nabla\phi_{n} \right|^{2}_{L^{2}} \right) \\
                         &\leq \frac{1}{2}\left|\mu_{n} \right|^{2}_{L^{2}}+C \left\|\textbf{a} \right\|^{2}\left|\nabla\phi_{n} \right|^{2}_{L^{2}}+ C\left|\textbf{a} \right|_{L^{2}} \left\|\textbf{a} \right\|\left\|\phi_{n} \right\|_{2} \left|\nabla\phi_{n} \right|_{L^{2}}\\
                         &\leq  \frac{1}{2}\left|\mu_{n} \right|^{2}_{L^{2}}+\frac{1}{4}\left\|\phi_{n} \right\|_{2}^{2}+C \Lambda\left|\nabla\phi_{n} \right|^{2}_{L^{2}}+C\left|\textbf{a} \right|^{2}_{L^{2}} \left\|\textbf{a} \right\|^{2} \left|\nabla\phi_{n} \right|^{2}_{L^{2}},
                         \end{aligned}
                         \end{equation}
                         and
                         \begin{equation}\label{est4}
                         \begin{aligned}
                         I_{4}&=\sum\limits^{n}_{i=1}\left|(g(s,\text{v}_{n}(s)), w_{i}) \right|^{2}_{L^{2}}\leq C\left|g(s,\text{v}_{n}(s)) \right|^{2}_{L^{2}} \leq C(1+\left|\text{v}_{n}(s) \right|^{2}_{L^{2}})\\
                         &\leq C(1+\left|\text{u}_{n}(s) \right|^{2}_{L^{2}}+\left|\textbf{a} \right|^{2}_{L^{2}})\leq C(\left|\text{u}_{n}(s) \right|^{2}_{L^{2}}+\Lambda).
                         \end{aligned}
                         \end{equation}
    It follows from $\eqref{es0}-\eqref{est4}$ that there exists a positive constant $C_{0}$ such that:
                         \begin{equation}\label{est5}
                         \begin{aligned}
                         \int_{0}^{t}\left[ I_{1}(s) +I_{2}(s)+I_{3}(s)+I_{4}(s)\right]ds&\leq\int_{0}^{t} 2C_{0}\Lambda(\left|\text{u}_{n}(s) \right|^{2}_{L^{2}}+\left|\nabla\phi_{n}(s) \right|^{2}_{L^{2}}+1)ds\\
                         &+C\int_{0}^{t}\left|\textbf{a} \right|^{2}_{L^{2}} \left\|\textbf{a} \right\|^{2} \left|\nabla\phi_{n} \right|^{2}_{L^{2}} ds.
                         \end{aligned}
                         \end{equation}
    Taking the function $\mathcal{G}$ as in $\eqref{g4}$ and  applying the It\^{o} formula, we derive using  $\eqref{es0}-\eqref{est5}$ that
                        \begin{equation}\label{A11a}
                          \begin{aligned}
                        &\mathcal{G}^{2}(s)\mathcal{E}(\text{u}_{n},\phi_{n})(s)+2\int_{0}^{s}\mathcal{G}^{2}(r)\left[ \left\|\text{u}_{n}(r)\right\|^{2}+\left|\mu_{n}(r) \right|^{2}_{L^{2}}\right] dr\\
                          &=\mathcal{E}(\text{u}_{n},\phi_{n})(0)
                         -2C_{0}\int_{0}^{s} \mathcal{G}^{2}(r)\tilde{f}(r) \mathcal{E}(\text{u}_{n},\phi_{n})(r)dr-C_{0}\int_{0}^{s}\mathcal{G}^{2}(r)\mathcal{E}(\text{u}_{n},\phi_{n})(r)dr\\
                         &+\int_{0}^{s}\mathcal{G}^{2}(r)\left[ I_{1}(r) +I_{2}(r)+I_{3}(r)+I_{4}(r)+I_{5}(r)\right]dr +2\int_{0}^{t}\mathcal{G}^{2}(r)(g(r,\text{v}_{n}(r)), \text{u}_{n})dW_{r}, %
                                        \end{aligned}
                                        \end{equation}
                         which implies
                                         \begin{equation}\label{A110}
                                         \begin{aligned}
                                         &\mathcal{G}^{2}(s)\mathcal{E}(\text{u}_{n},\phi_{n})(s)+\int_{0}^{s}\mathcal{G}^{2}(r)\left[ \left\|\text{u}_{n}(r)\right\|^{2}+\left|\mu_{n}(r) \right|^{2}_{L^{2}}\right] dr\\
                                         &\leq\mathcal{E}(\text{u}_{n,0},\phi_{n,0})
                                          +C\int_{0}^{s} \mathcal{G}^{2}(r)\Lambda(r)dr+\frac{1}{2}\int_{0}^{s} \mathcal{G}^{2}(r)\left\|\phi_{n}(r) \right\|_{2}^{2}dr\\
                                          &+2\int_{0}^{s}\mathcal{G}^{2}(r)(g(r,\text{v}_{n}(r)), \text{u}_{n})dW_{r}.
                                         \end{aligned}
                                          \end{equation} 
                       From $\eqref{g2}_{3}$, we note that
                       \begin{equation}\label{lm}
                       \text{A}_{\vartheta}\phi_{n}=  \mu_{n}- f_{\vartheta}(\phi_{n}). 
                       \end{equation}           
                       From $\eqref{kl}$, we also have 
                       \begin{equation}\label{lm0}
                     - \left(  f_{\vartheta}(\phi_{n}),\text{A}_{\vartheta}\phi_{n}\right) \leq \vartheta\left|\nabla\phi_{n}\right|^{2}_{L^{2}}.
                       \end{equation}
              Multipying $\eqref{lm}$ by $\text{A}_{\vartheta}\phi_{n}$, using $\eqref{lm0}$ and Young's inequality,  we derive that
              \begin{equation}\label{lm1}
              \left\|\phi_{n}\right\|^{2}_{2}\leq \left| \mu_{n}\right|^{2}_{L^{2}}  + 2\vartheta\left|\nabla\phi_{n}\right|^{2}_{L^{2}}.
              \end{equation}            
             Using $\eqref{lm1}$ and $\eqref{A110}$, we obtain 
            \begin{equation}\label{A110a}
              \begin{aligned}
             &\mathcal{G}^{2}(s)\mathcal{E}(\text{u}_{n},\phi_{n})(s)+\frac{1}{2}\int_{0}^{s}\mathcal{G}^{2}(r) \left\|(\text{u}_{n},\phi_{n})(r)\right\|^{2}_{\mathbb{V}} dr\\
             &\leq\mathcal{E}(\text{u}_{n,0},\phi_{n,0})
             +C\int_{0}^{s} \mathcal{G}^{2}(r)\Lambda(r)dr+2\vartheta\int_{0}^{s} \mathcal{G}^{2}(r)\left|\nabla\phi_{n}(r) \right|^{2}_{L^{2}}dr\\
                &+2\int_{0}^{s}\mathcal{G}^{2}(r)(g(r,\text{v}_{n}(r)), \text{u}_{n})dW_{r}.
                 \end{aligned}
                \end{equation}                                                    
           For $n\in\mathbb{N}$, let us define the stopping time by:
                                    \begin{equation}\label{g7}
                                  \tau_{N}^{n}=\inf\{t\geq 0: h(t)\geq N\}\wedge T,
                                      \end{equation}
                 where
                                     \begin{equation}\label{g8}
                                       h(t)=\mathcal{G}^{2}(t)\mathcal{E}((\text{u}_{n},\phi_{n})(t))+\int_{0}^{t}\mathcal{G}^{2}(s)\left\|(\text{u}_{n},\phi_{n})(r)\right\|^{2}_{\mathbb{V}} ds.
                                       \end{equation}  
       Let us consider the sequence $\tau^{n}_{N}$ of the stopping times introduce in $\eqref{g7}$. Thanks to \textbf{(H1)} and Burholder-Davis-Gundy's inequality, we infer that
                              \begin{equation*}
                              \begin{aligned}
                              &\mathbb{E}\sup\limits_{s\in[0,\tau^{n}_{N}\wedge t]}\left|\int_{0}^{s}\mathcal{G}^{2}(r)(g(r,\text{v}_{n}(r)), \text{u}_{n})dW_{r} \right|\leq\mathbb{E} \left(\int_{0}^{\tau^{n}_{N}\wedge t}\mathcal{G}^{4}(r)\left| (g(r,\text{v}_{n}(r)), \text{u}_{n})\right|^{2} ds \right)^{1/2} \\
                              &\leq \mathbb{E}\sup\limits_{s\in[0,\tau^{n}_{N}\wedge t]}\mathcal{G}(s)\left| \text{u}_{n}(s)\right|_{L^{2}}\left(\int_{0}^{\tau^{n}_{N}\wedge t}\mathcal{G}^{2}(r)\left| (g(r,\text{v}_{n}(r)))\right|^{2}_{L^{2}} ds \right)^{1/2}\\
                              &\leq \frac{1}{2}\mathbb{E}\sup\limits_{s\in[0,\tau^{n}_{N}\wedge t]}\mathcal{G}^{2}(s)\mathcal{E}(\text{u}_{n},\phi_{n})(s)+C\mathbb{E}\int_{0}^{\tau^{n}_{N}\wedge t}\mathcal{G}^{2}(s)\left(\mathcal{E}(\text{u}_{n},\phi_{n})(s)+\Lambda(s) \right)ds.
                              \end{aligned}
                              \end{equation*}
           For $t\in[0,T]$, we first take the supremum of the relation $\eqref{A110}$ for $s\in[0,\tau^{n}_{N}\wedge t]$, next we take the expectation, we conclude that
                              \begin{equation}\label{A11b}
                               \begin{aligned}
                                &\mathbb{E}\sup\limits_{s\in[0,\tau^{n}_{N}\wedge t]}\left[ \mathcal{G}^{2}(s)\mathcal{E}(\text{u}_{n},\phi_{n})(s)\right] +\mathbb{E}\int_{0}^{\tau^{n}_{N}\wedge t}\mathcal{G}^{2}(s)\left\|(\text{u}_{n},\phi_{n})(s)\right\|^{2}_{\mathbb{V}}ds\\
                                &\leq\mathbb{E}\mathcal{E}(\text{u}_{n,0},\phi_{n,0})
                                +C\mathbb{E}\int_{0}^{\tau^{n}_{N}\wedge t} \mathcal{G}^{2}(s)\Lambda(r)dr
                                +C\mathbb{E}\int_{0}^{\tau^{n}_{N}\wedge t}\mathcal{G}^{2}(s)\mathcal{E}(\text{u}_{n},\phi_{n})(s)ds.
                                 \end{aligned}
                                  \end{equation}
                It follows from the Gronwall lemma that
                                   \begin{equation}\label{A11b1}
                                    \begin{aligned}
                                    \mathbb{E}\sup\limits_{s\in[0,\tau^{n}_{N}\wedge t]}\left[ \mathcal{G}^{2}(s)\mathcal{E}(\text{u}_{n},\phi_{n})(s)\right] &+\mathbb{E}\int_{0}^{\tau^{n}_{N}\wedge t}\mathcal{G}^{2}(s)\left\|(\text{u}_{n},\phi_{n})(s)\right\|^{2}_{\mathbb{V}}  ds\\
                                     &\leq C\mathbb{E}\mathcal{E}(\text{u}_{n,0},\phi_{n,0})
                                     +C\mathbb{E}\int_{0}^{t} \mathcal{G}^{2}(s)\Lambda(s)ds.
                                   \end{aligned}
                                  \end{equation}
              As in \cite{Ta1} or in \cite{NSC}, using the fact the $\tau^{n}_{N}\uparrow T$ as $N$ goes to $\infty$, we obtain
                                  \begin{equation}\label{Q0}
                                  \begin{aligned}
                                   \mathbb{E}\sup\limits_{s\in[0,t]}\left[ \mathcal{G}^{2}(s)\mathcal{E}(\text{u}_{n},\phi_{n})(s)\right] &+\mathbb{E}\int_{0}^{ t}\mathcal{G}^{2}(s)\left\|(\text{u}_{n},\phi_{n})(s)\right\|^{2}_{\mathbb{V}}  ds\\
                                  &\leq C\mathbb{E}\mathcal{E}(\text{u}_{n,0},\phi_{n,0})
                                 +C\mathbb{E}\int_{0}^{t} \mathcal{G}^{2}(s)\Lambda(s)ds.
                                \end{aligned}
                                \end{equation}
         This completes the proof of Lemma $\ref{La1}$ .
             \end{prev}   
             
           In the next lemma, by assuming a better integrability for the initial data, we improve the integrability properties for the solution $(\text{v}_{n},\phi_{n})$ of the problem $\eqref{g2}$.
            \begin{lem}\label{La2}
         Suppose that $(a,b)$ and $(\text{u}_{0},\phi_{0})$  satisfy the regularity $\eqref{l2}$ and $\eqref{n1}$.
           Then, the solution $(\text{v}_{n},\phi_{n})$ with $\text{v}_{n}=\text{u}_{n}+\textbf{a}$ of problem $\eqref{g2}$ has the regularity:
                                    \begin{equation*}
                                   (\text{u}_{n},\phi_{n})\in\mathcal{C}([0,T]; \mathbb{Y})\cap L^{4}(0,T;\mathbb{V}), \ \ \mathbb{P}- \mbox{a.e. in}\ \ \Omega,
                                     \end{equation*}
                    such that
                                          \begin{equation}\label{g5a}
                                         \begin{aligned}
                                        \mathbb{E}\sup\limits_{s\in[0,t]}\mathcal{G}^{4}(s)\left\|(\text{u}_{n},\phi_{n})(s) \right\|^{4}_{\mathbb{Y}}&+\mathbb{E}\left( \int_{0}^{t}\mathcal{G}^{2}(s)\left\|(\text{u}_{n},\phi_{n})(s)\right\|^{2}_{\mathbb{V}} ds\right)^{2}  \\
                                        &\leq C\left(\mathbb{E}\left[ \mathcal{E}(\text{u}_{0},\phi_{0})\right] ^{2}+\mathbb{E}\int_{0}^{t}\mathcal{G}^{4}(s)B(s)ds \right),
                                       \end{aligned}
                                        \end{equation}
                where the function $\mathcal{G}$ is defined in $\eqref{g4}$,
                                        \begin{equation}\label{g6a}
                                        B=\left\|(a,b) \right\|^{4}_{\mathcal{H}_{p}(\Gamma)}+1\in L^{1}(\Omega\times(0,T)),
                                       \end{equation}
           and  $C$ is a positive constant, being independent of $n$. 
           \end{lem}
            \begin{prev}
          Taking the square on both side of the inequality $\eqref{A110a}$ and the supremum on  $s\in[0,\tau^{n}_{N}\wedge t]$ with $\tau^{n}_{N}$ of the stopping times introduce in $\eqref{g7}$, we see that
                              \begin{equation}\label{A1101}
                                \begin{aligned}
                                &\sup\limits_{s\in[0,\tau^{n}_{N}\wedge t]}\left[\mathcal{G}^{2}(s)\mathcal{E}(\text{u}_{n},\phi_{n})(s)\right]^{2} +\frac{1}{4}\left( \int_{0}^{\tau^{n}_{N}\wedge t}\mathcal{G}^{2}(s)\left\|(\text{u}_{n},\phi_{n})(s)\right\|^{2}_{\mathbb{V}} ds\right) ^{2}\\
                                &\leq\left[ \mathcal{E}(\text{u}_{n,0},\phi_{n,0})\right]^{2}
                                +C\int_{0}^{\tau^{n}_{N}\wedge t} \mathcal{G}^{4}(s)B(s)ds+4\vartheta^{2}\int_{0}^{\tau^{n}_{N}\wedge t} \mathcal{G}^{4}(s)\left|\nabla\phi_{n}(s) \right|^{4}_{L^{2}}ds\\
                                &+4\sup\limits_{s\in[0,\tau^{n}_{N}\wedge t]}\left|\int_{0}^{s}\mathcal{G}^{2}(r)(g(r,\text{v}_{n}(r)), \text{u}_{n})dW_{r} \right|^{2}.
                               \end{aligned}
                               \end{equation}
             Taking the expectation in this inequality and applying the Burkholder-Davis-Gundy inequality, we infer that
                                   \begin{equation}\label{k01}
                                    \begin{aligned}
                                   &\mathbb{E}\sup\limits_{s\in[0,\tau^{n}_{N}\wedge t]}\left|\int_{0}^{s}\mathcal{G}^{2}(r)(g(r,\text{v}_{n}(r)), \text{u}_{n})dW_{r} \right|^{2}\leq\mathbb{E} \left(\int_{0}^{\tau^{n}_{N}\wedge t}\mathcal{G}^{4}(r)\left| (g(r,\text{v}_{n}(r)), \text{u}_{n})\right|^{2} ds \right) \\
                                   &\leq \mathbb{E}\sup\limits_{s\in[0,\tau^{n}_{N}\wedge t]}\mathcal{G}^{2}(s)\left| \text{u}_{n}(s)\right|^{2}_{L^{2}}\left(\int_{0}^{\tau^{n}_{N}\wedge t}\mathcal{G}^{2}(r)\left| (g(r,\text{v}_{n}(r)))\right|^{2}_{L^{2}} ds \right)\\
                                  &\leq \frac{1}{2}\mathbb{E}\sup\limits_{s\in[0,\tau^{n}_{N}\wedge t]}\left[ \mathcal{G}^{2}(s)\mathcal{E}(\text{u}_{n},\phi_{n})(s)\right] ^{2}+C\mathbb{E}\int_{0}^{\tau^{n}_{N}\wedge t}\mathcal{G}^{4}(s)\left(\left[ \mathcal{E}(\text{u}_{n},\phi_{n})(s)\right]^{2} +B(s) \right)ds.
                                 \end{aligned}
                                 \end{equation}
            Inserting now $\eqref{k01}$ in $\eqref{A1101}$, it follows that
                             \begin{equation}\label{n2}
                              \begin{aligned}
                              &\mathbb{E}\sup\limits_{s\in[0,\tau^{n}_{N}\wedge t]}\left[\mathcal{G}^{2}(s)\mathcal{E}(\text{u}_{n},\phi_{n})(s)\right]^{2} +\frac{1}{2}\mathbb{E}\left( \int_{0}^{\tau^{n}_{N}\wedge t}\mathcal{G}^{2}(s)\left\|(\text{u}_{n},\phi_{n})(s)\right\|^{2}_{\mathbb{V}} ds\right) ^{2}\\
                             &\leq 2\mathbb{E}\left[ \mathcal{E}(\text{u}_{n,0},\phi_{n,0})\right]^{2}
                             +2C\mathbb{E}\int_{0}^{\tau^{n}_{N}\wedge t} \mathcal{G}^{4}(s)B(s)ds
                            +2C\mathbb{E}\int_{0}^{\tau^{n}_{N}\wedge t}\mathcal{G}^{4}(s)\left[ \mathcal{E}(\text{u}_{n},\phi_{n})(s)\right]^{2}ds.
                              \end{aligned}
                              \end{equation}
           Using  Gronwall's inequality, we derive that
                           \begin{equation}\label{n2a}
                           \begin{aligned}
                           \mathbb{E}\sup\limits_{s\in[0,\tau^{n}_{N}\wedge t]}\left[\mathcal{G}^{2}(s)\mathcal{E}(\text{u}_{n},\phi_{n})(s)\right]^{2} &+\frac{1}{2}\mathbb{E}\left( \int_{0}^{\tau^{n}_{N}\wedge t}\mathcal{G}^{2}(s)\left\|(\text{u}_{n},\phi_{n})(s)\right\|^{2}_{\mathbb{V}} ds\right) ^{2}\\
                          &\leq C\mathbb{E}\left[ \mathcal{E}(\text{u}_{n,0},\phi_{n,0})\right]^{2}
                         +C\mathbb{E}\int_{0}^{\tau^{n}_{N}\wedge t} \mathcal{G}^{4}(s)B(s)ds.
                           \end{aligned}
                         \end{equation}
         As in \cite{Ta1} or in \cite{NSC}, using the fact the $\tau^{n}_{N}\uparrow T$ as $N$ goes to $\infty$, we infer that
                           \begin{equation}\label{n2b}
                            \begin{aligned}
                             \mathbb{E}\sup\limits_{s\in[0, t]}\left[\mathcal{G}^{2}(s)\mathcal{E}(\text{u}_{n},\phi_{n})(s)\right]^{2} &+\frac{1}{2}\mathbb{E}\left( \int_{0}^{ t}\mathcal{G}^{2}(s)\left\|(\text{u}_{n},\phi_{n})(s)\right\|^{2}_{\mathbb{V}} ds\right) ^{2}\\
                               &\leq C\mathbb{E}\left[ \mathcal{E}(\text{u}_{n,0},\phi_{n,0})\right]^{2}
                               +C\mathbb{E}\int_{0}^{t} \mathcal{G}^{4}(s)B(s)ds.
                               \end{aligned}
                              \end{equation}
         This completes the proof of Lemma $\ref{La2}$.
                               \end{prev}  
                               
            While exploiting the above lemmas results, we are actually proving Theorem $\ref{Theo2}$.
                           \subsection{Proof of Theorem $\ref{Theo2}$}
                 The proof is splitted into three steps.\\
                              \textbf{Step 1:Convergence related to the projection operator}
                              Let $\Pi_{n}=(\mathcal{P}_{1}^{n},\mathcal{P}_{2}^{n}):\mathbb{V}\rightarrow\mathbb{V}_{n}$ be the orthogonal projection defined by
                              \begin{equation*}
                              \Pi_{n}\varPhi=\sum_{j=1}^{n}\widetilde{\beta}_{j}\widetilde{e}_{j}=\sum_{j=1}^{n}\beta_{j}e_{j}\ \ \mbox{with}\ \ \widetilde{\beta}_{j}=(\varPhi,\widetilde{e}_{j})_{\mathbb{V}} \ \ \mbox{and}\ \  \beta_{j}=(\varPhi,e_{j}),  \ \ \forall \varPhi\in \mathbb{V},
                              \end{equation*}
                              where $\{\widetilde{e}_{j}=\frac{1}{\sqrt{\lambda}_{j}}\}_{j=1}^{\infty}$ is the orthonormal basis of $\mathbb{V}$. Using Parseval's identity, it follows that for any $\varPhi\in \mathbb{V}$ we have (see \cite{EVAN})
                              \begin{equation}\label{p1}
                              \begin{aligned}
                              \left\|\Pi_{n}\varPhi \right\|_{\mathbb{Y}}&\leq\left\|\varPhi \right\|_{\mathbb{Y}},  \ \   \left\|\Pi_{n}\varPhi \right\|_{\mathbb{V}}\leq\left\|\varPhi \right\|_{\mathbb{V}},\ \ 
                              \Pi_{n}\varPhi&\rightarrow \varPhi \ \ \mbox{strongly in } \mathbb{V}. 
                              \end{aligned}
                              \end{equation}
              Considering an arbitrary $z\in L^{s}(\Omega\times(0,T); \mathbb{V})$ for some $s\geq 1$, it follows that
                              \begin{equation}\label{p2}
                              \begin{aligned}
                              \left\|\Pi_{n}z\right\|_{\mathbb{V}}\leq\left\|z \right\|_{\mathbb{V}} \ \ \mbox{and}\ \ \Pi_{n}z(\omega,t) \rightarrow z(\omega,t)  \ \ \mbox{strongly in } \mathbb{V},\\
                              \end{aligned}
                              \end{equation}
             which are valid $\mathbb{P}-$a.s. and $a.e.\ \ t\in (0,T)$. Hence, Lebesgue's dominated convergence theorem implies that for any $z\in L^{s}(\Omega\times(0,T); \mathbb{V})$,  we have 
                              \begin{equation}\label{p3}
                              \begin{aligned}
                              \Pi_{n}z(\omega,t) &\rightarrow z(\omega,t)  \ \ \mbox{strongly in } \ \ L^{s}(\Omega\times(0,T); \mathbb{V}). 
                              \end{aligned}
                              \end{equation}
                              \textbf{Step 2: Passage to the limit in the weak sense.}
                              Let us define
                              \begin{math}
                              \mathcal{H}_{0}(t)=C_{0}\left(\left\|(a,b) \right\|^{2}_{\mathcal{H}_{p}(\Gamma)}+1  \right).
                              \end{math}
                              Since
                              \begin{equation*}
                              \int_{0}^{T} \mathcal{H}_{0}(t)dt\leq C(\omega)<\infty \ \ \mbox{for all }\ \ \omega\in \Omega\setminus A, \ \ \mbox{where}\ \ \mathbb{P}(A)=0,
                              \end{equation*}
           it follows from $\eqref{l2}$ that  there exists a positive constant $K(\omega)$, which depend only on $\Omega\setminus A$ and satisfies
                              \begin{equation}\label{p4}
                            0 <K(\omega)\leq \mathcal{G}(t)=\exp\left( -\int_{0}^{t} \mathcal{H}_{0}(s)ds\right) \leq 1 \ \ \mbox{for all }\ \ \omega\in \Omega\setminus A, \ \ t\in[0,T].
                              \end{equation}
                    From $\eqref{f}$, we deduce that
                                 \begin{equation}\label{p7}
                                 \mathbb{E}\sup\limits_{t\in[0,T]}\mathcal{G}^{4}(t)\left|f_{\vartheta}(\phi_{n}) \right|^{2}_{L^{2}}\leq C, \ \ \forall n\in \mathbb{N}.
                                \end{equation}
               The estimates $\eqref{g5}$, $\eqref{n2b}$ and $\eqref{p7}$ imply that
                              \begin{equation}\label{p5}
                              \begin{aligned}
                              \mathbb{E}\sup\limits_{t\in[0,T]}\mathcal{G}^{2}(t)\left\|(\text{u}_{n},\phi_{n})(t) \right\|^{2}_{\mathbb{Y}}&\leq C, \ \ \mathbb{E}\int_{0}^{T}\ \mathcal{G}^{2}(t) \left\|(\text{u}_{n},\phi_{n})(t) \right\|^{2}_{\mathbb{V}}dt\leq C,\\
                              \mathbb{E}\sup\limits_{t\in[0,T]}\mathcal{G}^{4}(t)\left\|(\text{u}_{n},\phi_{n})(t) \right\|^{4}_{\mathbb{Y}}&\leq C, \ \
                              \mathbb{E}\left( \int_{0}^{T}\mathcal{G}^{2}(t) \left\|(\text{u}_{n},\phi_{n})(t) \right\|^{2}_{\mathbb{V}}dt\right)^{2}\leq C, \\
                               \mathbb{E}\int_{0}^{T}\mathcal{G}^{2}(t)\left| \mu_{n}(s) \right|^{2}_{L^{2}}dt\leq &C, \ \ \mathbb{E}\left(\int_{0}^{T}\mathcal{G}^{2}(t)\left| \mu_{n}(s)\right|^{2}_{L^{2}}dt\right)^{2}\leq C,
                              \end{aligned}
                              \end{equation}
                 for some constant $C$ that is independent of the index $n$. \\
                 Using the H\"{o}lder, Ladyzhenskaya inequalities, Sobolev embedding and  $\eqref{p5}$, we can prove that
                              \begin{equation}\label{p6}
                              \begin{aligned}
                              \left\|\mathcal{G}^{2}(\text{v}_{n}\cdot\nabla)\text{v}_{n} \right\|_{L^{2}(\Omega\times(0,T); \mathbb{V}^{*}_{div})} &\leq C, \ \ \forall n\in \mathbb{N},\\
                              \left\|\mathcal{G}^{2}(\mu_{n}\nabla\phi_{n}) \right\|_{L^{2}(\Omega\times(0,T); \mathbb{V}^{*}_{div})} &\leq C, \ \ \forall n\in \mathbb{N},\\
                              \left\|\mathcal{G}^{2}(\text{v}_{n}\cdot\nabla)\phi_{n} \right\|_{L^{2}(\Omega\times(0,T); \text{V}_{1}^{*})} &\leq C, \ \ \forall n\in \mathbb{N}.
                              \end{aligned}
                              \end{equation}
               By $\eqref{p5}$ and along with the Banach Alaogly theorem, we can abstract a subsequence still denoted by $\{(\text{u}_{n},\phi_{n})\}$ to simplify the notation which converges to the following limits such that
                              \begin{equation}\label{p8}
                              \begin{aligned}
                              \mathcal{G}(\text{u}_{n},\phi_{n})\rightharpoonup \mathcal{G}(\text{u},\phi) \ \ &\mbox{weakly in }\ \ L^{2}(\Omega\times(0,T); \mathbb{V})\cap L^{4}(\Omega, L^{2}(0,T; \mathbb{V})),\\
                              \mathcal{G}(\text{u}_{n},\phi_{n})\rightharpoonup \mathcal{G}(\text{u},\phi) \ \ &\mbox{*-weakly in }\ \ L^{2}(\Omega,L^{\infty}(0,T; \mathbb{Y})\cap  L^{4}(\Omega,L^{\infty}(0,T; \mathbb{Y})).
                              \end{aligned}
                              \end{equation}
                Moreover, we have (see $\eqref{p3}$)
                              \begin{equation}\label{p9}
                              \begin{aligned}
                              \mathcal{G}\Pi_{n}(\text{u},\phi)&\rightarrow \mathcal{G}(\text{u},\phi) \ \ \mbox{strongly in }\ \ L^{2}(\Omega\times(0,T); \mathbb{V})\cap L^{4}(\Omega, L^{2}(0,T; \mathbb{V})).
                              \end{aligned}
                            \end{equation}
           The limit function $(\text{u},\phi)$ satisfies the estimates $\eqref{n2b0}$, $\eqref{g50}$ by the lower semicontinuity of the integral in $L^{4}$ and $L^{2}$ spaces.\\
          Using \textbf{(H1)}, $\eqref{p6}$ and $\eqref{p7}$, we deduce that there exists some operators $\text{B}^{*}_{0}$, $\text{B}^{*}_{1}$, $\text{R}^{*}_{0}$ and $\text{g}^{*}$
                              \begin{equation}\label{p10}
                              \begin{aligned}
                              \mathcal{G}^{2}(\text{v}_{n}\cdot\nabla)\text{v}_{n}\rightharpoonup \mathcal{G}^{2}\text{B}^{*}_{0}(t) \ \ &\mbox{weakly in }\ \ L^{2}(\Omega\times(0,T); \mathbb{V}^{*}_{div}),\\
                              \mathcal{G}^{2}(\mu_{n}\nabla\phi_{n})\rightharpoonup \mathcal{G}^{2}\text{R}^{*}_{0}(t) \ \ &\mbox{weakly in }\ \ L^{2}(\Omega\times(0,T); \mathbb{V}^{*}_{div}),\\
                              \mathcal{G}^{2}(\text{v}_{n}\cdot\nabla)\phi_{n}\rightharpoonup \mathcal{G}^{2}\text{B}^{*}_{1}(t) \ \ &\mbox{weakly in }\ \ L^{2}(\Omega\times(0,T); \text{V}_{1}^{*}),\\
                              \mathcal{G}^{2}\text{g}(t,\text{v}_{n})\rightharpoonup \mathcal{G}^{2}\text{g}^{*}(t) \ \ &\mbox{weakly in }\ \ L^{2}(\Omega\times(0,T); \mathbb{H}_{div}^{m}),
                              \end{aligned}
                              \end{equation}
                              and
                              \begin{equation}\label{p11}
                              \mathcal{G}^{2}f_{\vartheta}(\phi_{n})\rightharpoonup \mathcal{G}^{2}f^{*}(t) \ \ \mbox{weakly in }\ \  L^{2}(\Omega,L^{\infty}(0,T; \text{H})).
                              \end{equation}   
         Now, let us define by
                            \begin{equation*}
                            d\left(\mathcal{G}^{2}(t)\mathcal{Z}_{n}(t),\varphi \right)= d\left(\mathcal{G}^{2}(t)(\text{v}_{n},\phi_{n})(t),\varphi \right)=d\left(\mathcal{G}^{2}(t)\text{v}_{n}(t),\varphi_{1} \right)+d\left(\mathcal{G}^{2}(t)\phi_{n}(t),\varphi_{2} \right).
                            \end{equation*}
            Since $(\text{v}_{n},\phi_{n})$ solves the equation $\eqref{g1}$, applying the It\^{o} formula to $\mathcal{G}^{2}(t)\text{v}_{n}(t)$ and adding with $d(\mathcal{G}^{2}(t)\phi_{n}(t),\varphi_{2}) $  we obtain
                            \begin{equation*}
                            \begin{aligned}
                            d\left(\mathcal{G}^{2}(t)\mathcal{Z}_{n}(t),\varphi \right)&=\mathcal{G}^{2}(t)\left[- ((\text{v}_{n},\varphi_{1}))+\int_{\Gamma}b(\varphi_{1}\cdot\tau)d\gamma-\left\langle (\text{v}_{n}\cdot\nabla)\text{v}_{n},\varphi_{1} \right\rangle  \right]dt\\
                             &+\mathcal{G}^{2}(t)\left[\left\langle \mu_{n}\nabla\phi_{n},\varphi_{1}\right\rangle +\left\langle -( \text{v}_{n}\cdot\nabla)\phi_{n}-\mu_{n},\varphi_{2}\right\rangle \right]dt\\
                             &-2\mathcal{H}_{0}(t)\mathcal{G}^{2}(t)(\mathcal{Z}_{n}(t),\varphi)dt+\mathcal{G}^{2}(t)\left( \text{g}(t,\text{v}_{n}),\varphi_{1}\right) dW_{t},
                            \end{aligned}
                            \end{equation*}
                            which implies
                           \begin{equation}\label{p12}
                            \begin{aligned}
                           &\left(\mathcal{G}^{2}(t)\mathcal{Z}_{n}(t),\varphi \right)-\left(\mathcal{Z}_{n,0},\varphi \right)\\
                           &=\int_{0}^{t}\mathcal{G}^{2}(s)\left[- ((\text{v}_{n}(s),\varphi_{1}))+\int_{\Gamma}b(s)(\varphi_{1}\cdot\tau)d\gamma-\left((\text{v}_{n}(s)\cdot\nabla)\text{v}_{n}(s),\varphi_{1} \right) \right]ds\\
                           &+\int_{0}^{t}\mathcal{G}^{2}(s)\left[ \left\langle \mu_{n}(s)\nabla\phi_{n}(s),\varphi_{1}\right\rangle -\left\langle ( \text{v}_{n}(s)\cdot\nabla)\phi_{n}(s),\varphi_{2}\right\rangle -(\mu_{n}(s),\varphi_{2})\right] ds\\
                          &-2\int_{0}^{t}\mathcal{H}_{0}(s)\mathcal{G}^{2}(s)(\mathcal{Z}_{n}(s),\varphi)ds+\int_{0}^{t}\mathcal{G}^{2}(s)\left( \text{g}(s,\text{v}_{n}(s)),\varphi_{1}\right) dW_{s},
                          \end{aligned}
                        \end{equation}
                        for all $t\in[0,T]$, $\mathbb{P}-$a.s.  Let us note by:
                        \begin{equation*}
                        \text{h}_{n}(t)=\mathcal{G}^{2}(t)\mathcal{Z}_{n}(t)-\int_{0}^{t}\mathcal{G}^{2}(s) \text{g}(s,\text{v}_{n}(s)) dW_{s}
                        \end{equation*}
                        the following differential equation holds 
                        \begin{equation}\label{p13}
                        \begin{aligned}
                        \frac{\partial}{\partial t}(\text{h}_{n}(t),\varphi)&=\mathcal{G}^{2}(t)\left[- ((\text{v}_{n}(t),\varphi_{1}))+\int_{\Gamma}b(\varphi_{1}\cdot\tau)d\gamma-\left\langle (\text{v}_{n}(t)\cdot\nabla)\text{v}_{n}(t),\varphi_{1} \right\rangle  \right]dt\\
                         &+\mathcal{G}^{2}(t)\left[\left\langle \mu_{n}(t)\nabla\phi_{n}(t),\varphi_{1}\right\rangle +\left\langle -( \text{v}_{n}(t)\cdot\nabla)\phi_{n}(t)-\mu_{n}(t),\varphi_{2}\right\rangle \right]dt\\
                         &-2\mathcal{H}_{0}(t)\mathcal{G}^{2}(t)(\mathcal{Z}_{n}(t),\varphi)dt
                        \end{aligned}
                        \end{equation}
                        for all $t\in[0,T]$, $\mathbb{P}-$a.e.  in $\Omega.$ \\
                        Due to the properties of the stochastic integral and the assumption \textbf{(H1)}, we deduce that
                        \begin{equation*}
                        \text{h}_{n}(t)\rightharpoonup \text{h}(t)=\mathcal{G}^{2}(t)\mathcal{Z}(t)-\int_{0}^{t}\mathcal{G}^{2}(s) \text{g}^{*}(s) dW_{s}\ \ \mbox{weakly in }\ \ L^{2}(\Omega\times(0,T); \mathbb{H}^{1}(\text{D})\times\text{V}_{2}).
                        \end{equation*}
                        We passe to the limit in the equation $\eqref{p12}$ in the distributed sense. Namely multiplying the equation $\eqref{p13}$ by the test function $\theta(t)\chi(t)$, with $\theta\in\mathcal{C}^{\infty}([0,T])$ with compact support and $\chi\in \mathbb{G}(\Omega)=\mathbb{L}^{2}(\Omega)\times\text{H}^{1}(\Omega)$, and passing to limit, we obtain
                         \begin{equation*}
                         \begin{aligned}
                         &\mathbb{E}\int_{0}^{T}(\text{h}(t),\varphi)\theta^{'}(t)\chi(t)dt\\
                         &=\mathbb{E}\int_{0}^{T}\mathcal{G}^{2}(t)\left[- ((\text{v}(t),\varphi_{1}))+\int_{\Gamma}b(\varphi_{1}\cdot\tau)d\gamma-\left\langle (\text{v}(t)\cdot\nabla)\text{v}(t),\varphi_{1} \right\rangle  \right]\theta(t)\chi(t)dt\\
                          &+\mathbb{E}\int_{0}^{T}\mathcal{G}^{2}(t)\left[\left\langle \mu(t)\nabla\phi(t),\varphi_{1}\right\rangle +\left\langle -( \text{v}(t)\cdot\nabla)\phi(t)-\mu(t),\varphi_{2}\right\rangle \right]\theta(t)\chi(t)dt\\
                          &-2\int_{0}^{T}\mathcal{H}_{0}(t)\mathcal{G}^{2}(t)(\mathcal{Z}(t),\varphi)\theta(t)\chi(t)dt.
                          \end{aligned}
                         \end{equation*}
                         Consequently, we have $\frac{\partial \text{h}}{\partial t}\in L^{2}(\Omega\times(0,T); \mathbb{H}^{-1}(\text{D})\times\text{V}^{*}_{2})$. Since $L^{2}(\Omega\times(0,T); \mathbb{H}^{1}(\text{D})\times\text{V}_{2})$, then, we conclude that $h\in L^{2}(\Omega, \mathcal{C}([0,T]); \mathbb{L}^{2}(\text{D})\times\text{V}_{1})$ by the Aubin-Lions embedding result. Taking into account the continuity property of the stochastic integral, we conclude that $\mathcal{G}^{2}(\text{v},\phi)\in  L^{2}(\Omega, \mathcal{C}([0,T]); \mathbb{L}^{2}(\text{D})\times\text{V}_{1})$. In addition
                         \begin{equation*}
                         \mathcal{G}^{2}(\text{v}_{n},\phi_{n})\rightharpoonup\mathcal{G}^{2}(\text{v},\phi)\ \ \mbox{in}\ \ \mathcal{C}_{w}([0,T]; \mathbb{G}(\Omega)\times\mathbb{U}),
                         \end{equation*}
                         where $\mathbb{G}(\Omega)=\mathbb{L}^{2}(\Omega)\times\text{H}^{1}(\Omega)$, $\mathbb{U}=\mathbb{L}^{2}(\text{D})\times\text{H}^{1}(\text{D})$ and the index $w$ means that we are considering $\mathbb{G}(\Omega)\times\mathbb{U}$ endowed with the weak topology. Thus, we have
                         \begin{equation}\label{b0}
                          \mathbb{E}\left[ \left(\mathcal{G}^{2}(t)(\text{v}_{n},\phi_{n})(t),\varphi \right)\chi \right] \rightarrow\mathbb{E}\left[ \left(\mathcal{G}^{2}(t)(\text{v},\phi)(t),\varphi \right)\chi \right], \ \ \forall t\in[0,T].
                         \end{equation}
                         We multiply the equation $\eqref{p12}$ by an arbitrary $\chi\in \mathbb{G}(\Omega)=\mathbb{L}^{2}(\Omega)\times\text{H}^{1}(\Omega)$, and take the expectation, we obtain
                         \begin{equation}\label{p12a}
                           \begin{aligned}
                            &\mathbb{E}\chi\{\left(\mathcal{G}^{2}(t)\mathcal{Z}_{n}(t),\varphi \right)-\left(\mathcal{Z}_{n,0},\varphi \right)\}\\
                            &=\mathbb{E}\chi\{\int_{0}^{t}\mathcal{G}^{2}(s)\left[- ((\text{v}_{n}(s),\varphi_{1}))+\int_{\Gamma}b(\varphi_{1}\cdot\tau)d\gamma-\left\langle (\text{v}_{n}(s)\cdot\nabla)\text{v}_{n}(s),\varphi_{1} \right\rangle  \right]ds\}\\
                             &+\mathbb{E}\chi\{\int_{0}^{t}\mathcal{G}^{2}(s)\left[ \left\langle \mu_{n}(s)\nabla\phi_{n}(s),\varphi_{1}\right\rangle -\left\langle ( \text{v}_{n}(s)\cdot\nabla)\phi_{n}(s),\varphi_{2}\right\rangle -(\mu_{n}(s),\varphi_{2})\right] ds\}\\
                             &\mathbb{E}\chi\{-2\int_{0}^{t}\mathcal{H}_{0}(s)\mathcal{G}^{2}(s)(\mathcal{Z}_{n}(s),\varphi)ds+\int_{0}^{t}\mathcal{G}^{2}(s)\left( \text{g}(s,\text{v}_{n}(s)),\varphi_{1}\right) dW_{s}\}.
                            \end{aligned}
                            \end{equation}
                           Using $\eqref{p8}-\eqref{p10}$ and $\eqref{b0}$, we pass to the limit $n\rightarrow\infty$ in this equality and we deduce that
                           \begin{equation}\label{p12b}
                            \begin{aligned}
                             &\mathbb{E}\chi\{\left(\mathcal{G}^{2}(t)\mathcal{Z}(t),\varphi \right)-\left(\mathcal{Z}_{0},\varphi \right)\}\\
                             &=\mathbb{E}\chi\{\int_{0}^{t}\mathcal{G}^{2}(s)\left[- ((\text{v}(s),\varphi_{1}))+\int_{\Gamma}b(\varphi_{1}\cdot\tau)d\gamma-\left(\text{B}^{*}_{0}(s),\varphi_{1} \right) \right]ds\}\\
                             &+\mathbb{E}\chi\{\int_{0}^{t}\mathcal{G}^{2}(\text{R}^{*}_{0}(s),\varphi_{1})ds-\int_{0}^{t}( \text{B}^{*}_{1},\varphi_{2})ds-\int_{0}^{t}(\mu(s),\varphi_{2})ds\}\\
                             &\mathbb{E}\chi\{-2\int_{0}^{t}\mathcal{H}_{0}(s)\mathcal{G}^{2}(s)(\mathcal{Z}(s),\varphi)ds+\int_{0}^{t}\mathcal{G}^{2}(s)\left( \text{g}^{*},\varphi_{1}\right) dW_{s}\}.
                             \end{aligned}
                             \end{equation}
                             Since  $\chi\in \mathbb{G}(\Omega)=\mathbb{L}^{2}(\Omega)\times\text{H}^{1}(\Omega)$ is arbitrary, it follows that
                             \begin{equation}\label{p12c}
                              \begin{aligned}
                                &\left(\mathcal{G}^{2}(t)\mathcal{Z}(t),\varphi \right)-\left(\mathcal{Z}_{0},\varphi \right)\\
                                 &=\int_{0}^{t}\mathcal{G}^{2}(s)\left[- ((\text{v}(s),\varphi_{1}))+\int_{\Gamma}b(\varphi_{1}\cdot\tau)d\gamma-\left(\text{B}^{*}_{0}(s),\varphi_{1} \right) \right]ds\\
                                 &+\int_{0}^{t}\mathcal{G}^{2}(\text{R}^{*}_{0}(s),\varphi_{1})ds-\int_{0}^{t}( \text{B}^{*}_{1}(s),\varphi_{2})ds-\int_{0}^{t}(\mu(s),\varphi_{2})ds\\
                                 &-2\int_{0}^{t}\mathcal{H}_{0}(s)\mathcal{G}^{2}(s)(\mathcal{Z}(s),\varphi)ds+\int_{0}^{t}\mathcal{G}^{2}(s)\left( \text{g}^{*}(s),\varphi_{1}\right) dW_{s}.
                                 \end{aligned}
                               \end{equation}
                             for any $t\in[0,T]$ and $\mathbb{P}-$a.e. in $\Omega,$ 
                              that is
                             \begin{equation*}
                              \begin{aligned}
                              d\left(\mathcal{G}^{2}\mathcal{Z},\varphi \right)&=\mathcal{G}^{2}\left[- ((\text{v},\varphi_{1}))+\int_{\Gamma}b(\varphi_{1}\cdot\tau)d\gamma-\left(\text{B}^{*}_{0}(t),\varphi_{1} \right) \right]dt\\
                               &+\mathcal{G}^{2}\left[(\text{R}^{*}_{0}(t),\varphi_{1})+\left(-\text{B}^{*}_{1}(t)+\Delta\mu,\varphi_{2}\right)\right]dt\\
                                &-2\mathcal{H}_{0}(t)\mathcal{G}^{2}(t)(\mathcal{Z},\varphi)dt+\mathcal{G}^{2}\left( \text{g}^{*}(t),\varphi_{1}\right) dW_{t}, \ \ \mathcal{Z}_{0}=(\text{v}_{0},\phi_{0}).
                                \end{aligned}
                                \end{equation*}
                                If we apply the It\^{o} formula, we obtain
                                \begin{equation}
                               d\left(\mathcal{Z},\varphi \right)=d[\mathcal{G}^{-2}\mathcal{G}^{2} (\mathcal{Z},\varphi)]=\mathcal{G}^{2}(\mathcal{Z},\varphi)d(\mathcal{G}^{-2})+\mathcal{G}^{-2}d[\left(\mathcal{G}^{2}\mathcal{Z},\varphi \right)],
                                \end{equation}
                                we derive that the limit function $\mathcal{Z}=(\text{v},\phi)$ in the form $\mathcal{Z}=(\text{u}+\textbf{a},\phi)$ with
                                \begin{equation*}
                                (\text{u},\phi)\in L^{\infty}(0,T;\mathbb{Y})\cap L^{2}(0,T; \mathbb{V}), \ \ \mathbb{P}-\mbox{a.e. in } \Omega,  \mbox{a.e. on} \ \ (0,T),
                                \end{equation*}
                                satisfies $\mathbb{P}-\mbox{a.e. in } \Omega$, the stochastic differential equation
                                \begin{equation}\label{p37}
                               \begin{aligned}
                               d\left(\mathcal{Z}(t),\varphi \right)&=\left[- ((\text{v},\varphi_{1}))+\int_{\Gamma}b(\varphi_{1}\cdot\tau)d\gamma-\left(\text{B}^{*}_{0}(t),\varphi_{1} \right) \right]dt\\
                               &+\left[(\text{R}^{*}_{0}(t),\varphi_{1})+\left(-\text{B}^{*}_{1}(t)-\mu,\varphi_{2}\right)\right]dt\\
                               &+\left( \text{g}^{*}(t),\varphi_{1}\right) dW_{t}, \ \ \forall t\in[0,T], \ \ \varphi\in\mathbb{V},\ \ \mathcal{Z}_{0}=(\text{y}_{0},\phi_{0}).
                               \end{aligned}
                              \end{equation}                                                               
          \textbf{Step 3: Deduction of strong convergence as $n\rightarrow\infty$.}
                                Exploiting the methods  used in \cite{Ta1}, \cite{H.B2}, we prove that the limit process $(\text{v},\phi)$ satisfies the equation $\eqref{Ton2}$. Now, we write $\text{v}=\text{u}+\textbf{a}$, $\text{v}_{n}=\text{u}_{n}+\textbf{a}$ and taking the difference of the equations $\eqref{g2}$ and $\eqref{p37}$ with $(\varphi,\psi)=(e_{i},\psi)\in\mathbb{V}_{n}, i=1,...,n$, we obtain the following system
                                \begin{equation}\label{g2a}
                                 \begin{cases}
                                 d(\widetilde{\text{u}}_{n}-\text{u}_{n},e_{i})=\left[-((\widetilde{\text{u}}_{n}-\text{u}_{n},e_{i}))+((\text{v}_{n}\cdot\nabla)\text{v}_{n}-\text{B}^{*}_{0}(t),e_{i}) \right]dt \\
                                 -(\mu_{n}\nabla\phi_{n}-\text{R}^{*}_{0}(t),e_{i})dt-(g(t,\text{v}_{n})-\text{g}^{*}(t), e_{i})dW(s),\\
                                 d(\widetilde{\phi}_{n}-\phi_{n},\psi)=\left((\widetilde{\mu}_{n}-\mu_{n})+(\text{v}_{n}\cdot\nabla\phi_{n}-\text{B}^{*}_{1}(t)), \psi \right)dt,\\
                                 (\widetilde{\mu}_{n}-\mu_{n})=\text{A}_{\vartheta}(\widetilde{\phi}_{n}-\phi_{n})+f_{\vartheta}(\phi_{n})-f^{*}(t),
                                \end{cases}
                                \end{equation}
                              where $(\widetilde{\text{u}}_{n},\widetilde{\phi}_{n},\widetilde{\mu}_{n} )=\Pi_{n}(\text{u},\phi, \mu)$  is the orthogonal projection of $\mathbb{Y}$ onto $\mathbb{Y}_{n}$. Let us set
                              \begin{math}
                              \theta_{n}=\widetilde{\text{u}}_{n}-\text{u}_{n}, \ \ \Psi_{n}=(\widetilde{\phi}_{n}-\phi_{n}) \ \ \mbox{and}\ \ \xi_{n}=(\widetilde{\mu}_{n}-\mu_{n}).\\
                              \end{math}
                              Applying the  It\^{o} formula, we derive that
                                \begin{equation}\label{k0}
                                  \begin{aligned}
                                   d\left( (\theta_{n},e_{i})^{2}\right) &=2(\theta_{n},e_{i})\left[-((\theta_{n}, e_{i}))+(\text{v}_{n}\cdot\nabla)(\text{v}_{n}-\text{B}^{*}_{0}(t)),e_{i}) \right]dt \\
                                    &-2(\theta_{n},e_{i})(\mu_{n}\nabla\phi_{n}-\text{R}^{*}_{0}(t),e_{i})dt-2(\theta_{n},e_{i})(g(t,\text{v}_{n})-\text{g}^{*}(t), e_{i})dW(s)\\
                                    &+ \left|(g(t,\text{v}_{n})-\text{g}^{*}(t), e_{i}) \right|^{2}_{L^{2}}dt.
                                      \end{aligned}
                                        \end{equation}
                                        Summing over i=1,...,n, we obtain
                               \begin{equation}\label{k1}
                                    \begin{aligned}
                                   d\left| \theta_{n}(t)\right|^{2}_{L^{2}} &+2\left\| \theta_{n}\right\|^{2}dt=2((\text{v}_{n}\cdot\nabla)(\text{v}_{n}-\text{B}^{*}_{0}(t)),\theta_{n}) dt \\
                                   &-2(\mu_{n}\nabla\phi_{n}-\text{R}^{*}_{0}(t),\theta_{n})dt-2(g(t,\text{v}_{n})-\text{g}^{*}(t), \theta_{n})dW(t)\\
                                   &+ \sum_{i=1}^{n}\left|(g(t,\text{v}_{n})-\text{g}^{*}(t), e_{i}) \right|^{2}_{L^{2}}dt.
                                   \end{aligned}
                                  \end{equation}
             Replacing $\psi$ in $\eqref{g2a}_{2}$ by $\text{A}_{\vartheta}\Psi_{n}$ and multiplying the equation $\eqref{g2a}_{3}$ with $\xi_{n}-\text{A}_{\vartheta}\Psi_{n}$ and adding these equations, we obtain
                                  \begin{equation}\label{k3}
                                  \begin{aligned}
                                  d\left|\nabla\Psi_{n}(t)\right|^{2}_{L^{2}} &+2\left\|\Psi_{n}\right\|^{2}_{2}dt+2\left|\xi_{n} \right|^{2}_{L^{2}}dt = 2((\text{v}_{n}\cdot\nabla\phi_{n}-\text{B}^{*}_{1}(t)), \text{A}_{\vartheta}\Psi_{n})dt\\
                                  &-\left(f_{\vartheta}(\phi_{n})-f^{*}, \text{A}_{\vartheta}\Psi_{n}\right)dt+\left(f_{\vartheta}(\phi_{n})-f^{*}, \xi_{n}\right)dt. 
                                  \end{aligned}
                                  \end{equation}
          It follows from $\eqref{k1}-\eqref{k3}$ that
                                   \begin{equation}\label{k1a}
                                     \begin{aligned}
                                    &d\left( \left| \theta_{n}(t)\right|^{2}_{L^{2}}+ \left|\nabla\Psi_{n}(t)\right|^{2}_{L^{2}}\right) +2\left( \left\| \theta_{n}\right\|^{2}+\left\|\Psi_{n}\right\|^{2}_{2}+\left|\xi_{n} \right|^{2}_{L^{2}}\right) dt=2((\text{v}_{n}\cdot\nabla)\text{v}_{n}-\text{B}^{*}_{0}(t)),\theta_{n}) dt \\
                                    &-2(\mu_{n}\nabla\phi_{n}-\text{R}^{*}_{0}(t),\theta_{n})dt
                                    + 2((\text{v}_{n}\cdot\nabla\phi_{n}-\text{B}^{*}_{1}(t)), \text{A}_{\vartheta}\Psi_{n})dt\\
                                    &+\left(f_{\vartheta}(\phi_{n})-f^{*}, \text{A}_{\vartheta}\Psi_{n}\right)dt+ \left(f_{\vartheta}(\phi_{n})-f^{*}, \xi_{n}\right)dt\\
                                    &-2(g(t,\text{v}_{n})-\text{g}^{*}(t), \theta_{n})dW(t)
                                    + \sum_{i=1}^{n}\left|(g(t,\text{v}_{n})-\text{g}^{*}(t), e_{i}) \right|^{2}_{L^{2}}dt.
                                      \end{aligned}
                                   \end{equation}
                                   Now, we observe that  
                      \begin{equation}\label{h1}
                                              \begin{aligned}
                                              &((\text{v}_{n}\cdot\nabla)(\text{v}_{n}-\text{B}^{*}_{0}(t)),\theta_{n})\\
                                              &=-(((\text{u}_{n}+\textbf{a})\cdot\nabla\theta_{n}),\theta_{n})-(\theta_{n}\cdot\nabla(\widetilde{\text{u}}_{n}+\textbf{a}),\theta_{n})\\
                                              &+((\widetilde{\text{u}}_{n}-\text{u})\cdot\nabla(\widetilde{\text{u}}_{n}+\textbf{a}),\theta_{n})+((\text{u}+\textbf{a})\nabla(\widetilde{\text{u}}_{n}-\text{u}),\theta_{n})\\
                                              &+((\text{v}\cdot\nabla)\text{v}-\text{B}^{*}_{0}(t),\theta_{n})=(I_{0}+I_{1},  \theta_{n} ) +I_{2}+I_{3}+I_{4},
                                              \end{aligned}
                                              \end{equation}
                                              \begin{equation}\label{k5}
                                              \begin{aligned}
                                              &(\mu_{n}\nabla\phi_{n}-\text{R}^{*}_{0}(t),\theta_{n})
                                              =-(\mu_{n}\nabla\Psi_{n},\theta_{n})-(\xi_{n}\nabla\widetilde{\phi}_{n},\theta_{n})+((\widetilde{\mu}_{n}-\mu)\nabla\widetilde{\phi}_{n},\theta_{n})\\
                                              &+(\mu\nabla(\widetilde{\phi}_{n}-\phi))
                                              +(\mu\nabla\phi-\text{R}^{*}_{0}(t),\theta_{n})=I_{5}+I_{6}+I_{7}+I_{8}+I_{9},
                                              \end{aligned}
                                              \end{equation}
                                              \begin{equation}\label{k5a}
                                               \begin{aligned}
                                               &((\text{v}_{n}\cdot\nabla\phi_{n}-\text{B}^{*}_{1}(t)), \Psi_{n}-\Delta\Psi_{n})
                                               =-(((\text{u}_{n}+\textbf{a})\cdot\nabla\Psi_{n}),\Psi_{n}-\Delta\Psi_{n})\\
                                               &-(\theta_{n}\cdot\nabla\widetilde{\phi}_{n},\Psi_{n}-\Delta\Psi_{n})
                                               +((\widetilde{\text{u}}_{n}-\text{u})\cdot\nabla\widetilde{\phi}_{n},\Psi_{n}-\Delta\Psi_{n})\\
                                               &+((\text{u}+\textbf{a})\nabla(\widetilde{\phi}_{n}-\phi),\Psi_{n}-\Delta\Psi_{n})+((\text{v}\cdot\nabla\phi-\text{B}^{*}_{1}(t)), \Psi_{n}-\Delta\Psi_{n})\\
                                               &=I_{10}+I_{11}+I_{12}+I_{13}+I_{14},
                                               \end{aligned}
                                               \end{equation}
                                               and
                                               \begin{equation}\label{h2}
                                               f_{\vartheta}(\phi_{n})-f^{*}(t)=-(f_{\vartheta}(\widetilde{\phi}_{n})-f_{\vartheta}(\phi_{n}))+(f_{\vartheta}(\widetilde{\phi}_{n})-f_{\vartheta}(\phi))+(f_{\vartheta}(\phi)-f^{*}(t)).
                                               \end{equation}
                       Thanks to the H\"{o}lder inequality, Ladyzhenskaya, Young, Gagliardo-Nirenberg inequalities and Sobolev embedding, we estimate each term of the right-hand side of the above equality individually.    
                   \begin{equation}\label{e1}
                                 \begin{aligned}
                                & \left| (I_{0}+I_{1},  \theta_{n} )\right|\leq \left|\int_{\Gamma}a \left(\theta_{n}.\tau \right)^{2}d\gamma  \right|+\left|(\theta_{n}\cdot\nabla(\widetilde{\text{u}}_{n}+\textbf{a}),\theta_{n}) \right|\\
                                  &\leq \left|a \right|_{L^{\infty}(\Gamma)} \left|\theta_{n} \right|^{2}_{L^{2}(\Gamma)} +\left\|\widetilde{\text{u}}_{n}+\textbf{a} \right\| \left|\theta_{n} \right|^{2}_{\mathbb{L}^{4}} \\
                                  &\leq \left(\left|a \right|_{L^{\infty}(\Gamma)}+\left\|\textbf{a} \right\|_{\mathbb{H}^{1}}+\left\|\widetilde{\text{u}}_{n} \right\| \right) \left|\theta_{n} \right|_{L^{2}}\left\|\theta_{n} \right\|\\
                                  &\leq C_{2}\left(\left\|(a,b) \right\|^{2} _{\mathcal{H}_{p}(\Gamma)}+\left\|\widetilde{\text{u}}_{n} \right\|^{2} \right)\left|\theta_{n} \right|^{2}_{L^{2}} +\frac{1}{4}\left\|\theta_{n} \right\|^{2}.
                                 \end{aligned}
                                 \end{equation}
                                  As in \cite{NSC}, we derive that
                                  \begin{equation}\label{e2}
                                  \left|  I_{2}\right|\leq C\left\|\widetilde{\text{u}}_{n}-\text{u} \right\|_{\mathbb{L}^{4}} \left(\left\|\text{u} \right\| +\left\|\textbf{a} \right\|_{\mathbb{H}^{1}}  \right)\left(\left\|\widetilde{\text{u}}_{n} \right\|_{\mathbb{L}^{4}} +\left\|\text{u}_{n} \right\|_{\mathbb{L}^{4}} \right),
                                  \end{equation}
                                  and
                                   \begin{equation}\label{e3}
                                   \left|  I_{3}\right|\leq C\left\|\widetilde{\text{u}}_{n}-\text{u} \right\| \left(\left\|\text{u} \right\|_{\mathbb{L}^{4}} +\left\|\textbf{a} \right\|_{\mathbb{L}^{4}}  \right)\left(\left\|\widetilde{\text{u}}_{n} \right\|_{\mathbb{L}^{4}} +\left\|\text{u}_{n} \right\|_{\mathbb{L}^{4}} \right).
                                  \end{equation}
                                  Setting   
                   \begin{equation}\label{e4}
                                   \text{g}_{n}=\text{g}(t,\text{v}_{n}),\ \ \text{g}=\text{g}(t,\text{v}), \ \ \text{g}^{*}=\text{g}^{*}(t),
                                   \end{equation}
                    and we have
                                   \begin{equation*}
                                   \sum_{i=1}^{n}\left|(g(t,\text{v}_{n})-\text{g}^{*}(t), e_{i}) \right|^{2}_{L^{2}}=\sum_{i=1}^{n}\left|(\text{g}_{n}-\text{g}^{*}, e_{i}) \right|^{2}_{L^{2}}=\left|\mathcal{
                                   P}_{1}^{n}\text{g}_{n}-\mathcal{P}_{1}^{n}\text{g}^{*} \right|^{2}_{L^{2}}.
                                   \end{equation*}
                 Using the standard relation $x^{2}=(x-y)^{2}-y^{2}+2xy$, we see that
                                   \begin{equation*}
                                   \begin{aligned}
                                   \left|\mathcal{P}_{1}^{n}\text{g}_{n}-\mathcal{P}_{1}^{n}\text{g}^{*} \right|^{2}_{L^{2}}&=\left|\mathcal{P}_{1}^{n}\text{g}_{n}-\mathcal{P}_{1}^{n}\text{g} \right|^{2}_{L^{2}}-\left|\mathcal{P}_{1}^{n}\text{g}-\mathcal{P}_{1}^{n}\text{g}^{*} \right|^{2}_{L^{2}}\\
                                   &-2(\mathcal{P}_{1}^{n}\text{g}_{n}-\mathcal{P}_{1}^{n}\text{g}^{*},\mathcal{P}_{1}^{n}\text{g}-\mathcal{P}_{1}^{n}\text{g}^{*}).
                                   \end{aligned}
                                   \end{equation*}
                  From  assumption \textbf{(H1)} and $\eqref{p2}$, we obtain
                                   \begin{equation*}
                                   \left|\mathcal{P}_{1}^{n}\text{g}_{n}-\mathcal{P}_{1}^{n}\text{g} \right|^{2}_{L^{2}}\leq \left|\text{g}_{n}-\text{g} \right|^{2}_{L^{2}}\leq K\left|\text{u}_{n}-\text{u} \right|^{2}_{L^{2}},
                                   \end{equation*}
             then for the fixed constant $C_{3}=2K$ we deduce that
                                   \begin{equation}\label{e5}
                                    \begin{aligned}
                                   &\sum_{i=1}^{n}\left|(g(t,\text{v}_{n})-\text{g}^{*}(t), e_{i}) \right|^{2}=\left|\mathcal{P}_{1}^{n}\text{g}_{n}-\mathcal{P}_{1}^{n}\text{g}^{*} \right|^{2}_{L^{2}}\\
                                   &\leq K\left|\text{u}_{n}-\text{u} \right|^{2}_{L^{2}}-\left|\mathcal{P}_{1}^{n}\text{g}-\mathcal{P}_{1}^{n}\text{g}^{*} \right|^{2}_{L^{2}}
                                   +2(\mathcal{P}_{1}^{n}\text{g}_{n}-\mathcal{P}_{1}^{n}\text{g}^{*},\mathcal{P}_{1}^{n}\text{g}-\mathcal{P}_{1}^{n}\text{g}^{*})\\
                                   &\leq C_{3}\left|\theta_{n} \right|^{2}_{L^{2}}+C\left|\widetilde{\text{u}}_{n}-\text{u} \right|^{2}_{L^{2}}- \left|\mathcal{P}_{1}^{n}\text{g}-\mathcal{P}_{1}^{n}\text{g}^{*} \right|^{2}_{L^{2}}
                                   +2(\mathcal{P}_{1}^{n}\text{g}_{n}-\mathcal{P}_{1}^{n}\text{g}^{*},\mathcal{P}_{1}^{n}\text{g}-\mathcal{P}_{1}^{n}\text{g}^{*}).
                                   \end{aligned}
                                \end{equation}
              The positive constants $C_{2}$ and $C_{3}$ in $\eqref{e5}$ and $\eqref{e1}$ are independent of $n$.\\
                We have
                                \begin{equation}\label{e6}
                                \begin{aligned}
                                \left| I_{5}\right| &\leq\left\|\theta_{n} \right\|_{\mathbb{L}^{4}} \left|\mu_{n} \right|_{L^{2}} \left\|\nabla\Psi_{n} \right\|_{L^{4}} \leq C\left\|\theta_{n} \right\|\left|\mu_{n} \right|_{L^{2}}\left\|\nabla\Psi_{n} \right\|_{L^{4}} \\
                                & \leq \frac{1}{4}\left\|\theta_{n} \right\|^{2}+C\left|\mu_{n} \right|^{2}_{L^{2}}\left\|\nabla\Psi_{n} \right\|_{L^{4}} ^{2}\\
                                &\leq \frac{1}{4}\left\|\theta_{n} \right\|^{2}+C\left|\mu_{n} \right|^{2}_{L^{2}}\left(\left\|\Psi_{n} \right\|_{2} \left|\nabla\Psi_{n} \right|_{L^{2}} +\left| \nabla\Psi_{n}\right|^{2}_{L^{2}}  \right) \\
                                &\leq \frac{1}{4}\left\|\theta_{n} \right\|^{2}+\frac{1}{4}\left\| \Psi_{n}\right\|^{2}_{2}+C\left|\mu_{n} \right|^{2}_{L^{2}}\left| \nabla\Psi_{n}\right|^{2}_{L^{2}} +C\left|\mu_{n} \right|^{4}_{L^{2}}\left| \nabla\Psi_{n}\right|^{2}_{L^{2}}. 
                                \end{aligned}
                                \end{equation}
                                Using a similarly argument as in $\eqref{e6}$, we arrive at
                               \begin{equation}\label{e7}
                              \begin{aligned}
                              \left| I_{6}\right| &\leq\left\|\theta_{n} \right\|_{\mathbb{L}^{4}} \left|\xi_{n} \right|_{L^{2}} \left\|\nabla\widetilde{\phi}_{n} \right\|_{L^{4}}\leq C\left\|\theta_{n} \right\|^{1/2}\left|\theta_{n} \right|^{1/2}\left|\xi_{n} \right|_{L^{2}}\left\|\nabla\widetilde{\phi}_{n} \right\|_{L^{4}}\\
                              &\leq \frac{1}{3}\left|\xi_{n} \right|^{2}_{L^{2}}+\frac{1}{4}\left\|\theta_{n} \right\|^{2}+C\left(\left\|\widetilde{\phi}_{n} \right\|_{2}^{2} \left|\nabla\widetilde{\phi}_{n} \right|^{2}_{L^{2}} +\left| \nabla\widetilde{\phi}_{n}\right|^{4}_{L^{2}}  \right)\left|\theta_{n} \right|^{2}_{L^{2}}, 
                               \end{aligned}
                               \end{equation}
                               \begin{equation}\label{e8}
                                \begin{aligned}
                                \left| I_{7}\right| &\leq \left\|\theta_{n} \right\|_{\mathbb{L}^{4}} \left|(\widetilde{\mu}_{n}-\mu) \right|_{L^{2}} \left\|\nabla\widetilde{\phi}_{n} \right\|_{L^{4}}
                                \leq C\left(\left\|\text{u}_{n} \right\|_{\mathbb{L}^{4}}+\left\|\widetilde{\text{u}}_{n} \right\|_{\mathbb{L}^{4}}\right) \left|(\widetilde{\mu}_{n}-\mu) \right|_{L^{2}} \left\|\nabla\widetilde{\phi}_{n} \right\|_{L^{4}},
                               \end{aligned}
                                \end{equation}
                                \begin{equation}\label{e9a}
                                  \begin{aligned}
                                  \left| I_{8}\right| &\leq \left\|\theta_{n} \right\|_{\mathbb{L}^{4}} \left\|\nabla(\widetilde{\phi}_{n}-\phi) \right\|_{L^{4}} \left|\mu\right|_{L^{2}}\leq C\left(\left\|\text{u}_{n} \right\|_{\mathbb{L}^{4}}+\left\|\widetilde{\text{u}}_{n} \right\|_{\mathbb{L}^{4}}\right)\left|\mu \right|_{L^{2}}  \left\|\nabla(\widetilde{\phi}_{n}-\phi) \right\|_{L^{4}},
                                   \end{aligned}
                                 \end{equation}
                               \begin{equation}\label{e9}
                              \begin{aligned}
                             \left|I_{10}\right| &\leq \left\|\text{u}_{n}+\textbf{a} \right\|_{\mathbb{L}^{4}} \left\| \nabla\Psi_{n}\right\|_{L^{4}} \left\| \Psi_{n}\right\|_{2}
                             \leq \frac{1}{4} \left\| \Psi_{n}\right\|^{2}_{2}+C\left\|\text{u}_{n}+\textbf{a} \right\|^{4}_{\mathbb{L}^{4}} \left| \nabla\Psi_{n}\right|^{2}_{L^{2}} \\
                             &\leq \frac{1}{4} \left\| \Psi_{n}\right\|^{2}_{2}+C_{4}\left(\left\|(a,b) \right\|^{4} _{\mathcal{H}_{p}(\Gamma)}+\left\|\text{u}_{n} \right\|^{4} \right)\left|\nabla \Psi_{n}\right|^{2}_{L^{2}},
                             \end{aligned}
                             \end{equation}
                             \begin{equation}\label{e10}
                             \begin{aligned}
                             \left|I_{11}\right| &\leq \frac{1}{4} \left\| \Psi_{n}\right\|^{2}_{2}+\frac{1}{4}\left\|\theta_{n} \right\|^{2}+C\left|\nabla\widetilde{\phi}_{n} \right|^{2}_{L^{2}}\left\|\widetilde{\phi}_{n} \right\|^{2}_{2}\left| \theta_{n}\right|^{2}_{L^{2}}+C\left\|\widetilde{\phi}_{n} \right\|^{2}_{2}\left\| \Psi_{n}\right\|^{2}_{1},
                             \end{aligned}
                             \end{equation}
                             \begin{equation}\label{e11}
                             \begin{aligned}
                              \left|I_{13}\right| &\leq C\left|\nabla(\widetilde{\phi}_{n}-\phi) \right|^{1/2}_{L^{2}}\left\|(\widetilde{\phi}_{n}-\phi) \right\|^{1/2}_{2} \left(\left\|\text{u} \right\|_{\mathbb{L}^{4}} +\left\|\textbf{a} \right\|_{\mathbb{L}^{4}}  \right)\left( \left\|\Psi_{n} \right\|_{2}+\left|\Psi_{n} \right|_{L^{2}}\right),
                              \end{aligned}
                              \end{equation}
                              \begin{equation}\label{e12}
                               \begin{aligned}
                               \left|I_{12}\right|\leq C\left\|\widetilde{\text{u}}_{n}-\text{u}\right\|_{\mathbb{L}^{4}} \left|\nabla\widetilde{\phi}_{n} \right|^{1/2}_{L^{2}}\left\|\widetilde{\phi}_{n} \right\|^{1/2}_{2}\left( \left\|\Psi_{n} \right\|_{2}+\left|\Psi_{n} \right|_{L^{2}}\right).
                                \end{aligned}
                                \end{equation}
                Note that
                             \begin{equation}\label{e13}
                             -(f_{\vartheta}(\widetilde{\phi}_{n})-f_{\vartheta}(\phi_{n}), \text{A}_{\vartheta}\Psi_{n})\leq \frac{1}{4} \left\| \Psi_{n}\right\|^{2}_{2}+ C\left|\nabla\Psi_{n} \right|^{2}_{L^{2}},\\
                             \end{equation}
                             \begin{equation}\label{e14}
                             \left| (f_{\vartheta}(\widetilde{\phi}_{n})-f_{\vartheta}(\phi),\text{A}_{\vartheta}\Psi_{n} )\right| \leq \frac{1}{4}\left\|\Psi_{n} \right\|^{2}_{2}+C_{f}\left|\nabla(\widetilde{\phi}_{n}- \phi)\right|^{2}_{L^{2}},
                             \end{equation}
                           \begin{equation}\label{e15}
                           \left| (f_{\vartheta}(\widetilde{\phi}_{n})-f_{\vartheta}(\phi_{n}), \xi_{n})\right| \leq \frac{1}{3}\left|\xi_{n} \right|^{2}_{L^{2}}+C_{f} \left|\nabla\Psi_{n} \right|^{2}_{L^{2}},
                          \end{equation}
                             \begin{equation}\label{e16}
                             \left| (f_{\vartheta}(\widetilde{\phi}_{n})-f_{\vartheta}(\phi),\xi_{n} )\right| \leq \frac{1}{3}\left|\xi_{n} \right|^{2}_{L^{2}}+C_{f}\left|\nabla(\widetilde{\phi}_{n}- \phi)\right|^{2}_{L^{2}}.
                              \end{equation}
                 Let
                               \begin{equation*}
                               \begin{aligned}
                               \mathbf{Z}_{1}(t)&= \left|\theta_{n} \right|^{2}_{L^{2}}+\left| \nabla\Psi_{n}\right|^{2}_{L^{2}}, \ \ \mathbf{Z}_{2}(t)= \left\|\theta_{n} \right\|^{2}+\left\| \Psi_{n}\right\|^{2}_{2}+\left|\xi_{n} \right|^{2}_{L^{2}}\\
                               \widetilde{f}(s)&=C(1+\left|\mu_{n} \right|^{2}_{L^{2}}+\left|\mu_{n} \right|^{4}_{L^{2}}+\left\|\widetilde{\phi}_{n} \right\|_{2}^{2} \left|\nabla\widetilde{\phi}_{n} \right|^{2}_{L^{2}} +\left| \nabla\widetilde{\phi}_{n}\right|^{4}_{L^{2}}+\left\|\widetilde{\phi}_{n} \right\|^{2}_{2})\\
                               &+\max\left(3C_{0},C_{2},C_{4} \right)\left(1+\left\|(a,b) \right\|^{2} _{\mathcal{H}_{p}(\Gamma)}+\left\|\text{u} \right\|^{2}+\left\|(a,b) \right\|^{4} _{\mathcal{H}_{p}(\Gamma)}+\left\|\text{u}_{n} \right\|^{4} \right)\\
                               \sigma(t)&=\exp\left( -\int_{0}^{t}\widetilde{f}(s)ds\right) .
                               \end{aligned}
                               \end{equation*}     
           Applying the It\^{o} formula to the process $\sigma(t)\mathbf{Z}_{1}(t)$,  we derive using $\eqref{e1}-\eqref{e16}$ that
                      \begin{equation}\label{e17}
                      \begin{aligned}
                      &\sigma(t)\mathbf{Z}_{1}(t)+\int_{0}^{t}\sigma(s)\mathbf{Z}_{2}(s)ds+\int_{0}^{t}\sigma(s)\left|\mathcal{P}_{1}^{n}\text{g}-\mathcal{P}_{1}^{n}\text{g}^{*} \right|^{2}_{L^{2}}ds\\
                      &\leq \int_{0}^{t}\sigma(s)I_{2}ds+\int_{0}^{t}\sigma(s)I_{3}ds+ \int_{0}^{t}\sigma(s)I_{7}ds+\int_{0}^{t}\sigma(s)I_{8}ds+\int_{0}^{t}\sigma(s)I_{12}ds\\
                      &++\int_{0}^{t}\sigma(s)I_{13}ds+\int_{0}^{t}\sigma(s)I_{4}ds+\int_{0}^{t}\sigma(s)I_{9}ds+\int_{0}^{t}\sigma(s)I_{14}ds\\
                      &+C\int_{0}^{t}\sigma(s)\left|\widetilde{\text{u}}_{n}-\text{u} \right|^{2}_{L^{2}}ds
                       +2\int_{0}^{t}\sigma(s)(\mathcal{P}_{1}^{n}\text{g}_{n}-\mathcal{P}_{1}^{n}\text{g}^{*},\mathcal{P}_{1}^{n}\text{g}-\mathcal{P}_{1}^{n}\text{g}^{*})ds\\
                       &+C\int_{0}^{t}\sigma(s)\left|\nabla(\widetilde{\phi}_{n}- \phi)\right|^{2}_{L^{2}}ds+
                       +C\int_{0}^{t}\sigma(s)(f_{\vartheta}(\phi)-f^{*}(t),\text{A}_{\vartheta}\Psi_{n} )ds\\
                       &+C\int_{0}^{t}\sigma(s)(f_{\vartheta}(\phi)-f^{*}(t),\xi_{n} )ds
                       -2\int_{0}^{t}\sigma(s)(g(t,\text{v}_{n})-\text{g}^{*}(t), \theta_{n})dW_{s}.
                      \end{aligned}
                      \end{equation}
            It follows from $\eqref{e17}$ that
                     \begin{equation}\label{e18}
                      \begin{aligned}
                       &\mathbb{E}\sigma(t)\mathbf{Z}_{1}(t)+\mathbb{E}\int_{0}^{t}\sigma(s)\mathbf{Z}_{2}(s)ds+\mathbb{E}\int_{0}^{t}\sigma(s)\left|\mathcal{P}_{1}^{n}\text{g}-\mathcal{P}_{1}^{n}\text{g}^{*} \right|^{2}_{L^{2}}ds\\
                                  &\leq \mathbb{E}\int_{0}^{t}\sigma(s)I_{2}ds+\mathbb{E}\int_{0}^{t}\sigma(s)I_{3}ds+ \int_{0}^{t}\sigma(s)I_{7}ds+\mathbb{E}\int_{0}^{t}\sigma(s)I_{8}ds+\mathbb{E}\int_{0}^{t}\sigma(s)I_{12}ds\\
                                  &+\mathbb{E}\int_{0}^{t}\sigma(s)I_{13}ds+\mathbb{E}\int_{0}^{t}\sigma(s)I_{4}ds+\mathbb{E}\int_{0}^{t}\sigma(s)I_{9}ds+\mathbb{E}\int_{0}^{t}\sigma(s)I_{14}ds\\
                                  &+C\mathbb{E}\int_{0}^{t}\sigma(s)\left|\widetilde{\text{u}}_{n}-\text{u} \right|^{2}_{L^{2}}ds
                                   +2\mathbb{E}\int_{0}^{t}\sigma(s)(\mathcal{P}_{1}^{n}\text{g}_{n}-\mathcal{P}_{1}^{n}\text{g}^{*},\mathcal{P}_{1}^{n}\text{g}-\mathcal{P}_{1}^{n}\text{g}^{*})ds\\
                                   &+C\mathbb{E}\int_{0}^{t}\sigma(s)\left|\nabla(\widetilde{\phi}_{n}- \phi)\right|^{2}_{L^{2}}ds
                                   +C\mathbb{E}\int_{0}^{t}\sigma(s)(f_{\vartheta}(\phi)-f^{*}(t),\text{A}_{\vartheta}\Psi_{n} )ds\\
                                   &+C\mathbb{E}\int_{0}^{t}\sigma(s)(f_{\vartheta}(\phi)-f^{*}(t),\xi_{n} )ds.
                                 \end{aligned}
                                 \end{equation}
         Now, we will show that the right-hand side of this inequality tends to zero as $n\rightarrow\infty$.
       Note that $\sigma(s)\leq \mathcal{G}^{3}(t) \leq\mathcal{G}^{2}(t)$. From the estimate $\eqref{e2}$, we deduce that             
        \begin{equation}\label{gty}
                  \begin{aligned}
                  \mathbb{E}\int_{0}^{t}\sigma(s)I_{2}ds&\leq C\left(\mathbb{E}\int_{0}^{T}\mathcal{G}^{3}(s)\left|\widetilde{\text{u}}_{n}-\text{u} \right|^{2}_{\mathbb{L}^{4}} \left(\left\|\text{u} \right\| +\left\|\textbf{a} \right\|_{\mathbb{H}^{1}}\right) \right)^{1/2} \\
                  &\times \left(\mathbb{E}\int_{0}^{T}\mathcal{G}^{3}(s)\left(\left\|\text{u} \right\| +\left\|\textbf{a} \right\|_{\mathbb{H}^{1}}  \right)\left(\left\|\widetilde{\text{u}}_{n} \right\|^{2}_{\mathbb{L}^{4}} +\left\|\text{u}_{n} \right\|^{2}_{\mathbb{L}^{4}} \right)ds\right)^{1/2}.
                  \end{aligned}
                  \end{equation}
                  Using the estimates $\eqref{g50a}-\eqref{n2b0}$, we deduce  that
                  \begin{equation}\label{mo1}
                  \begin{aligned}
                  &\mathbb{E}\int_{0}^{T}\mathcal{G}^{3}(s)\left|\widetilde{\text{u}}_{n}-\text{u} \right|^{2}_{\mathbb{L}^{4}} \left(\left\|\text{u} \right\| +\left\|\textbf{a} \right\|_{\mathbb{H}^{1}}\right)ds\\
                  &\leq \left(\mathbb{E}\sup_{t\in [0,T]}\mathcal{G}^{2}(s)\left|\widetilde{\text{u}}_{n}-\text{u} \right|^{2}\int_{0}^{T}\mathcal{G}^{2}(s) \left(\left\|\text{u} \right\|^{2} +\left\|\textbf{a} \right\|^{2}_{\mathbb{H}^{1}}\right) ds\right)^{1/2} \\
                  &\times \left(\mathbb{E}\int_{0}^{T}\mathcal{G}^{2}(s)\left\|\widetilde{\text{u}}_{n}-\text{u} \right\|^{2} ds\right)^{1/2} \leq C\left(\mathbb{E}\int_{0}^{T}\mathcal{G}^{2}(s)\left\|\widetilde{\text{u}}_{n}-\text{u} \right\|^{2}ds\right)^{1/2}.
                  \end{aligned}
                  \end{equation} 
                   Now, we can prove that there exists a constant $C>0$ such that
                             \begin{equation}\label{m0}
                             \left(\mathbb{E}\int_{0}^{T}\mathcal{G}^{3}(s)\left(\left\|\text{u} \right\| +\left\|\textbf{a} \right\|_{\mathbb{H}^{1}}  \right)\left(\left\|\widetilde{\text{u}}_{n} \right\|^{2}_{\mathbb{L}^{4}} +\left\|\text{u}_{n} \right\|^{2}_{\mathbb{L}^{4}} \right)ds\right)^{1/2}\leq C.
                             \end{equation}
             We derive from $\eqref{mo1}-\eqref{m0}$
                             \begin{equation}\label{k2a}
                             J_{1}=\mathbb{E}\int_{0}^{t}\sigma(s)I_{2}ds\leq C\left(\mathbb{E}\int_{0}^{T}\mathcal{G}^{2}(s)\left\|\widetilde{\text{u}}_{n}-\text{u} \right\|^{2}ds\right)^{1/4}.
                             \end{equation} 
                  Similarly, we can prove that
                                           \begin{equation}\label{k2b}
                                    J_{2}= \mathbb{E}\int_{0}^{t}\sigma(s)I_{3}ds\leq C\left(\mathbb{E}\int_{0}^{T}\mathcal{G}^{2}(s)\left\|\widetilde{\text{u}}_{n}-\text{u} \right\|^{2}ds\right)^{1/2}.
                                \end{equation}
                                Note that
                                \begin{equation}\label{k2c}
                               \begin{aligned}
                                  C\mathbb{E}\int_{0}^{t}\sigma(s)I_{7}(s)ds
                                  &\leq C \left( \mathbb{E}\int_{0}^{T}\mathcal{G}^{2}(s)\left|(\widetilde{\mu}_{n}-\mu) \right|^{2}_{L^{2}}ds\right) ^{1/2}\\
                                &\times \left( \mathbb{E}\int_{0}^{T}\mathcal{G}^{2}(s)\left(\left\|\text{u}_{n} \right\|^{2}_{\mathbb{L}^{4}}+\left\|\widetilde{\text{u}}_{n} \right\|^{2}_{\mathbb{L}^{4}}\right)  \left\|\nabla\widetilde{\phi}_{n} \right\|^{2}_{L^{4}}ds\right)^{1/2}.
                              \end{aligned}
                          \end{equation}
                          Thanks to the Ladyzhenkaya, Galiardo-Nirenberg-sobolev inequalities and using the estimates $\eqref{g50a}-\eqref{n2b0}$,  we find that
                          \begin{equation}\label{po}
                          \left( \mathbb{E}\int_{0}^{T}\mathcal{G}^{3}(s)\left(\left\|\text{u}_{n} \right\|^{2}_{\mathbb{L}^{4}}+\left\|\widetilde{\text{u}}_{n} \right\|^{2}_{\mathbb{L}^{4}}\right)  \left\|\nabla\widetilde{\phi}_{n} \right\|^{2}_{L^{4}}ds\right)\leq C.
                          \end{equation}
                    Hence we have 
                      \begin{equation}
                           J_{3}= C\mathbb{E}\int_{0}^{t}\sigma(s)I_{7}(s)ds
                          \leq C \left( \mathbb{E}\int_{0}^{T}\mathcal{G}^{2}(s)\left|(\widetilde{\mu}_{n}-\mu) \right|^{2}_{L^{2}}ds\right) ^{1/2}.
                         \end{equation}  
                     Thanks to the Ladyzhenkaya, Galiardo-Nirenberg-Sobolev inequalities and using the estimates $\eqref{g50a}-\eqref{n2b0}$,  we find that
                          \begin{equation}\label{k2aa}
                          \begin{aligned}
                   J_{4}= \mathbb{E}\int_{0}^{t}\sigma(s)I_{8}ds&\leq C\left(\mathbb{E}\int_{0}^{T}\mathcal{G}^{2}(s)\left|\nabla(\widetilde{\phi}_{n}-\phi) \right|^{2}_{L^{2}}ds\right)^{1/2}\\
                   &+C\left(\mathbb{E}\int_{0}^{T}\mathcal{G}^{2}(s)\left\|(\widetilde{\phi}_{n}-\phi) \right\|^{2}_{2}\left|\nabla(\widetilde{\phi}_{n}-\phi) \right|^{2}_{L^{2}}ds\right)^{1/4}.
                  \end{aligned} 
                  \end{equation}
               Using a similarly argument as  $\eqref{gty}$, we arrive at
                                   \begin{equation}\label{k2a0}
                                              J_{5}=\mathbb{E}\int_{0}^{t}\sigma(s)I_{12}ds\leq C\left(\mathbb{E}\int_{0}^{T}\mathcal{G}^{2}(s)\left\|\widetilde{\text{u}}_{n}-\text{u} \right\|^{2}ds\right)^{1/4},
                                              \end{equation}
                 and
                  \begin{equation}\label{k2aa0}
                                       J_{6}= \mathbb{E}\int_{0}^{t}\sigma(s)I_{13}ds\leq C\left(\mathbb{E}\int_{0}^{T}\mathcal{G}^{2}(s)\left|\nabla(\widetilde{\phi}_{n}-\phi) \right|^{2}_{L^{2}}\left\|(\widetilde{\phi}_{n}-\phi) \right\|^{2}_{2}ds\right)^{1/4}.
                                                \end{equation}
                It follows from $\eqref{k2a}-\eqref{k2aa0}$ that the terms $J_{i}, i=1,...,5$ converge to zero as $n\rightarrow\infty$ by $\eqref{p3}$. \\
                It follows from $\eqref{p8}$ and $\eqref{p9}$ that
                             \begin{equation}\label{qsd}
                              \begin{aligned}
                      \mathcal{G}^{2}(\Pi_{n}(\text{u},\phi)-(\text{u}_{n},\phi_{n}))&\rightharpoonup 0\ \ \mbox{weakly in }\ \ L^{2}(\Omega\times(0,T), \mathbb{V}), \ \ n\rightarrow\infty,\\
                       \mathcal{G}^{2}(\widetilde{\mu}_{n}- \mu_{n})&\rightharpoonup 0\ \ \mbox{weakly in }\ \ L^{2}(\Omega\times(0,T), L^{2}(D)), \ \ n\rightarrow\infty.
                       \end{aligned}
                      \end{equation}  
          Thanks to $\eqref{qsd}$, we derive  as in \cite{Zhang} that
                                     \begin{equation*}
                                     C\mathbb{E}\int_{0}^{t}\sigma(s)(f_{\vartheta}(\phi)-f^{*}(t),\text{A}_{\vartheta}\Psi_{n} )ds+C\mathbb{E}\int_{0}^{t}\sigma(s)(f_{\vartheta}(\phi)-f^{*}(t),\xi_{n} )ds=0 \ \ \mbox{as}\ \ n\rightarrow\infty.
                                     \end{equation*}
                Now, the operators  $\mathcal{G}^{2}((\text{y}\cdot\nabla)\text{y}-\text{B}^{*}_{0})$ and $\mathcal{G}^{2}(\mu\nabla\phi-\text{R}^{*}_{0})$ belong to  $L^{2}(\Omega\times(0,T), \mathbb{V}^{*}_{div})$. Thanks to $\eqref{p10}$ and $\eqref{e4}$, it follows that
                                 \begin{equation*}
                                 \mathbb{E}\int_{0}^{t}\sigma(s)I_{4}ds+\mathbb{E}\int_{0}^{t}\sigma(s)I_{9}ds=0\ \ \mbox{as}\ \ n\rightarrow\infty.
                                 \end{equation*}
             Similarly, we deduce that the operators $\mathcal{G}^{2}((\text{y}\cdot\nabla)\phi-\text{B}^{*}_{0})$ belongs to  $L^{2}(\Omega\times(0,T), \text{V}^{*}_{1})$. Thanks to  $\eqref{p6}$ and $\eqref{p10}$, we deduce that
                                 \begin{equation*}
                                 \mathbb{E}\int_{0}^{T}\sigma(s)I_{4}ds=0\ \ \mbox{as}\ \ n\rightarrow\infty.
                                 \end{equation*}
                     Thanks to $\eqref{p9}$, we deduce that
                                 \begin{equation*}
                                 C\mathbb{E}\int_{0}^{T}\sigma(s)\left\|(\widetilde{\text{u}}_{n},\widetilde{\phi}_{n})-(\text{u},\phi) \right\|^{2}_{\mathbb{Y}}ds\rightarrow 0.
                                 \end{equation*}
                  Due to the convergence results $\eqref{p6}$, $\eqref{p8}$, $\eqref{p9}$, $\eqref{p10}$ and $\eqref{p11}$, we derive that
                                 \begin{equation}\label{ft}
                                 \begin{aligned}
                                 \mathcal{G}^{2}\mathcal{P}_{1}^{n}(\text{g}_{n}-\text{g}^{*})&\rightharpoonup 0\ \ \mbox{weakly in }\ \ L^{2}(\Omega\times(0,T), \mathbb{H}^{m}_{div}),\\
                                 \mathcal{G}^{2}\mathcal{P}_{1}^{n}(\text{g}-\text{g}^{*})&\rightarrow \text{g}-\text{g}^{*}\ \ \mbox{strongly in }\ \ L^{2}(\Omega\times(0,T), \mathbb{H}^{m}_{div}),
                                 \end{aligned}
                                 \end{equation}
                                 which gives
                                 \begin{equation*}
                                 2\mathbb{E}\int_{0}^{t}\sigma(s)(\mathcal{P}_{1}^{n}\text{g}_{n}-\mathcal{P}_{1}^{n}\text{g}^{*},\mathcal{P}_{1}^{n}\text{g}-\mathcal{P}_{1}^{n}\text{g}^{*})ds\rightarrow 0\ \ \mbox{as}\ \ n\rightarrow\infty.
                                 \end{equation*}
                Hence, we have the following strong convergences
                                 \begin{equation}\label{fc}
                                 \lim\limits_{n\rightarrow\infty}\mathbb{E}\sigma(t)\mathbf{Z}_{1}(t)=\lim\limits_{n\rightarrow\infty}\mathbb{E}\int_{0}^{t}\sigma(s)\mathbf{Z}_{2}(s)ds=0,
                                 \end{equation}
                  for all $t\in(0,T)$. We derive from $\eqref{p8}$ and $\eqref{fc}$ that
                                 \begin{equation}\label{hg}
                                 \begin{aligned}
                             \lim\limits_{n\rightarrow\infty}\mathbb{E}\left( \sigma(t)\left\|(\text{u}_{n},\phi_{n})(t)-(\text{u},\phi)(t) \right\|^{2}_{\mathbb{Y}}\right) &=0,\\
                                 \mathbb{E}\int_{0}^{t}\sigma(s)\left\|(\text{u}_{n},\phi_{n})(s)-(\text{u},\phi)(s) \right\|^{2}_{\mathbb{V}}ds &=0,\\
                                 \mathbb{E}\int_{0}^{t}\sigma(s)\left|\text{g}(s,\text{v})-\text{g}^{*}(s) \right|^{2}_{L^{2}}ds=0.
                                 \end{aligned}
                                 \end{equation}
        Now, using the fact that $0< \sigma(s)\leq 1$, we derive that
                                 \begin{equation}\label{ght}
                                 \text{g}(t,\text{v})=\text{g}^{*}(t) \ \ \mbox{a.e. in}\ \ (\omega,t)\in\Omega\times(0,T).
                                 \end{equation}  
       It follows from $\eqref{e4}$ and $\eqref{hg}$, that $\sigma(t)(\text{v}\cdot\nabla)\text{v}=\sigma(t)\text{B}^{*}_{0}(t)$ \ \ $\mbox{a.e. in}\ \ (\omega,t)\in\Omega\times(0,T)$, that implies
                               \begin{equation}\label{ght1}
                                  (\text{v}\cdot\nabla)\text{v}=\text{B}^{*}_{0}(t) \ \ \mbox{a.e. in}\ \ (\omega,t)\in\Omega\times(0,T).
                                   \end{equation}
                                    Similarly, we can prove that the operators $\sigma(t)(\mu\nabla\phi)=\sigma(t)\text{R}^{*}_{0}(t)$ and $\sigma(t)(\text{v}\cdot\nabla)\phi=\sigma(t)\text{B}^{*}_{0}(t)$ \ \ $\mbox{a.e. in }\ \ (\omega,t)\in\Omega\times(0,T)$, we infer that
                                \begin{equation}\label{ght2}
                                \begin{aligned}
                                 \mu\nabla\phi&=\text{R}^{*}_{0}(t)\ \ \mbox{a.e. in}\ \ (\omega,t)\in\Omega\times(0,T),\\
                                 (\text{v}\cdot\nabla)\phi&=\text{B}^{*}_{0}(t) \ \ \mbox{a.e. in}\ \ (\omega,t)\in\Omega\times(0,T).
                                 \end{aligned}
                               \end{equation}
                               From $\eqref{ght}, \eqref{ght1}, \eqref{ght2}$, we obtain the following system
                               \begin{equation}\label{Ton2a}
                               \begin{cases}
                              (\text{v}(t),\text{w})=(\text{v}_{0},\text{w}) +\int_{0}^{t}\left[  -((\text{v}(s),\text{w}))+\int_{\Gamma}b(\text{w}\cdot\tau)d\gamma -\left\langle  (\text{v}(s)\cdot\nabla\text{v}(s)),\text{w}\right\rangle  \right]   ds\\
                              +\int_{0}^{t}\left\langle (\mu\nabla\phi(s), \text{w}\right\rangle ds
                              +\int^{t}_{0}(\text{g}(s,\text{v}(s)),\text{w})dW_{s},\\
                              (\phi(t),\psi)=(\phi_{0},\psi) -\int_{0}^{t}\left( \mu(s),\psi\right)  ds- \int_{0}^{t}\left\langle (\text{v}(s)\cdot\nabla\phi(s)),\psi\right\rangle   ds,\\
                             \mu=\text{A}_{\vartheta}\phi + f_{\vartheta}(\phi),
                             \end{cases}
                              \end{equation}
                    for  all $t\in[0,T]$, $\text{w}\in\mathbb{V}_{div}$, $\psi\in\text{V}_{1}$ and $\mathbb{P}$-a.e. in $\Omega.$ This completes the proof of Theorem $\ref{Theo2}$.\newline
                    The uniqueness of the solution $(\text{v}, \phi)$ follows the stability result established in the next theorem. Now, we assume that $(\text{v}_{1},\phi_{1}), (\text{v}_{2},\phi_{1})$ are two strong solutions to $\eqref{NS1}$. We define by
                                        \begin{math}
                                       \text{v}=\text{v}_{1}-\text{v}_{2}, \ \ \phi=\phi_{1}-\phi_{2}
                                         \end{math}
                          the difference of two strong solutions to $\eqref{NS1}$.    
            \begin{thm}\label{12}
                                    Let  $(\text{v}_{1},\phi_{1})=(\text{u}_{1}+\textbf{a}_{1},\phi_{1})$, $(\text{v}_{2},\phi_{2})=(\text{u}_{2}+\textbf{a}_{2},\phi_{2})$ with
                                       \begin{equation*}
                                     (\text{u}_{1},\phi_{1}),(\text{u}_{2},\phi_{2})\in \mathcal{C}([0,T], \mathbb{Y})\cap L^{4}(0,T; \mathbb{V}), \ \ \mathbb{P}-a.e. \ \ in \ \ \Omega,
                                       \end{equation*}
                                two solutions of $\eqref{NS1}$, satisfying the estimates $\eqref{n2b0}, \eqref{g50a}$ with two corresponding boundary conditions $a_{1},a_{2},b_{1},b_{2}$ and the initial conditions
                                     \begin{equation*}
                                (\text{v}_{1,0},\phi_{1,0})=(\text{u}_{1,0}+\textbf{a}_{1}(0),\phi_{1,0}), \ \  (\text{v}_{2,0},\phi_{2,0})=(\text{u}_{2,0}+\textbf{a}_{2}(0),\phi_{2,0}).
                                 \end{equation*}
                               Then there exist a strictly positive function $h(t)\in L^{1}(0,T)$ $\mathbb{P}$-a.e. in $\Omega$, depending only on the data, such that the following estimate
                                   \begin{equation}\label{g50}
                                   \begin{aligned}
                                \mathbb{E}\sup\limits_{s\in[0,t]}\mathcal{H}^{2}(s)\left\|(\text{v},\phi)(s) \right\|^{2}_{\mathbb{Y}}&+2\int_{0}^{t}\mathcal{H}^{2}(s)\left[ \left\|(\text{v},\phi)(s) \right\|^{2}_{\mathbb{V}}+\left|\mu(s) \right|^{2}_{L^{2}} \right] ds \\
                                 &\leq C\left(\mathbb{E}\mathcal{E}(\text{u}_{0},\phi_{0})+\mathbb{E}\int_{0}^{t}\mathcal{H}^{2}(s)\left\|(a,b) \right\|^{2}_{\mathcal{H}_{p}(\Gamma)} ds \right),
                               \end{aligned}
                                \end{equation}
                           is valid with the function $\mathcal{H}$ defined as
                                \begin{equation}
                                \mathcal{H}(t)=\exp\left(-\int_{0}^{t}h(s)ds \right) \ \ \mbox{with} \ \  h(t)\in L^{1}(0,T)\ \ \mathbb{P}-a.e. \, in \ \ \Omega.
                                 \end{equation}
                                \end{thm}
                                \begin{prev}
                         Let us denote  $\textbf{a}=\textbf{a}_{1}-\textbf{a}_{2}$, $p=p_{1}-p_{2}$ where $\textbf{a}_{1}$ and $\textbf{a}_{2}$ are the solutions of system $\eqref{S1}$ with two corresponding boundary conditions $(a_{1},b_{1})$ and  $(a_{2},b_{2})$. We easily verify that the functions $\text{u}=\text{v}-\textbf{a}, \phi=\phi_{1}-\phi_{2} $ and $\mu=\mu_{1}-\mu_{2}$ satisfy the following system 
                         \begin{equation}\label{g2aa0}
                         \begin{cases}
                         d(\text{u},\varphi)=\left[-((\text{u},\varphi))-\left\langle (\text{v}_{1}\cdot\nabla\right\rangle \text{v}_{1}-\text{v}_{2}\cdot\nabla\text{v}_{2},\varphi) \right] dt
                         -((\textbf{a},\varphi))dt-\left(\partial_{t}(\textbf{a}) ,\varphi \right)\\ +\int_{\Gamma}b(\varphi\cdot\tau)d\gamma dt
                           +\left\langle \mu_{1}\nabla\phi_{1}-\mu_{2}\nabla\phi_{2},\varphi\right\rangle dt+(g(t,\text{v}_{1})-g(t,\text{v}_{2}), \varphi)dW(t),\\
                           d(\phi,\psi)=-\left(  \mu,\psi \right)  dt -\left\langle(\text{v}_{1}\cdot\nabla\phi_{1}-\text{v}_{2}\cdot\nabla\phi_{2}), \psi \right\rangle dt,\\
                            \mu=\text{A}_{\vartheta}\phi+f_{\vartheta}(\phi_{1})-f_{\vartheta}(\phi_{2}),
                            \end{cases}
                           \end{equation}
                           for  all $t\in[0,T]$, $\varphi\in\mathbb{V}_{div}$, $\psi\in\text{V}_{1}$. and $\mathbb{P}$-a.e. in $\Omega$.\newline 
                By applying It\^{o}'s formula, the equation $\eqref{g2aa0}_{1}$, we obtain 
                                         \begin{equation}\label{k10a}
                                          \begin{aligned}
                                           d\left| \text{u}\right|^{2}_{L^{2}} &+2\left\| \text{u}\right\|^{2}dt=-2\left\langle  \text{v}_{1}\cdot\nabla\text{v}_{1}-\text{v}_{2}\cdot\nabla\text{v}_{2},\text{u}\right\rangle  dt \\
                                          &-2\left[ ((\textbf{a},\text{u}))-\left(\partial_{t}(\textbf{a}) ,\text{u} \right) +\int_{\Gamma}b(\text{u}\cdot\tau)d\gamma\right]dt \\
                                           &+2\left\langle \mu_{1}\nabla\phi_{1}-\mu_{2}\nabla\phi_{2},\text{u}\right\rangle dt+2(g(t,\text{v}_{1})-g(t,\text{v}_{2}), \text{u})dW(t)\\
                                         &+ \sum_{i=1}^{n}\left|(g(t,\text{v}_{1})-g(t,\text{v}_{2})), e_{i}) \right|^{2}dt.
                                             \end{aligned}
                                           \end{equation}
           Replacing $\psi$ in $\eqref{g2aa0}_{2}$ by $\text{A}_{\vartheta}\phi$ and multiplying the equation $\eqref{g2aa0}_{3}$ with $\mu-\text{A}_{\vartheta}\phi$ and adding these equations, we derive that  
          \begin{equation}\label{k3aa}
                          \begin{aligned}
                     d\left|\nabla\phi\right|^{2}_{L^{2}} &+2(\left\|\phi\right\|^{2}_{2}+\left|\mu \right|^{2}_{L^{2}})dt = -2\left\langle(\text{v}_{1}\cdot\nabla\phi_{1}-\text{v}_{2}\cdot\nabla\phi_{2}), \text{A}_{\vartheta}\phi \right\rangle dt\\
                        &-\left(f_{\vartheta}(\phi_{1})-f_{\vartheta}(\phi_{2})), \text{A}_{\vartheta}\phi\right)dt +\left(f_{\vartheta}(\phi_{1})-f_{\vartheta}(\phi_{2}),\mu \right) dt.
                      \end{aligned}
                      \end{equation}
                   Adding $\eqref{k3aa}$ and $\eqref{k10a}$, we obtain  
             \begin{equation}\label{gh0a}
                          \begin{aligned}
                     &d\left( \left|\text{u}\right|^{2}_{L^{2}}+ \left\|\nabla\phi\right\|^{2}_{L^{2}}\right) +2\left\| \text{u}\right\|^{2}dt+2(\left\|\phi\right\|^{2}_{2}+\left|\mu \right|^{2}_{L^{2}})dt\\
                          &=-2\left\langle  \text{v}_{1}\cdot\nabla\text{v}_{1}-\text{v}_{2}\cdot\nabla\text{v}_{2},\text{u}\right\rangle dt
                         -2\left\langle(\text{v}_{1}\cdot\nabla\phi_{1}-\text{v}_{2}\cdot\nabla\phi_{2}), \text{A}_{\vartheta}\phi \right\rangle dt\\
                          &+\left(f_{\vartheta}(\phi_{1})-f_{\vartheta}(\phi_{2})), \text{A}_{\vartheta}\phi\right)dt -\left(f_{\vartheta}(\phi_{1})-f_{\vartheta}(\phi_{2}),\mu \right) dt+2\left\langle \mu_{1}\nabla\phi_{1}-\mu_{2}\nabla\phi_{2},\text{u}\right\rangle dt \\
                         &-2\left[((\textbf{a},\text{u}))-\left(\partial_{t}(\textbf{a}) ,\text{u} \right) +\int_{\Gamma}b(\text{u}\cdot\tau)d\gamma\right]dt\\
                        &+2(g(t,\text{v}_{1})-g(t,\text{v}_{2}), \text{u})dW(t)
                       + \sum_{i=1}^{n}\left|(g(t,\text{v}_{1})-g(t,\text{v}_{2})), e_{i}) \right|^{2}dt.
                                         \end{aligned}
                                          \end{equation}  
           Note that     
           \begin{equation*}
                       \begin{aligned}
                       \left\langle  \text{v}_{1}\cdot\nabla\text{v}_{1}-\text{v}_{2}\cdot\nabla\text{v}_{2},\text{u}\right\rangle
                       &= \left\langle(\text{u}+\textbf{a})\cdot\nabla\text{u}_{1},\text{u} \right\rangle +\left\langle(\text{u}+\textbf{a}) \cdot\nabla\textbf{a}_{1},\text{u}\right\rangle +\left\langle (\text{u}_{2}+\textbf{a}_{2})\cdot\nabla\text{u},\text{u} \right\rangle\\ &+\left\langle(\text{u}_{2}+\textbf{a}_{2}) \cdot\nabla\textbf{a},\text{u}\right\rangle =I_{1}+I_{2}+I_{3}+I_{4},
                       \end{aligned}
                       \end{equation*} 
                       \begin{equation*}
                       \begin{aligned}
                       \left\langle(\text{v}_{1}\cdot\nabla\phi_{1}-\text{v}_{2}\cdot\nabla\phi_{2}), \text{A}_{\vartheta}\phi \right\rangle&=\left\langle(\text{u}+\textbf{a})\cdot\nabla\phi_{1},  \text{A}_{\vartheta}\phi \right\rangle+\left\langle(\text{u}_{2}+\textbf{a}_{2})\cdot\nabla\phi, \text{A}_{\vartheta}\phi \right\rangle\\
                       &=I_{5}+I_{6},
                       \end{aligned}
                       \end{equation*}
                       and 
                       \begin{equation*}
                       \left\langle \mu_{1}\nabla\phi_{1}-\mu_{2}\nabla\phi_{2},\text{u}\right\rangle=\left\langle\mu\nabla\phi_{1}, \text{u}\right\rangle +\left\langle\mu_{2}\nabla\phi, \text{u}\right\rangle=I_{7}+I_{8}. 
                       \end{equation*}
          Using the H\"{o}lder, Ladyzhenskaya, Young, Poincaré, Gagliardo-Nirenberg, Agmon inequalities, Sobolev embedding and Theorem $\ref{th01}$, we obtain the following estimates:  
         \begin{equation}\label{pm1}
                     \left| I_{7}\right|\leq \left\|\text{u} \right\|_{\mathbb{L}^{4}} \left|\mu \right| \left\|\nabla\phi_{1} \right\|_{L^{4}} \leq \frac{1}{2} \left|\mu \right|^{2}_{L^{2}}+\frac{1}{8}\left\| \text{u}\right\|^{2}+C\left\|\nabla\phi_{1}\right\|^{4}_{L^{4}}\left|\text{u}\right|^{2}_{L^{2}}, 
                     \end{equation}
                     \begin{equation}\label{pm2}
                      \left| I_{8}\right|\leq \left\|\text{u} \right\|_{\mathbb{L}^{4}} \left|\mu_{2} \right| \left\|\nabla\phi \right\|_{L^{4}} \leq \frac{1}{4}\left\|\text{u} \right\|^{2}+\frac{1}{4}\left\| \phi\right\|^{2}_{2}+C\left|\mu_{2} \right|^{2}_{L^{2}}\left| \nabla\phi\right|^{2}_{L^{2}} +C\left|\mu_{2} \right|^{4}_{L^{2}}\left| \nabla\phi\right|^{2}_{L^{2}},
                      \end{equation}
                     \begin{equation}\label{pm3}
                     \begin{aligned}
                     \left|I_{6}\right| &\leq \left\|\text{u}_{2}+\textbf{a}_{2} \right\|_{\mathbb{L}^{4}} \left| \nabla\phi\right|^{1/2}_{L^{2}}\left\|\phi\right|^{1/2}_{2} \left\|\phi \right\|_{2}
                      \leq \frac{1}{3} \left\| \phi\right\|^{2}_{2}+C\left\|\text{u}_{2}+\textbf{a}_{2} \right\|^{4}_{\mathbb{L}^{4}} \left| \nabla\phi\right|^{2}_{L^{2}} \\
                      &\leq \frac{1}{3} \left\| \phi\right\|^{2}_{2}+C\left( \left\|\text{u}_{2}\right\|^{4}+\left\|\textbf{a}_{2} \right\|^{4}\right)\left| \nabla\phi\right|^{2}_{L^{2}}\\
                      &\leq \frac{1}{3} \left\| \phi\right\|^{2}_{2}+C\left( \left\|\text{u}_{2}\right\|^{4}+\varUpsilon^{2} \right)\left| \nabla\phi\right|^{2}_{L^{2}}.
                      \end{aligned}
                     \end{equation}
                     We observe that $I_{5}=J_{1}+J_{2}$ and we estimate the terms $J_{1}$ and $J_{2}$:
                     \begin{equation}\label{pm4}
                     \begin{aligned}
                   \left|J_{1}\right| &\leq\left\|\text{u} \right\|_{\mathbb{L}^{4}} \left\|\nabla\phi_{1} \right\|_{L^{4}}\left\|\phi \right\|_{2} 
                   \leq\frac{1}{3} \left\| \phi\right\|^{2}_{2}+ \frac{1}{8}\left\| \text{u}\right\|^{2}+C\left|\nabla\phi_{1} \right|^{2}_{L^{2}}\left\|\phi_{1} \right\|^{2}_{2}\left|\text{u} \right|^{2}_{L^{2}},
                       \end{aligned}
                     \end{equation}
                    and 
                              \begin{equation}\label{pm5}
                              \begin{aligned}
                             \left|J_{2}\right| &\leq \left\|\textbf{a} \right\|_{\mathbb{L}^{4}} \left\|\nabla\phi_{1} \right\|_{L^{4}} \left\|\phi \right\|_{2}\leq \frac{1}{3} \left\| \phi\right\|^{2}_{2}+\left\|\textbf{a} \right\|^{2}_{\mathbb{L}^{4}} \left\|\nabla\phi_{1} \right\|^{2}_{L^{4}}\leq \frac{1}{3} \left\| \phi\right\|^{2}_{2}+\varUpsilon^{2} \left\|\nabla\phi_{1} \right\|^{2}_{L^{4}},
                             \end{aligned}
                              \end{equation}   
              \begin{equation}
                        \begin{aligned}
                        \left|I_{1}\right| &\leq \frac{1}{8}\left\| \text{u}\right\|^{2}+\left\|\text{u}_{1} \right\|^{2} \left|\text{u} \right|^{2}_{L^{2}} +\left|\text{u}\right|^{2}_{L^{2}} \left\| \text{u}_{1}\right\|^{2} \left\|\textbf{a} \right\| ^{2}\\
                        &\leq  \frac{1}{8}\left\| \text{u}\right\|^{2}+\left\|\text{u}_{1} \right\|^{2} \left|\text{u} \right|^{2}_{L^{2}} +\left|\text{u}\right|^{2}_{L^{2}} \left\| \text{u}_{1}\right\|^{2}\varUpsilon,
                        \end{aligned}
                        \end{equation} 
                        \begin{equation}
                        \begin{aligned}
                                  \left|I_{2}\right| &\leq \frac{1}{8}\left\| \text{u}\right\|^{2}+\left\|\textbf{a}_{1}\right\|^{2} \left|\text{u} \right|^{2}_{L^{2}} +\left|\text{u}\right|^{2}_{L^{2}} \left\| \textbf{a}_{1}\right\|^{2} \left\|\textbf{a} \right\| ^{2}\\
                                  &\leq\frac{1}{8}\left\| \text{u}\right\|^{2}+\varUpsilon \left|\text{u} \right|^{2}_{L^{2}} +\left|\text{u}\right|^{2}_{L^{2}} \varUpsilon^{2},
                                  \end{aligned}
                                  \end{equation}  
                              \begin{equation}
                              \begin{aligned}
                           \left|I_{4}\right| &\leq \frac{1}{8}\left\| \text{u}\right\|^{2}+\left|\text{u}\right|^{2}_{L^{2}} \left\| \text{u}_{2}\right\|^{2} \left\|\textbf{a} \right\| ^{2} +\left|\text{u}\right|^{2}_{L^{2}} \left\| \textbf{a}_{2}\right\|^{2} \left\|\textbf{a} \right\|^{2}\\
                           &\leq  \frac{1}{8}\left\| \text{u}\right\|^{2}+\left|\text{u}\right|^{2}_{L^{2}} \left\| \text{u}_{2}\right\|^{2} \varUpsilon +\left|\text{u}\right|^{2}_{L^{2}} \varUpsilon^{2},
                           \end{aligned}
                         \end{equation} 
                           \begin{equation}
                           \begin{aligned}
                           \left|((\textbf{a},\text{u})) \right| &\leq \frac{1}{8}\left\| \text{u}\right\|^{2}+C\left\| \textbf{a}\right\|^{2}\leq  \frac{1}{8}\left\| \text{u}\right\|^{2}+C\varUpsilon,\\
                           \left|(\partial_{t}\textbf{a},\text{u}) \right| &\leq  \frac{1}{2}\left| \text{u}\right|^{2}_{L^{2}}+\frac{1}{2}\left| \partial_{t}\textbf{a}\right|^{2}_{L^{2}}\leq \frac{1}{2}\left| \text{u}\right|^{2}_{L^{2}}+C\varUpsilon.
                           \end{aligned}
                           \end{equation} 
                            We have
                              \begin{equation}\label{e5aq}
                              \begin{aligned}
                           &\sum_{i=1}^{n}\left|(g(t,\text{v}_{1})-g(t,\text{v}_{2})), e_{i}) \right|^{2}_{L^{2}}
                           &\leq K\left|\text{u}\right|^{2}_{L^{2}}+K\left|\textbf{a} \right|^{2}_{L^{2}}\leq K\left|\text{u}\right|^{2}_{L^{2}}+K\varUpsilon,
                          \end{aligned}
                        \end{equation} 
                      \begin{equation}\label{e14aa}
                        \left| (f_{\vartheta}(\phi_{1})-f_{\vartheta}(\phi_{2}),\text{A}_{\vartheta}\phi )\right| \leq \frac{1}{3}\left\|\phi \right\|^{2}_{2}+C_{f}\left|\nabla\phi\right|_{L^{2}}^{2},
                        \end{equation}
                          \begin{equation}\label{e150}
                          \left| (f_{\vartheta}(\phi_{2})-f_{\vartheta}(\phi_{1}), \mu)\right| \leq \frac{1}{2}\left|\mu \right|^{2}_{L^{2}}+C_{f} \left|\nabla\phi \right|^{2}_{L^{2}}.
                                        \end{equation}                      
                  Since $\text{v}_{2}$ is divergence free, $\text{v}_{2}\cdot n=a_{2}$ and $\text{u}_{2}\cdot n=0$, we deduce that
                                 \begin{equation}\label{mp0}
                                 \int_{D}I_{3}dx= \int_{D}\left[ (\text{u}_{2}+\textbf{a}_{2})\cdot\nabla\text{u}\right] \text{u}dx=\int_{\Gamma}\frac{a_{2}}{2}(\text{u}\cdot\tau)^{2}d\gamma.
                                 \end{equation} 
              From $\eqref{mp0}$, we derive that
                                \begin{equation}\label{mp6}
                                \left|I_{3} \right|\leq \left|a_{2} \right|_{L^{2}(\Gamma)}  \left|\text{u} \right|_{L^{2}(\Gamma)}\leq C \left\|a_{2} \right\|^{2}_{W^{1-1/p}(\Gamma)}  \left|\text{u} \right|^{2}_{L^{2}}+ \frac{1}{8}\left\| \text{u}\right\|^{2}. 
                                \end{equation} 
             Let us define the functions  
              \begin{equation*}
                             \begin{aligned}
                           \mathbf{Z}_{1}(t)&= \left|\text{u}\right|^{2}_{L^{2}}+\left| \nabla\phi\right|^{2}_{L^{2}}, \ \ \mathbf{Z}_{2}(t)= \left\|\text{u} \right\|^{2}+\left\| \phi\right\|^{2}_{2}+\left|\mu \right|^{2}_{L^{2}},\\
                          \widetilde{f}(s)&=C\left(1+\left\|\nabla\phi_{1}\right\|^{4}_{L^{4}}+\left|\mu_{2} \right|^{4}_{L^{2}}+\left|\mu_{2} \right|^{2}_{L^{2}}+ \left\|\text{u}_{2}\right\|^{4}+\varUpsilon^{2}+\left\|\nabla\phi_{1} \right\|^{2}_{L^{2}}\left\|\phi_{1} \right\|^{2}_{2}+\varUpsilon^{2} \left\|\nabla\phi_{1} \right\|^{2}_{L^{4}} \right)\\
                          &+C\left(\left\| \text{u}_{1}\right\|^{2}\varUpsilon+\left\| \text{u}_{1}\right\|^{2}+\varUpsilon+\varUpsilon^{2}+\left\| \text{u}_{2}\right\|^{2}\varUpsilon \right),\\ 
                             \mathcal{H}(t)&=\exp(-\int_{0}^{t}\widetilde{f}(s)ds), \ \ \varUpsilon=\max\left( \left\| (a_{1},b_{1})\right\|^{2}_{\mathcal{H}_{p}(\Gamma)}, \left\| (a_{2},b_{2})\right\|^{2}_{\mathcal{H}_{p}(\Gamma)}\right) .
                                \end{aligned}
                               \end{equation*}
                          Now, we consider 
                           \begin{equation}\label{g7aa}
                          \tau^{1}_{N}=\inf\{t\geq 0: h_{1}(t)\geq N\}\wedge T,\ \ \tau^{2}_{N}=\inf\{t\geq 0: h_{2}(t)\geq N\}\wedge T, 
                          \end{equation}
                          where
                           \begin{equation*}
                           \begin{aligned}
                 h_{1}(t)&=\mathcal{G}^{2}(t)\mathcal{E}((\text{u}_{1},\phi_{1})(t))+2\int_{0}^{t}\mathcal{G}^{2}(s)\left[ \left\|(\text{u}_{1},\phi_{1})(s) \right\|^{2}_{\mathbb{V}}+\left|\mu_{1}(s) \right|^{2}_{L^{2}} \right] ds.\\
                 h_{2}(t)&=\mathcal{G}^{2}(t)\mathcal{E}((\text{u}_{2},\phi_{2})(t))+2\int_{0}^{t}\mathcal{G}^{2}(s)\left[ \left\|(\text{u}_{2},\phi_{2})(s) \right\|^{2}_{\mathbb{V}}+\left|\mu_{2}(s) \right|^{2}_{L^{2}} \right] ds,
                           \end{aligned}
                           \end{equation*}
                          and
                             \begin{equation}\label{g7m}
                             \tau_{N}=\tau^{1}_{N}\wedge \tau^{2}_{N}.
                           \end{equation}    
            Applying It\^{o}'s formula to the process $\mathcal{H}^{2}(s)\mathbf{Z}_{1}(s)$ and using $\eqref{pm1}-\eqref{mp6}$, we derive that   
                        \begin{equation}\label{A110m}
                                              \begin{aligned}
                                              \mathcal{H}^{2}(s)\mathbf{Z}_{1}(s)&+2\int_{0}^{s}\mathcal{H}^{2}(r)\mathbf{Z}_{2}(r) dr
                                              \leq\mathcal{E}(\text{u}_{0},\phi_{0})
                                               +C\int_{0}^{s} \mathcal{H}^{2}(r)\varUpsilon(r)dr\\
                                              &+2\int_{0}^{s}\mathcal{H}^{2}(r)(g(r,\text{v}_{1}(r))-g(r,\text{v}_{2}(r)), \text{u}(r))dW_{r}.
                                              \end{aligned}
                                               \end{equation}
                   We consider the sequence $\tau_{N}$ of the stopping times introduce in $\eqref{g7m}$. Thanks to \textbf{(H1)} and Burholder-Davis-Gundy's inequality, we infer that
                                            \begin{equation*}
                                               \begin{aligned}
                                               &\mathbb{E}\sup\limits_{s\in[0,\tau_{N}\wedge t]}\left|\int_{0}^{s}\mathcal{H}^{2}(r)(g(r,\text{v}_{1}(r))-g(r,\text{v}_{2}(r))), \text{u})dW_{r} \right|\\
                                               &\leq\mathbb{E} \left(\int_{0}^{\tau_{N}\wedge t}\mathcal{H}^{4}(s)\left| (g(s,\text{v}_{1}(s))-g(s,\text{v}_{2}(s))), \text{u}(s))\right|^{2} ds \right)^{1/2} \\
                                               &\leq \mathbb{E}\sup\limits_{s\in[0,\tau_{N}\wedge t]}\mathcal{H}(s)\left| \text{u}(s)\right|_{L^{2}}\left(\int_{0}^{\tau_{N}\wedge t}\mathcal{H}^{2}(s)\left| g(s,\text{v}_{1}(s))-g(s,\text{v}_{2}(s))\right|^{2}_{L^{2}} ds \right)^{1/2}\\
                                               &\leq \frac{1}{2}\mathbb{E}\sup\limits_{s\in[0,\tau_{N}\wedge t]}\mathcal{H}^{2}(s)\mathbf{Z}_{1}(s)+C\mathbb{E}\int_{0}^{\tau_{N}\wedge t}\mathcal{H}^{2}(s)\left(\mathbf{Z}_{1}(s)+\varUpsilon(s) \right)ds.
                                               \end{aligned}
                                               \end{equation*}
                        Taking the supremum of the relation $\eqref{A110m}$ for $s\in[0,\tau_{N}\wedge t]$, next we take the expectation, we derive that  
                         \begin{equation}\label{A11bm}
                            \begin{aligned}
                            &\mathbb{E}\sup\limits_{s\in[0,\tau_{N}\wedge t]}\left[ \mathcal{H}^{2}(s)\mathbf{Z}_{1}(s)\right] +2\mathbb{E}\int_{0}^{\tau_{N}\wedge t}\mathcal{H}^{2}(s)\mathbf{Z}_{2}(s) ds\\
                             &\leq\mathbb{E}\mathcal{E}(\text{u}_{0},\phi_{0})
                              +C\mathbb{E}\int_{0}^{\tau_{N}\wedge t} \mathcal{H}^{2}(s)\varUpsilon(s)ds
                             +C\mathbb{E}\int_{0}^{\tau_{N}\wedge t}\mathcal{H}^{2}(s)\mathbf{Z}_{1}(s)ds.
                             \end{aligned}
                               \end{equation}
                      By Gronwall's lemma and the fact that $\tau_{N}\uparrow T$ as $N$ goes to $\infty$, we derive that
                                              \begin{equation}\label{Q0m}
                                               \begin{aligned}
                                  \mathbb{E}\sup\limits_{s\in[0,t]}\left[ \mathcal{H}^{2}(s)\mathbf{Z}_{1}(s)\right] +2\mathbb{E}\int_{0}^{ t}\mathcal{H}^{2}(s)\mathbf{Z}_{2}(s) ds
                                   \leq C\mathbb{E}\mathcal{E}(\text{u}_{0},\phi_{0})
                                       +C\mathbb{E}\int_{0}^{t} \mathcal{H}^{2}(s)\varUpsilon(s)ds.
                                      \end{aligned}
                                       \end{equation}
            Therefore, taking into account that $\text{u}=\text{v}-\textbf{a}$ and  $\eqref{hj1}-\eqref{hj}$, we derive $\eqref{g50}$. This completes the proof of Theorem $\ref{12}$.                              
           \end{prev}   
           \section{SOLUTION TO THE CONTROL PROBLEM}\label{N03}   
          The  main goal of this article is to control the solution $\eqref{NS1}$ by boundary values $(a,b)$, which belongs to the space $\mathcal{A}$ of admissible controls  defined as a compact subset of $L^{2}(\Omega\times (0,T);\mathcal{H}_{p}(\Gamma))$. 
               In particular, there exists a constant $\delta>0$ such that:
                                  \begin{equation}\label{jk1}
                                  \mathbb{E}  \exp\left(4C_{0}\int_{0}^{T}\left\|(a,b) \right\|^{2}_{\mathcal{H}_{p}(\Gamma)} ds \right)< \delta, \ \ (a,b)\in \mathcal{A}.
                                  \end{equation}
              The cost functional is given by
                                  \begin{equation}\label{Ns2a}
                                              \begin{aligned}
                                              \mathcal{J}(a,b,\text{v},\phi)=\mathbb{E}\int_{0}^{T}\int_{\text{D}}\frac{1}{2}\left( \left\|\text{v}-\text{v}_{d} \right\|^{2}+\left\|\phi-\phi_{d} \right\|^{2}_{1}\right)dxdt
                                              +\mathbb{E}\int_{\Gamma_{T}}\frac{1}{2}\left(\lambda_{1}\left| a\right|^{2}+\lambda_{2}\left| b\right|^{2} \right) d\gamma dt.
                                              \end{aligned}
                                              \end{equation}
                   Our objective is to minimize $(\ref{Ns2a})$ over $\mathcal{A}$. More precisely, our goal is to solve the following problem:
                                              \begin{equation*}
                                              \textbf{(OCP)} \ \ \min_{(a,b)\in \mathcal{U}}\mathcal{J}(a,b,\text{v},\phi)
                                              \end{equation*}
              subject to $\eqref{NS1}$. Here $\text{v}_{d}$ and $\phi_{d}$ are the target functions such that $\text{v}_{d}, \phi_{d}\in L^{2}(\Omega\times \text{D}\times(0,T))$ and $\lambda_{1}, \lambda_{2}\geq 0.$\\
           The main result of this section is stated as follows:
                \begin{thm}\label{Theo3}
              Suppose that $(a,b)$ and $(\text{u}_{0},\phi_{0})$  satisfy the regularity $\eqref{l2}$ and $\eqref{n1}$ such that $(a,b)$ belongs to the space $\mathcal{A}$. Then, there exists at least one solution for the optimal control problem \textbf{(OCP)}.
            \end{thm}  
            \begin{prev}
            Let us consider a minimizing sequence
                          \begin{equation*}
                          ((a_{n},b_{n}),(\text{v}_{n},\phi_{n}))\in \mathcal{A}\times L^{2}(\Omega;L^{\infty}(0,T;\mathbb{Y}_{1})\cap L^{2}(0,T;\mathbb{V}_{1}))
                          \end{equation*}
                         of the cost functional $\mathcal{J}$, namely
                         \begin{equation*}
                         \mathcal{J}((a_{n},b_{n}),(\text{v}_{n},\phi_{n}))\rightarrow d=\inf(\textbf{OCP})\ \ \mbox{as}\ \ n\rightarrow\infty,
                         \end{equation*}
                         and $(\text{v}_{n},\phi_{n})$ is the solution $\eqref{NS1}$ for the sequence $(a_{n},b_{n})\in \mathcal{A}.$
                         \begin{equation}\label{g20}
                          \begin{cases}
                             d(\text{v}_{n},\varphi)=\left[-((\text{v}_{n},\varphi))+\int_{\Gamma}b_{n}(\varphi\cdot\tau)d\gamma-\left\langle (\text{v}_{n}\cdot\nabla)\text{v}_{n},\varphi\right\rangle  \right]dt \\
                             \left\langle \mu_{n}\nabla\phi_{n},\varphi\right\rangle dt+(g(t,\text{v}_{n}), \varphi)dW(s),\\
                             d(\phi_{n},\psi)=\left\langle-\mu_{n}-\text{v}_{n}\cdot\nabla\phi_{n} ,\psi\right\rangle dt, \ \
                             \mu_{n}=\text{A}_{\vartheta}\phi_{n}+f(\phi_{n}), \\ (\text{u}_{n},\phi_{n})(0)=(\text{u}_{0},\phi_{0})\in \mathbb{Y},
                               \end{cases}
                                \end{equation}
                     for  all $t\in[0,T]$, $\varphi\in\mathbb{V}_{div}$, $\psi\in\text{V}_{1}$ and $\mathbb{P}$-a.s.\\ 
             Since $\mathcal{A}$ is compact, then there exists a subsequence, still indexed by $n$ such that
                     \begin{equation}\label{k2}
                     (a_{n},b_{n})\rightarrow (a,b)\ \ \mbox{strongly in }\ \ L^{2}(\Omega\times(0,T); \mathcal{H}_{p}(\Gamma)).
                     \end{equation}
              It follows from  \cite[ Theorem 4.9, p. 94]{EVAN}, that there exists a subsequence of $(a_{n},b_{n})$, still denoted by $(a_{n},b_{n})$ a function $h\in L^{2}(\Omega\times(0,T))$ such that
                      \begin{equation}\label{k20}
                               \left\| (a_{n},b_{n})\right\|_{\mathcal{H}_{p}(\Gamma)} \leq h, \ \ \left\| (a,b)\right\|_{\mathcal{H}_{p}(\Gamma)}\leq h , \ \ \forall n\in \mathbb{N},\ \ \mbox{a.e. in }\ \  \Omega\times(0,T).
                               \end{equation}
             Let us take the function $h=h(t)$ and we define the following function
             \begin{equation}\label{k21}
                           \sigma_{h}(t)=\exp\left(-C_{0}t-\int_{0}^{t}h^{2}(s)ds \right) \ \ \mathbb{P}-\mbox{a.e. in} \, \Omega.
                           \end{equation}
                           Now, we replace $(a,b), \textbf{a}$ by  $(a_{n},b_{n}), \textbf{a}_{n}$, respectively in the relations  $\eqref{S1}$, $\eqref{hj}$, then  we derive from  $\eqref{g50a}$, $\eqref{n2b0}$, that the sequence $(\text{u}_{n},\phi_{n})=(\text{v}_{n}-\textbf{a}_{n},\phi_{n})$, $n\in \mathbb{N}$, satisfies the following estimates
                            \begin{equation}\label{g50aa}
                            \begin{aligned}
                     \mathbb{E}\sup\limits_{s\in[0,t]}\sigma_{h}^{2}(s)\left\|(\text{u}_{n},\phi_{n})(s) \right\|^{2}_{\mathbb{Y}}+\int_{0}^{t}\sigma_{h}^{2}(s)\left[ \left\|(\text{u}_{n},\phi_{n})(s) \right\|^{2}_{\mathbb{V}} \right] ds &\leq C, \\
                     \mathbb{E}\sup\limits_{s\in[0, t]}\left[\sigma_{h}^{2}(s)\left\| (\text{u}_{n},\phi_{n})(s)\right\|^{2}_{\mathbb{Y}} \right]^{2} +\mathbb{E}\left( \int_{0}^{ t}\sigma_{h}^{2}(s)\left[ \left\|(\text{u}_{n},\phi_{n})(s)\right\|^{2}_{\mathbb{V}}\right] ds\right) ^{2}&\leq C.
                            \end{aligned}
                             \end{equation}
                         It follows $\eqref{g50aa}$ that there exists a subsequence, still indexed by $n$, such that
                               \begin{equation}\label{p8a}
                                  \begin{aligned}
                                 \sigma_{h}(\text{u}_{n},\phi_{n})\rightharpoonup \sigma_{h}(\text{u},\phi) \ \ &\mbox{weakly in }\ \ L^{2}(\Omega\times(0,T); \mathbb{V})\cap L^{4}(\Omega, L^{2}(0,T; \mathbb{V})),\\
                                  \sigma_{h}(\text{u}_{n},\phi_{n})\rightharpoonup \sigma_{h}(\text{u},\phi) \ \ &\mbox{*-weakly in }\ \ L^{2}(\Omega,L^{\infty}(0,T; \mathbb{Y})\cap  L^{4}(\Omega,L^{\infty}(0,T; \mathbb{Y}).
                                  \end{aligned}
                                 \end{equation}
                                 Using \textbf{(H1)},  $\eqref{g50aa}$, we deduce that there exists some operators $\text{B}^{*}_{0}$, $\text{B}^{*}_{1}$, $\text{R}^{*}_{0}$ and $\text{g}^{*}$ such that
                                         \begin{equation}\label{p10a}
                                                     \begin{aligned}
                                                     \sigma_{h}^{2}(\text{y}_{n}\cdot\nabla)\text{y}_{n}\rightharpoonup \sigma_{h}^{2}\text{B}^{*}_{0}(t) \ \ &\mbox{weakly in }\ \ L^{2}(\Omega\times(0,T); \mathbb{V}^{*}_{div}),\\
                                                     \sigma_{h}^{2}(\mu_{n}\nabla\phi_{n})\rightharpoonup \sigma_{h}^{2}\text{R}^{*}_{0}(t) \ \ &\mbox{weakly in }\ \ L^{2}(\Omega\times(0,T); \mathbb{V}^{*}_{div}),\\
                                                     \sigma_{h}^{2}(\text{y}_{n}\cdot\nabla)\phi_{n}\rightharpoonup \sigma_{h}^{2}\text{B}^{*}_{1}(t) \ \ &\mbox{weakly in }\ \ L^{2}(\Omega\times(0,T); \text{V}_{1}^{*}),\\
                                                     \sigma_{h}^{2}\text{g}(t,\text{y}_{n})\rightharpoonup \sigma_{h}^{2}\text{g}^{*}(t) \ \ &\mbox{weakly in }\ \ L^{2}(\Omega\times(0,T); \mathbb{H}_{div}^{m}),
                                                     \end{aligned}
                                                     \end{equation}
                                                     and
                                                     \begin{equation}\label{p11a}
                                                     \mathcal{G}^{2}f_{\vartheta}(\phi_{n})\rightharpoonup \mathcal{G}^{2}f^{*}(t) \ \ \mbox{weakly in }\ \  L^{2}(\Omega,L^{\infty}(0,T; \text{H})).
                                                     \end{equation}
                    As in \textbf{Step 2 of the proof of Theorem $\ref{Theo2}$}, we pass to the limit equation $\eqref{g20}$ in the distributional sense, as $n\rightarrow\infty$, to obtain
                                     \begin{equation}\label{g20a}
                                                   \begin{cases}
                                                      d(\text{v},\varphi)=\left[-((\text{v},\varphi))+\int_{\Gamma}b(\varphi\cdot\tau)d\gamma-(\text{B}^{*}_{0}(t),\varphi) \right]dt \\
                                                      +(\text{R}^{*}_{0}(t),\varphi)dt+(\text{g}^{*}(t), \varphi)dW(s),\\
                                                      d(\phi_{n},\psi)=-\left( \mu, \psi\right)dt -\left( \text{B}^{*}_{1}(t), \psi\right)dt, \ \
                                                      \mu=\text{A}_{\vartheta}\phi+f^{*}(t), \\ (\text{u},\phi)(0)=(\text{u}_{0},\phi_{0})\in \mathbb{Y},
                                                        \end{cases}
                                                         \end{equation}
                                              for  all $t\in[0,T]$, $\varphi\in\mathbb{V}_{div}$, $\psi\in\text{V}_{1}$ and $\mathbb{P}$-a.e. in  $\Omega.$
            Let us write $\text{v}=\text{u}+\textbf{a}$,  $\text{v}_{n}=\text{u}_{n}+\textbf{a}_{n}$. It follows from $\eqref{g20}$ and $\eqref{g20a}$  with $\varphi=e_{i}$ that
             \begin{equation}\label{g2aa}
                \begin{cases}
                 d(\text{u}-\text{u}_{n},e_{i})=\left[-((\text{u}-\text{u}_{n},e_{i}))+((\text{v}_{n}\cdot\nabla)\text{v}_{n}-\text{B}^{*}_{0}(t),e_{i}) \right]dt \\
                  -((\textbf{a}-\textbf{a}_{n},e_{i}))-\left(\partial_{t}(\textbf{a}-\textbf{a}_{n}) ,e_{i} \right) +\int_{\Gamma}(b-b_{n})(e_{i}\cdot\tau)d\gamma\\
                   -(\mu_{n}\nabla\phi_{n}-\text{R}^{*}_{0}(t),e_{i})dt-(g(t,\text{v}_{n})-\text{g}^{*}(t), e_{i})dW(s),\\
                   d(\phi-\phi_{n},\psi)=\left((\mu-\mu_{n})+(\text{v}_{n}\cdot\nabla\phi_{n}-\text{B}^{*}_{1}(t)), \psi \right)dt,\\
                  (\mu-\mu_{n})=\text{A}_{\vartheta}(\phi-\phi_{n})+f^{*}(t)-f_{\vartheta}(\phi_{n}),
                  \end{cases}
                    \end{equation}
                                 Let  set
                                               \begin{math}
                                      \theta_{n}=\text{u}-\text{u}_{n}, \ \ \Psi_{n}=(\phi-\phi_{n}) \ \ \mbox{and}\ \ \xi_{n}=(\mu-\mu_{n}).\\
                                              \end{math}
              From the It\^{o} formula, we obtain
           \begin{equation}\label{k0a}
                                   \begin{aligned}
                            d\left( (\theta_{n},e_{i})^{2}\right) &=2(\theta_{n},e_{i})\left[-((\theta_{n},e_{i}))+((\text{v}_{n}\cdot\nabla)\text{v}_{n}-\text{B}^{*}_{0}(t),e_{i}) \right]dt \\
                             &-2(\theta_{n},e_{i})\left[ ((\textbf{a}-\textbf{a}_{n},e_{i}))-\left(\partial_{t}(\textbf{a}-\textbf{a}_{n}) ,e_{i} \right) +\int_{\Gamma}(b-b_{n})(e_{i}\cdot\tau)d\gamma\right]dt \\
                            &-2(\theta_{n},e_{i})\left[ (\mu_{n}\nabla\phi_{n}-\text{R}^{*}_{0}(t),e_{i})dt-(g(t,\text{v}_{n})-\text{g}^{*}(t), e_{i})dW(s)\right] \\
                             &+ \left|(g(t,\text{v}_{n})-\text{g}^{*}(t), e_{i}) \right|^{2}_{L^{2}}dt.
                             \end{aligned}
                              \end{equation}
                                 Summing over i=1,...,n, we obtain
                                \begin{equation}\label{k10}
                               \begin{aligned}
                               d\left| \theta_{n}\right|^{2}_{L^{2}} &+2\left\| \theta_{n}\right\|^{2}dt=2((\text{v}_{n}\cdot\nabla)(\text{v}_{n}-\text{B}^{*}_{0}(t)),\theta_{n}) dt \\
                               &-2\left[ ((\textbf{a}-\textbf{a}_{n},\theta_{n}))-\left(\partial_{t}(\textbf{a}-\textbf{a}_{n}) ,\theta_{n} \right) +\int_{\Gamma}(b-b_{n})(\theta_{n}\cdot\tau)d\gamma\right]dt \\
                              &-2(\mu_{n}\nabla\phi_{n}-\text{R}^{*}_{0}(t),\theta_{n})dt-2(g(t,\text{v}_{n})-\text{g}^{*}(t), \theta_{n})dW(t)\\
                              &+ \sum_{i=1}^{n}\left|(g(t,\text{v}_{n})-\text{g}^{*}(t), e_{i}) \right|^{2}_{L^{2}}dt.
                               \end{aligned}
                           \end{equation}
            Replacing $\psi$ in $\eqref{g2aa}_{2}$ by $\text{A}_{\vartheta}\Psi_{n}$ and multiplying the equation $\eqref{g2aa}_{3}$ with $\xi_{n}-\text{A}_{\vartheta}\Psi_{n}$ and adding these equations, we derive that  
            \begin{equation}\label{k3a}
                                 \begin{aligned}
                                  d\left|\nabla\Psi_{n}\right|^{2}_{L^{2}} &+2(\left\|\Psi_{n}\right\|^{2}_{2}+\left|\xi_{n} \right|^{2}_{L^{2}})dt = 2((\text{v}_{n}\cdot\nabla\phi_{n}-\text{B}^{*}_{1}(t)), \text{A}_{\vartheta}\Psi_{n})dt\\
                                &+\left(f_{\vartheta}(\phi_{n})-f^{*}(t), \text{A}_{\vartheta}\Psi_{n}\right)dt -\left(f_{\vartheta}(\phi_{n})-f^{*}(t),\xi_{n} \right) dt.
                               \end{aligned}
                               \end{equation}
                            Adding $\eqref{k3a}$ and $\eqref{k10}$, we obtain
                            \begin{equation}\label{gh0}
                            \begin{aligned}
                           &d\left( \left|\theta_{n}\right|^{2}_{L^{2}}+ \left\|\Psi_{n}\right\|^{2}_{1}\right) +2\left\| \theta_{n}\right\|^{2}dt+2(\left\|\Psi_{n}\right\|^{2}_{2}+\left|\xi_{n} \right|^{2}_{L^{2}})dt\\
                           &=2((\text{v}_{n}\cdot\nabla)\text{v}_{n}-\text{B}^{*}_{0}(t)),\theta_{n}) dt
                          + 2((\text{v}_{n}\cdot\nabla\phi_{n}-\text{B}^{*}_{1}(t)), \text{A}_{\vartheta}\Psi_{n})dt\\
                           &+\left(f_{\vartheta}(\phi_{n})-f^{*}(t), \text{A}_{\vartheta}\Psi_{n}\right)dt -\left(f_{\vartheta}(\phi_{n})-f^{*}(t),\xi_{n} \right) dt-2(\mu_{n}\nabla\phi_{n}-\text{R}^{*}_{0}(t),\theta_{n})dt \\
                           &-2\left[ ((\textbf{a}-\textbf{a}_{n},\theta_{n}))-\left(\partial_{t}(\textbf{a}-\textbf{a}_{n}) ,\theta_{n} \right) +\int_{\Gamma}(b-b_{n})(\theta_{n}\cdot\tau)d\gamma\right]dt\\
                           &-2(g(t,\text{v}_{n})-\text{g}^{*}(t), \theta_{n})dW(t)+ \sum_{i=1}^{n}\left|(g(t,\text{v}_{n})-\text{g}^{*}(t), e_{i}) \right|^{2}_{L^{2}}dt.
                           \end{aligned}
                            \end{equation}
                            Now, we see that  
              \begin{equation}\label{h1a}
                              \begin{aligned}
                             ((\text{v}_{n}\cdot\nabla)(\text{v}_{n}-\text{B}^{*}_{0}(t)),\theta_{n})
                             &=-(((\text{u}_{n}+\textbf{a}_{n})\cdot\nabla\theta_{n}),\theta_{n})-(\theta_{n}\cdot\nabla(\text{u}+\textbf{a}),\theta_{n})\\
                            &+((\text{u}_{n}+\textbf{a}_{n})\cdot\nabla(\textbf{a}_{n}-\textbf{a}),\theta_{n})+((\textbf{a}_{n}-\textbf{a})\nabla(\text{u}+\textbf{a}),\theta_{n})\\
                             &+((\text{v}\cdot\nabla)\text{v}-\text{B}^{*}_{0}(t),\theta_{n})=(B_{0}+B_{1},  \theta_{n} ) +B_{2}+B_{3}+B_{4},
                          \end{aligned}
                          \end{equation}
                            \begin{equation}\label{k5a0}
                            \begin{aligned}
                            (\mu_{n}\nabla\phi_{n}-\text{R}^{*}_{0}(t),\theta_{n})
                           &=-(\mu_{n}\nabla\Psi_{n},\theta_{n})-(\xi_{n}\nabla\phi,\theta_{n})
                            +(\mu\nabla\phi-\text{R}^{*}_{0}(t),\theta_{n})\\
                            &=B_{5}+B_{6}+B_{7},
                            \end{aligned}
                           \end{equation}
                           \begin{equation}\label{k5aa}
                          \begin{aligned}
                           &((\text{v}_{n}\cdot\nabla\phi_{n}-\text{B}^{*}_{1}(t)), \text{A}_{\vartheta}\Psi_{n})
                           =-(((\text{u}_{n}+\textbf{a}_{n})\cdot\nabla\Psi_{n}),\text{A}_{\vartheta}\Psi_{n})\\
                           &-(\theta_{n}\cdot\nabla\phi,\text{A}_{\vartheta}\Psi_{n})
                           +((\textbf{a}_{n}-\textbf{a})\cdot\nabla\phi,\text{A}_{\vartheta}\Psi_{n})\\
                           &+((\text{v}\cdot\nabla\phi-\text{B}^{*}_{1}(t)), \text{A}_{\vartheta}\Psi_{n})
                           =B_{8}+B_{9}+B_{10}+B_{11},
                           \end{aligned}
                           \end{equation}
                               and
                           \begin{equation}\label{h2a}
                          f_{\vartheta}(\phi_{n})-f^{*}(t)=-(f_{\vartheta}(\phi_{n})-f_{\vartheta}(\phi))+(f_{\vartheta}(\phi)-f^{*}(t)).
                          \end{equation}
                              Note that
                              \begin{equation}\label{e1a}
                                \begin{aligned}
                               \left| (B_{0}+B_{2},  \theta_{n} )\right|&\leq \left|\int_{\Gamma}a_{n} \left(\theta_{n}.\tau \right)^{2}d\gamma  \right|+\left|(\theta_{n}\cdot\nabla(\text{u}+\textbf{a}),\theta_{n}) \right|\\
                                &\leq C_{2}\left(\left\|(a,b) \right\|^{2} _{\mathcal{H}_{p}(\Gamma)}+\left\|(a_{n},b_{n}) \right\|^{2} _{\mathcal{H}_{p}(\Gamma)}+\left\|\text{u} \right\|^{2} \right)\left|\theta_{n} \right|^{2}_{L^{2}} +\frac{1}{4}\left\|\theta_{n} \right\|^{2}\\
                                 &  \leq C_{2}\left(2h^{2}+\left\|\text{u} \right\|^{2} \right)\left|\theta_{n} \right|^{2}_{L^{2}} +\frac{1}{4}\left\|\theta_{n} \right\|^{2},
                          \end{aligned}
                      \end{equation}  
                   \begin{equation}\label{e6a}
                                              \begin{aligned}
                                              \left| B_{5}\right| &\leq\left\|\theta_{n} \right\|_{\mathbb{L}^{4}} \left|\mu_{n} \right|_{L^{2}} \left\|\nabla\Psi_{n} \right\|_{L^{4}} \leq C\left\|\theta_{n} \right\|\left|\mu_{n} \right|_{L^{2}}\left\|\nabla\Psi_{n} \right\|_{L^{4}}\\
                                              & \leq \frac{1}{4}\left\|\theta_{n} \right\|^{2}+\frac{1}{4}\left\| \Psi_{n}\right\|^{2}_{2}+C\left|\mu_{n} \right|^{2}_{L^{2}}\left| \nabla\Psi_{n}\right|^{2}_{L^{2}} +C\left|\mu_{n} \right|^{4}_{L^{2}}\left| \nabla\Psi_{n}\right|^{2}_{L^{2}},
                                              \end{aligned}
                                              \end{equation}
                                             \begin{equation}\label{e7a}
                                            \begin{aligned}
                                            \left| B_{6}\right| &\leq\left\|\theta_{n} \right\|_{\mathbb{L}^{4}}\left|\xi_{n} \right|_{L^{2}} \left\|\nabla\phi \right\|_{L^{4}}\\
                                            &\leq \frac{1}{3}\left|\xi_{n} \right|^{2}_{L^{2}}+\frac{1}{4}\left\|\theta_{n} \right\|^{2}+C\left(\left\|\phi \right\|_{2}^{2} \left|\nabla\phi \right|^{2}_{L^{2}} +\left| \nabla\phi\right|^{4}_{L^{2}}  \right)\left|\theta_{n} \right|^{2}_{L^{2}},
                                             \end{aligned}
                                             \end{equation}
                                    \begin{equation}\label{e9a0}
                                              \begin{aligned}
                                             \left|B_{8}\right| &\leq \left\|\text{u}_{n}+\textbf{a}_{n} \right\|_{\mathbb{L}^{4}} \left| \nabla\Psi_{n}\right|^{1/2}_{L^{2}}\left\|\Psi_{n} \right\|^{1/2}_{2} \left\|\Psi_{n} \right\|_{2}
                                             \leq \frac{1}{3} \left\| \Psi_{n}\right\|^{2}_{2}+C\left\|\text{u}_{n}+\textbf{a}_{n} \right\|^{4}_{\mathbb{L}^{4}} \left| \nabla\Psi_{n}\right|^{2}_{L^{2}} \\
                                             &\leq \frac{1}{3} \left\| \Psi_{n}\right\|^{2}_{2}+C_{4}\left(\left\|(a_{n},b_{n}) \right\|^{4} _{\mathcal{H}_{p}(\Gamma)}+\left\|\text{u}_{n} \right\|^{4} \right)\left\|\nabla \Psi_{n}\right\|^{2}_{L^{2}}\\
                                             &\leq \frac{1}{3} \left\| \Psi_{n}\right\|^{2}_{2}+C_{4}\left(h^{4}+\left\|\text{u}_{n} \right\|^{4} \right)\left\|\nabla \Psi_{n}\right\|^{2}_{L^{2}},
                                             \end{aligned}
                                             \end{equation}  
        \begin{equation}\label{e10a}
                                  \begin{aligned}
                                  \left|B_{9}\right| &\leq \frac{1}{3}  \left\|\Psi_{n} \right\|_{2}^{2}+\frac{1}{4}\left\|\theta_{n} \right\|^{2}+C\left\|\nabla\phi \right\|^{2}_{L^{4}}\left| \theta_{n}\right|^{2}_{L^{2}},
                                  \end{aligned}
                                  \end{equation}
                   \begin{equation}\label{e14a}
                             \left| (f_{\vartheta}(\phi_{n})-f_{\vartheta}(\phi),\text{A}_{\vartheta}\Psi_{n} )\right| \leq \frac{1}{6}\left\|\Psi_{n} \right\|^{2}_{2}+C_{f}\left|\nabla\Psi_{n} \right|^{2}_{L^{2}},
                             \end{equation}
                           \begin{equation}\label{e15a}
                           \left| (f_{\vartheta}(\phi)-f_{\vartheta}(\phi_{n}), \xi_{n})\right| \leq \frac{1}{4}\left|\xi_{n} \right|^{2}_{L^{2}}+C_{f} \left|\nabla\Psi_{n} \right|^{2}_{L^{2}},
                          \end{equation}
                            Setting
                                            \begin{equation}\label{e4a}
                                            \text{g}_{n}=\text{g}(t,\text{v}_{n}),\ \ \text{g}=\text{g}(t,\text{v}), \ \ \text{g}^{*}=\text{g}^{*}(t).
                                            \end{equation}
                                            As in $\eqref{e5}$, we can show that
                         \begin{equation}\label{e5a}
                                          \begin{aligned}
                                         &\left|(g(t,\text{v}_{n})-\text{g}^{*}(t) \right|^{2}_{L^{2}}
                                         &\leq C_{3}\left|\theta_{n} \right|^{2}_{L^{2}}+C\left|\textbf{a}_{n}-\textbf{a} \right|^{2}_{L^{2}}- \left|\text{g}-\text{g}^{*} \right|^{2}_{L^{2}}
                                         +2(\text{g}_{n}-\text{g}^{*},\text{g}-\text{g}^{*}).
                                         \end{aligned}
                                      \end{equation}
                           The positive constant $C_{3}$ in $\eqref{e5a}$ is independent of $n$.\\
                               Let us define the functions:
                               \begin{equation*}
                               \begin{aligned}
                               \mathbf{Z}_{1}(t)&= \left|\theta_{n} \right|^{2}_{L^{2}}+\left\| \nabla\Psi_{n}\right\|^{2}_{L^{2}}, \ \ \mathbf{Z}_{2}(t)= \left\|\theta_{n} \right\|^{2}+\left\| \Psi_{n}\right\|^{2}_{2}+\left|\xi_{n} \right|^{2}_{L^{2}},\\
                               \widetilde{f}(s)&=C\left(1+\left|\mu_{n} \right|^{4}_{L^{2}}+\left|\mu_{n} \right|^{2}_{L^{2}}+\left\|\phi \right\|_{2}^{2} \left|\nabla\phi \right|^{2}_{L^{2}} +\left| \nabla\phi\right|^{4}_{L^{2}}+\left\|\nabla\phi \right\|^{2}_{L^{4}} \right) \\
                               &+\max\left(3C_{0},C_{2},C_{4} \right)\left(1+h^{2} +\left\|\text{u}_{n} \right\|^{2}+\left\|\text{u} \right\|^{2} \right),\ \
                               \sigma(t)=\exp(-\int_{0}^{t}\widetilde{f}(s)ds).
                               \end{aligned}
                               \end{equation*} 
            From $\eqref{e1a}-\eqref{e5a}$, using the It\^{o} formula and taking the expectation, we derive that
                               \begin{equation}\label{e18a}
                                           \begin{aligned}
                                            &\mathbb{E}\sigma(t)\mathbf{Z}_{1}(t)+2\mathbb{E}\int_{0}^{t}\sigma(s)\mathbf{Z}_{2}(s)ds+2\mathbb{E}\int_{0}^{t}\sigma(s)\left|\text{g}-\text{g}^{*} \right|^{2}_{L^{2}}ds\\
                                           &\leq 2\mathbb{E}\int_{0}^{t}\sigma(s)B_{1}ds+2\mathbb{E}\int_{0}^{t}\sigma(s)B_{3}ds+ 2\int_{0}^{t}\sigma(s)B_{10}ds+2\mathbb{E}\int_{0}^{t}\sigma(s)B_{4}ds\\
                                           &+2\mathbb{E}\int_{0}^{t}\sigma(s)B_{11}ds+2\mathbb{E}\int_{0}^{t}\sigma(s)B_{7}ds+2\mathbb{E}\int_{0}^{t}\int_{\Gamma}\sigma(s)(b-b_{n})(\theta_{n}\cdot\tau)d\gamma ds\\
                                          &+2C\mathbb{E}\int_{0}^{t}\sigma(s)\left\|(a_{n},b_{n})-(a,b) \right\|^{2}_{\mathcal{H}_{p}(\Gamma)}ds
                                          +2\mathbb{E}\int_{0}^{t}\sigma(s)(\text{g}_{n}-\text{g}^{*},\text{g}-\text{g}^{*})ds\\
                                         &+2C\mathbb{E}\int_{0}^{t}\sigma(s)\left|\nabla(\widetilde{\phi}_{n}- \phi)\right|^{2}_{L^{2}}ds
                                        +2C\mathbb{E}\int_{0}^{t}\sigma(s)(f_{\vartheta}(\phi)-f^{*}(t),\text{A}_{\vartheta}\Psi_{n} )ds\\
                                        &+2C\mathbb{E}\int_{0}^{t}\sigma(s)(f_{\vartheta}(\phi)-f^{*}(t),\xi_{n} )ds.
                                                      \end{aligned}
                                                      \end{equation}
             Now, we will prove that the right-hand side of this inequality tends to zero as $n\rightarrow\infty$. By $\eqref{l2}_{1}$, $\eqref{hj}$ and using the H\"{o}lder inequality and the fact that  $\sigma(t)\leq \sigma^{2}_{h}(t)$, we derive that    
          \begin{equation}\label{m1}
                      \begin{aligned}
                      I_{0}&=\mathbb{E} \int_{0}^{t}\int_{\Gamma}\sigma(s)(b-b_{n})(\theta_{n}\cdot\tau)d\gamma ds\\
                      & \leq C\left\|(a_{n},b_{n})-(a,b) \right\|_{L^{2}(\Omega\times (0,T);\mathcal{H}_{p}(\Gamma))}+\frac{1}{2}\mathbb{E}\int_{0}^{t}\sigma(s)\left\|\theta_{n}(s) \right\|^{2}ds
                      \end{aligned}
                      \end{equation}
                        Thanks to $\eqref{g50aa}$ and using the fact that $\sigma(t)\leq \sigma^{2}_{h}(t)$ on $(0,T)$, we derive as in \cite{NSC} that
                        \begin{equation}
                        \begin{aligned}
                        I_{1}&=\mathbb{E}\int_{0}^{T}\sigma(s)B_{1}ds\leq C\left(\mathbb{E}\int_{0}^{T}\left\|\textbf{a}_{n}-\textbf{a} \right\|^{2}_{\mathbb{H}^{1}}ds  \right)^{1/2} \\
                        &\times\left[\left(\mathbb{E}\int_{0}^{T}\sigma^{4}_{h}(s)\left\|\text{u}_{n} \right\|^{2}_{\mathbb{L}^{4}}\left\|\theta_{n} \right\|^{2}_{\mathbb{L}^{4}}ds\right)^{1/2} +\left(\mathbb{E}\int_{0}^{T}\sigma^{4}_{h}(s)\left\|\textbf{a}_{n} \right\|^{2}_{\mathcal{C}(\bar{D})} \left|\theta_{n} \right|^{2}_{L^{2}}  \right)^{1/2}  \right] \\
                        &\leq C\left\|(a_{n},b_{n})-(a,b) \right\|_{L^{2}(\Omega\times (0,T);\mathcal{H}_{p}(\Gamma))}.
                        \end{aligned}
                        \end{equation}
                        Similarly, we can prove that  
                                \begin{equation}
                               I_{2}=\mathbb{E}\int_{0}^{T}\sigma(s)B_{3}ds\leq C\left\|(a_{n},b_{n})-(a,b) \right\|_{L^{2}(\Omega\times (0,T);\mathcal{H}_{p}(\Gamma))}.
                                \end{equation}
                                Note also that 
                \begin{equation}\label{kl0}
                                \begin{aligned}
                                I_{3}&=\mathbb{E}\int_{0}^{T}\sigma(s)B_{10}ds\leq  \mathbb{E}\int_{0}^{T}\sigma(s)\left\|\textbf{a}_{n}-\textbf{a} \right\|_{\mathbb{H}^{1}}\left\|\nabla\phi \right\|_{L^{4}} \left\|\Psi_{n} \right\|_{2}ds\\
                                &\leq C\left( \mathbb{E}\int_{0}^{T}\left\|\textbf{a}_{n}-\textbf{a} \right\|^{2}_{\mathbb{H}^{1}}ds  \right)^{1/2}
                                \times \left(\mathbb{E}\int_{0}^{T}\sigma^{4}_{h}(s)\left\|\nabla\phi \right\|^{2}_{L^{4}} \left\|\Psi_{n} \right\|^{2}_{2}ds \right)^{1/2}\\
                                &\leq C\left\|(a_{n},b_{n})-(a,b) \right\|_{L^{2}(\Omega\times (0,T);\mathcal{H}_{p}(\Gamma))}.
                                \end{aligned}
                                \end{equation}                                            
                 The term $C\mathbb{E}\int_{0}^{t}\sigma(s)\left\|(a_{n},b_{n})-(a,b) \right\|^{2}_{\mathcal{H}_{p}(\Gamma)}ds$,  we give
                                                \begin{equation}\label{mp}
                                                I_{4}=C\mathbb{E}\int_{0}^{t}\sigma(s)\left\|(a_{n},b_{n})-(a,b) \right\|^{2}_{\mathcal{H}_{p}(\Gamma)}ds\leq C\left\|(a_{n},b_{n})-(a,b) \right\|_{L^{2}(\Omega\times (0,T);\mathcal{H}_{p}(\Gamma))}.
                                                \end{equation}
                 It follows from $\eqref{m1}-\eqref{mp}$, that the terms $I_{i}, i=1,2,3,4$ converge to zero as $n\rightarrow \infty$.\\   
          Using $\eqref{p8a}$, we derive that
                         \begin{equation}\label{qsda}
                         \begin{aligned}
                          \sigma_{h}^{2}((\text{u},\phi)-(\text{u}_{n},\phi_{n}))&\rightharpoonup 0\ \ \mbox{weakly in }\ \ L^{2}(\Omega\times(0,T), \mathbb{V}), \ \ n\rightarrow\infty,\\
                           \sigma_{h}^{2}(\mu- \mu_{n})&\rightharpoonup 0\ \ \mbox{weakly in }\ \ L^{2}(\Omega\times(0,T), L^{2}(D)), \ \ n\rightarrow\infty.
                           \end{aligned}
                          \end{equation}
                         The operators  $\sigma_{h}^{2}((\text{y}\cdot\nabla)\text{y}-\text{B}^{*}_{0})$  and $\sigma_{h}^{2}(\mu\nabla\phi-\text{R}^{*}_{0})$  belong to  $L^{2}(\Omega\times(0,T), \mathbb{V}^{*}_{div})$. Thanks to $\eqref{p6}$ and   $\eqref{e4}$, it follows that
                                                \begin{equation}
                                                \mathbb{E}\int_{0}^{t}\sigma(s)B_{4}ds+\mathbb{E}\int_{0}^{t}\sigma(s)B_{7}ds=0\ \ \mbox{as}\ \ n\rightarrow\infty.
                                                \end{equation}
                    Similarly, we deduce that the operators  $\sigma_{h}^{2}((\text{y}\cdot\nabla)\phi-\text{B}^{*}_{0})$ belongs to  $L^{2}(\Omega\times(0,T), \text{V}^{*}_{1})$. Thanks to $\eqref{p6}$ and   $\eqref{p10}$, it follows that
                             \begin{equation}
                            \mathbb{E}\int_{0}^{T}\sigma(s)B_{11}ds=0\ \ \mbox{as}\ \ n\rightarrow\infty.
                               \end{equation}
          From $\eqref{qsda}$, we deduce that
                             \begin{equation}\label{jk}
                             C\mathbb{E}\int_{0}^{t}\sigma(s)(f_{\vartheta}(\phi)-f^{*}(t),\text{A}_{\vartheta}\Psi_{n} )ds+C\mathbb{E}\int_{0}^{t}\sigma(s)(f_{\vartheta}(\phi)-f^{*}(t),\xi_{n} )ds=0 \ \mbox{as}\ \ n\rightarrow\infty.
                             \end{equation}
               Thanks to the convergence result $\eqref{p10a}$ and $\eqref{e4a}$, we obtain
                             \begin{equation}\label{fta}
                                 \begin{aligned}
                               \sigma_{h}(\text{g}_{n}-\text{g}^{*})\rightharpoonup 0\ \ \mbox{weakly in }\ \ L^{2}(\Omega\times(0,T), \mathbb{H}^{m}_{div}),
                                 \end{aligned}
                             \end{equation}
            which implies
                          \begin{equation}\label{fg}
                        2\mathbb{E}\int_{0}^{t}\sigma(s)(\text{g}_{n}-\text{g}^{*},\text{g}-\text{g}^{*})ds\rightarrow 0\ \ \mbox{as}\ \ n\rightarrow\infty.
                         \end{equation}
             From $\eqref{m1}-\eqref{fg}$ and passing to the limit in the inequality $\eqref{e18a}$, we obtain the following strong convergences
                         \begin{equation}\label{fca}
                          \lim\limits_{n\rightarrow\infty}\mathbb{E}\sigma(t)\mathbf{Z}_{1}(t)=\lim\limits_{n\rightarrow\infty}\mathbb{E}\int_{0}^{t}\sigma(s)\mathbf{Z}_{2}(s)ds=0,
                           \end{equation}
              for $t\in(0,T)$. Hence, we have
                          \begin{equation}\label{hga}
                           \begin{aligned}
                          \lim\limits_{n\rightarrow\infty}\mathbb{E}\left( \sigma(t)\left\|(\text{u}_{n},\phi_{n})(t)-(\text{u},\phi)(t) \right\|^{2}_{\mathbb{Y}}\right) &=0,\\
                           \lim\limits_{n\rightarrow\infty}\mathbb{E}\int_{0}^{t}\sigma(s)\left\|(\text{u}_{n},\phi_{n})(s)-(\text{u},\phi)(s) \right\|^{2}_{\mathbb{V}}ds &=0.
                           \end{aligned}
                          \end{equation}
            In addition, we derive that
             \begin{equation*}
              \mathbb{E}\int_{0}^{t}\sigma(s)\left|\text{g}(s,\text{v})-\text{g}^{*}(s) \right|^{2}_{L^{2}}ds=0.
                  \end{equation*}
                This implies that
                       \begin{equation}\label{ghta}
                        \text{g}(t,\text{v})=\text{g}^{*}(t), \ \ \mbox{a.e. in}\ \ (\omega,t)\in\Omega\times(0,T).
                            \end{equation}
            By $\eqref{p10a}$ and by $\eqref{hga}$,  we infer that $\sigma(t)(\text{v}\cdot\nabla)\text{v}=\sigma(t)\text{B}^{*}_{0}(t)$, \ \ $\mbox{a.e. in}\ \ (\omega,t)\in\Omega\times(0,T)$, that implies
                                                     \begin{equation}\label{ght1a}
                                                        (\text{v}\cdot\nabla)\text{v}=\text{B}^{*}_{0}(t), \ \ \mbox{a.e. in}\ \ (\omega,t)\in\Omega\times(0,T).
                                                         \end{equation}
                  Similarly, we can deduce that the operators $\sigma(t)(\mu\nabla\phi)=\sigma(t)\text{R}^{*}_{0}(t)$ and $\sigma(t)(\text{v}\cdot\nabla)\phi=\sigma(t)\text{B}^{*}_{1}(t)$, \ \ $\mbox{a.e. in}\ \ (\omega,t)\in\Omega\times(0,T)$, we infer that
                                                      \begin{equation}\label{ght2a}
                                                      \begin{aligned}
                                                       \mu\nabla\phi&=\text{R}^{*}_{0}(t),\ \ \mbox{a.e. in}\ \ (\omega,t)\in\Omega\times(0,T),\\
                                                       (\text{v}\cdot\nabla)\phi&=\text{B}^{*}_{1}(t), \ \ \mbox{a.e. in }\ \ (\omega,t)\in\Omega\times(0,T).
                                                       \end{aligned}
                                                     \end{equation}
              From $\eqref{ghta}, \eqref{ght1a}, \eqref{ght2a}$, we obtain the following system:
                                               \begin{equation}\label{Ton2aa}
                                               \begin{cases}
                                              (\text{v}(t),\text{w})=(\text{v}_{0},\text{w}) +\int_{0}^{t}\left[  -((\text{v}(s),\text{w}))+\int_{\Gamma}b(\text{w}\cdot\tau)d\gamma -\left\langle  (\text{v}(s)\cdot\nabla\text{v}(s)),\text{w}\right\rangle  \right]   ds\\
                                              +\int_{0}^{t}\left\langle (\mu\nabla\phi(s), \text{w}\right\rangle ds
                                              +\int^{t}_{0}(\text{g}(s,\text{v}(s)),\text{w})dW_{s},\\
                                              (\phi(t),\psi)=(\phi_{0},\psi) -\int_{0}^{t}\left( \mu(s),\psi\right)  ds- \int_{0}^{t}\left\langle (\text{v}(s)\cdot\nabla\phi(s)),\psi\right\rangle   ds,\\
                                             \mu=\text{A}_{\vartheta}\phi + f_{\vartheta}(\phi),
                                             \end{cases}
                                              \end{equation}
        for  all $t\in[0,T]$, $\text{w}\in\mathbb{V}_{div}$, $\psi\in\text{V}_{1}$ and $\mathbb{P}$-a.e. in $\Omega$. We deduce that $(\text{v},\phi)$ is the solution $\eqref{NS1}$, corresponding to the control pair $(a,b)$.\\
                 Now, using the lower semicontinuity of the cost functional, the strong convergence $\eqref{hga}$, $\eqref{jk1}$ and Fatou's Lemma, we deduce that
                                    \begin{equation*}
                                    \mathcal{J}((a,b),\text{v},\phi)\leq \lim\limits_{n\rightarrow\infty}\mathcal{J}((a_{n},b_{n}),\text{v}_{n},\phi_{n}),
                                    \end{equation*}
                        which gives
                                   \begin{equation*}
                                    \mathcal{J}((a,b),\text{v},\phi)=\inf(\textbf{OCP}).
                                                \end{equation*}
           Hence, the triplet $\left( (a,b),\text{v},\phi\right) $  is a solution to the control problem $(\textbf{OCP})$.
                      \end{prev}       
                   \section*{Acknowledgments}
                      The second author acknowledges the financial support from IMU (International Mathematical Union) and the GRAID Program (Graduate Research Assistantship in Developing Countries).
                 \begin{center}
                   
          \end{center}                                                            

\begin{thebibliography}{99}\addcontentsline{toc}{chapter}{References}
      \bibitem[1] {A30}  D. M Aderson, G. B Mefadden,  A.A Whecler, Diffuse-interface mehods in fluid mechanics, Annu. Rev. Fluid mech, \textbf{30}, 139-165.   
    \bibitem[2]{josue} R. D. Ayissi, G. Deugoué, J. Ngandjou Zangue and T. Tachim Medjo: Existence of optimal and $\varepsilon$-optimal controls for the stochastic Cahn-Hilliard Navier-Stokes system (Submitted).
    \bibitem[3] {PB1} P. Benner and C. Trautwein. Optimal control problems constrained by the stochastic Navier–Stokes equations with multiplicative L\'{e}vy noise. Mathematische Nachrichten, \textbf{292}(7),1444-1461, (2019).
     \bibitem[4] {PB2} P. Benner and C. Trautwein. A Stochastic Maximum Principle for Control Problems Constrained by the Stochastic Navier–Stokes Equations. Applied Mathematics and Optimization,\textbf{ 84}, 1001-1054, (2021).
      \bibitem[5] {PB3} P. Benner and C. Trautwein. Optimal Distributed and Tangential Boundary Control for the Unsteady Stochastic Stokes Equations. arXiv:1809.00911, (2018).
      \bibitem[6]{A.B} A. Bensoussan, Stochastic maximum principle for distributed parameter systems, J. Franklin Inst. \textbf{315}, 387–406 (1983).
       \bibitem[7] {TDM1} T. Biswas, S. Dharmatti and M. T. Mohan. Pontryagin maximum principle and second order optimality conditions for optimal control problems governed by 2D nonlocal Cahn-Hilliard Navier-Stokes equations. Analysis, \textbf{40} (3), 127-150, (2020).
      \bibitem[8] {TDM2} T. Biswas, S. Dharmatti and M. T. Mohan, Maximum principle and data assimilation problem for the optimal control problems governed by 2D nonlocal Cahn-Hilliard-Navier-Stokes equations. J. Math. Fluid Mech \textbf{22}, 1-42 (2020).
      \bibitem[9]{JM} J. M. Bismut, An introductory approach to duality in optimal stochastic control, SIAM Rev.,\textbf{20} , pp. 62–78 (1978).
      \bibitem[10]{ALB} A. L. Braslow, A history of suction-type laminar-flow control with emphasis on flight Research, NASA History Division, 1999.
     \bibitem[11]{TLB}T. L. Black and A. J. Sarnecki, The turbulent boundary layer with suction or injection Aeronautical Research Council Reports and Memoranda 3387 (October 1958), London, 1965.
     \bibitem[12]{S.H}  N. Chemetov and F. Cipriano, Injection-Suction control for two- dimensional Navier-Stokes equations with slippage, SIAM J.Control Optim, Vol. \textbf{56}, No. 2, pp 1253-1281, (2018).
     \bibitem[13]{KIS}  N. Chemetov and F. Cipriano, Well-posedness of stochastic second grade fluids, J. Math. Anal. Appl. \textbf{454}, 585-616, (2017).
      \bibitem[14]{NSC} N. Chemetov and F. Cipriano: A boundary control problem for stochastic 2D-Navier-Stokes equations. arXiv:2312.05935. (2023).
      \bibitem[15] {AD3}F. Coron, Deviation of slip boundary conditions for the Navier-Stokes system from Boltzmann equation, Journal of Statistical physics, vol \textbf{54} No. 314, (19840).
     \bibitem[16]{GTLD} G. Deugoué, and T. Tachim Medjo: Large deviation for a 2D Cahn-Hilliard Navier-Stokes model under random influences. J.Math.Anal.Appl. \textbf{486}(1):123863,34 pp (2020).
     \bibitem[17] {DN2} G. Deugoué and T. Tachim Medjo, Convergence of the solutions of the stochastic 3D globally modified Cahn-Hilliard-Navier-Stokes equations. J. Differential Equations.  \textbf{265},  545–592, (2018).
      \bibitem[18]{GTC} G. Deugoué, T. Tachim Medjo, Convergence of the solution of the  stochastic globally modified Cahn-Hilliard-Navier-Stokes equation, J. Differential. Equation. \textbf{460} (1), 140–163 (2018).
      \bibitem[19]{Ta33} G. Deugoué and T. Tachim Medjo, On a Stochastic 2D Cahn-Hilliard-Navier-Stokes System Driven by Jump Noise, Communications on Stochastic Analysis. Vol. \textbf{13} : No. 1 , Article 5,(2019) .
     \bibitem[20]{TD5} G. Deugoué and T. Tachim Medjo, The exponential behavior of a stochastic globally modified Cahn-Hilliard-Navier-Stokes. Journal of Mathematical Analysic and Applications. \textbf{460}(1), 140-163 (2018).
      \bibitem[21]{jid} G. Deugoué, B. Jidjou Moghomye and T. Tachim Medjo, Splitting-up scheme for the stochastic Cahn–Hilliard Navier–Stokes model. Stochastics and Dynamics.1-46 (2021).
    \bibitem[22] {A251} G. Deugoué and T. Tachim Medjo. Large Deviation for a 2D Allen–Cahn–Navier–Stokes Model Under Random Influences. In: Asymptot. Anal. 123.1-2 (2021), pp. 41–78.
    \bibitem[23] {EVAN}  L. C. Evans, Partial Differential Equations, AMS, Graduate Studies in Mathematics, Vol.\textbf{19}, 1998.
     \bibitem[24] {A25} Gurtin, M. E Polignone D., Vi\~{n}als, Two-phase binary fluids and immiscible fluids described by an order parameter. Math. Models Methods Appl. Sci \textbf{6}, no. 6, 815-831 (1996).
    \bibitem[25]{YB0} M Hintermuller and D. Wegner Optimal control of a semidiscrete Cahn-Hilliard-Navier-Stokes system.  SIAM J. Control Optim, \textbf{52}(1): 747-772, (2014).
     \bibitem[26] {AD1} D. Iftime and F. Sueur, Viscous boundary layers for the Navier-Stokes equations with the Navier slip conditions. Arch. Rational. Mech. Anal, 145-175 (2011).
     \bibitem[27] {YB}You B Li, F. Optimal distributed control for a model of homogeneous incomprissible two-phase flows. J. Dyn Control Syst, (2020).
      \bibitem[28] {H.B1}  H. Lisei, Approximation and optimal control of the stochastic Navier-Stokes equation. Dissertation, Martin-Luther University, Halle-Wittenberg, Germany, (1999).
      \bibitem[29] {H.B2} H. Lisei, Galerkin approximation and the strong solution of the stochastic Navier–Stokes equation, J. Appl. Math Stoch. Anal. \textbf{13}, 239–259 (2000).
       \bibitem[30] {H.B3} H. Lisei, Existence of optimal and $\epsilon$-optimal controls for the stochastic Navier-Stokes equations, Nonlinear Analysis. \textbf{51}, 95-119 (2002).
      \bibitem[31] {AD2}J.P. Kelliher, Navier-Stokes equations with Navier-boundary condition for a bounded domain in the plane. SIAM. J. MATH. ANAL, vol \textbf{38} No 1, 210-232.
       \bibitem[32] {TW} T. T. Medjo: Wong-Zakai approximation for a stochastic 2D Cahn-Hilliard Navier-Stokes model. Reseach Square:10.21203/rs.3.rs-2600198/v1 (2023).
       \bibitem[33] {TAL} T. T. Medjo: Asymptotic Log-harnack inequality for the 2D stochastic Cahn-Hilliard Navier-Stokes  system with degenerate noise. Reseach Square: 10.21203/rs.3.rs-2600062/v1.(2023).
       \bibitem[34]{Ta1} T. T. Medjo: On the existence and uniqueness of solution to a stochastic 2D Cahn-Hilliard-Navier-Stokes model. J. Differential Equations, \textbf{262}, 1028-1054 (2017).
       \bibitem[35] {A20} T. Tachim Medjo. On the existence and uniqueness of solution to a stochastic 2D
                 Allen-Cahn-Navier-Stokes model. In: Stoch. Dyn. 19 (2019), p. 1950007.
     \bibitem[36] {A21} T. Tachim Medjo. On weak martingale solutions to a stochastic Allen-Cahn-Navier-Stokes model with inertial effects. In: Discrete Contin. Dyn. Syst. Ser. B (2021).
        \bibitem[37] {Moh} M. T. Mohan, K. Sakthivel and S. S Sritharan, Dynamic programming of the stochastic 2D Navier-Stokes equations forced by L\'{e}vy noise. Analysis of PDEs. \textbf{57}(5) 3571-3602 (2019).
        \bibitem[38] {NA}  Navier, C.-L.-M.-H.: Mémoire sur les lois du mouvement des fluides. Mem. Acad. R. Sci. Paris 6, 389–416 (1823).
      \bibitem[39] {Peng1990} S. Peng, A general stochastic maximun principle for optimal problems, SIAM J. Control And Optimization, Vol. \textbf{28}, No. 4, pp. 966-979, July (1990) .
      \bibitem[40] {Zhang} Z. Qiu and H. wang: Large deviation principle for a 2D stochastic Cahn-Hilliard Navier Stokes equations. Z .Angew. Math.Phys. \textbf{71}, 88(2020).
    \bibitem[41] {MA10} L. Scarpa, Optimal distributed control of a stochastic Cahn-Hilliard  equation. SIAM J. Control Optim. \textbf{57}, no. \textbf{5}, 3571-3601, (2019).
    \bibitem[42] {A250} L. Scarpa. Analysis and optimal velocity control of a stochastic convective Cahn-Hilliard equation. In: J. Nonlinear Sci. 31.2, No. 45-57. (2021).
     \bibitem[43] {TM} R. Temam. Navier-Stokes equations. Theory and numerical analysis. North-Holland Publishing Co., Amsterdam/New York/Oxford (1977).
   \bibitem[44] {12B} B. Oksendal, Stochastic differential equations. Springer Verlag, Berlin-New York (1985).
    \bibitem[45] {EZ1} E. Zeidler, Nonlinear Functional Analysis and its Applications, Vol. \textbf{III}: Variational Methods and Optimization, Springer Verlag, New York-Berlin (1985).
              \end{thebibliography}
                      \end{document}